\documentclass[10pt,reqno]{amsart}

\usepackage{amsmath,amsfonts,amssymb,amsthm}
\usepackage[latin1]{inputenc}
\usepackage{verbatim}
\usepackage{hyperref}

\newtheorem{theorem}{Theorem}[section]
\newtheorem{lemma}[theorem]{Lemma}
\newtheorem{prop}[theorem]{Proposition}
\newtheorem{cor}[theorem]{Corollary}
\newtheorem{rmk}[theorem]{Remark}
\newtheorem{defn}[theorem]{Definition}

\newcommand{\sect}{\section}
\numberwithin{equation}{section}

\newcommand{\supp}{\mbox{\, \rm supp \, } }

\newcommand{\w}[1]{\langle {#1} \rangle}

\newcommand{\pf}{\noindent {\bf Proof. \hspace{2mm}}}
\newcommand{\ef}{ \hfill $ \Box $ \vskip 3mm}
\newcommand{\be}{\begin{equation}}
\newcommand{\ee}{\end{equation}}
\newcommand{\bea}{\begin{eqnarray}}
\newcommand{\eea}{\end{eqnarray}}
\newcommand{\beas}{\begin{eqnarray*}}
\newcommand{\eeas}{\end{eqnarray*}}
\newcommand{\ep}{{\epsilon}}
\newcommand{\f}{\frac}
\newcommand{\na}{\nabla}

\newcommand{\nn}{\nonumber}

\newcommand{\bC}{{\mathbb C}}
\newcommand{\bR}{{\mathbb R}}
\newcommand{\bN}{{\mathbb N}}

\newcommand{\bD}{{\bf D}}

\newcommand{\bJ}{{\bf J}}

\newcommand{\bL}{{\bf L}}

\newcommand{\bP}{{\bf P}}

\newcommand{\bT}{{\bf T}}

\newcommand{\bW}{{\bf W}}

\newcommand{\bX}{{\bf X}}
\newcommand{\bF}{{\bf F}}

\newcommand{\vA}{{\mathcal A}}
\newcommand{\vB}{{\mathcal B}}
\newcommand{\vC}{{\mathcal C}}
\newcommand{\vD}{{\mathcal D}}
\newcommand{\vE}{{\mathcal E}}

\newcommand{\vG}{{\mathcal G}}
\newcommand{\vI}{{\mathcal I}}

\newcommand{\vH}{{\mathcal H}}
\newcommand{\vK}{{\mathcal K}}
\newcommand{\vL}{{\mathcal L}}
\newcommand{\vM}{{\mathcal M}}
\newcommand{\vN}{{\mathcal N}}
\newcommand{\vO}{{\mathcal O}}
\newcommand{\vP}{{\mathcal P}}
\newcommand{\vQ}{{\mathcal Q}}
\newcommand{\vR}{{\mathcal R}}
\newcommand{\vT}{{\mathcal T}}
\newcommand{\vU}{{\mathcal U}}
\newcommand{\vV}{{\mathcal V}}
\newcommand{\vX}{{\mathcal X}}
\newcommand{\vW}{{\mathcal W}}
\newcommand{\vY}{{\mathcal Y}}
\newcommand{\vZ}{{\mathcal Z}}
\newcommand{\vS}{{\mathcal S}}

\newcommand{\gU}{{\mathfrak U}}
\newcommand{\gm}{{\mathfrak m}}

\newcommand{\Ran}{\mbox{\rm Ran }}

\newcommand{\p}{\partial}
\newcommand{\dG}{{\rm d}\Gamma}

\def\p{\partial}

\def\f{\frac}

\def\na{\nabla}

\usepackage{xcolor}
\usepackage{graphicx}

\hypersetup{hidelinks}

\begin{document}

\title[Resonances of the KFP operator]{
Eigenvalues and resonances \\ of the Kramers--Fokker--Planck operator \\ in the large-friction limit}

\author{Xue Ping WANG}
\address[X.P. Wang]{School of Mathematics and Key Laboratory of MIIT, Nanjing University of Aeronautics and Astronautics, Nanjing 211106, China, and Laboratoire de Math\'ematiques Jean Leray, Nantes Universit\'e, F-44322 Nantes Cedex 3, France}
\email{xue-ping.wang@univ-nantes.fr}

\date{\today}


\begin{abstract}
We study the eigenvalues and resonances of the Kramers--Fokker--Planck operator by the methods of analytic dilation and phase-space deformation. We prove that in the large-friction limit $\gamma \to +\infty$ there are neither eigenvalues nor resonances away from the thresholds $\gamma \bN$, and that near each threshold $\gamma k$, there exists a bijection between the eigenvalues and resonances of the Kramers--Fokker--Planck operator and those of the Witten Laplacian at the level $k$, $-\Delta^{(k)}_V=-\Delta_V\otimes 1_{\mu_k}-\dG_k(\mbox{\rm Hess}\,V)$.
We also show that each negative eigenvalue of $-\Delta^{(k)}_V$, $k\ge 1$, produces a pole of the resolvent at distance of order $\gamma^{-1}$ below $\gamma k$, whose imaginary part is $O(e^{-c\gamma})$ for every $c<2\delta_V$, $\delta_V$ being the width of a strip of analyticity of $V$.
\end{abstract}

\maketitle
\setcounter{tocdepth}{1}
\tableofcontents

\sect{Introduction}\label{sect1}

The Kramers--Fokker--Planck (KFP, for short) equation, also called the Kramers equation or the kinetic Fokker--Planck equation (\cite{risk,ne,ls}),
is the evolution equation for the distribution function describing the Brownian motion of particles in an external field:
\be \label{kfp1}
\frac{\partial W}{\partial t} = \left( -v\cdot \nabla_x  + \nabla_v \cdot(\gamma v  -F(x))
+ \f{\gamma \kappa T }{m} \Delta_v \right)W,
\ee
where $W=W(t, x,v)$ for $x, v \in \bR^n$ and $t > 0$, and $F(x) = -  \nabla V(x)$ is the  external force.
Here $x$ and $v$ are the position and the velocity of the particles, $m$ is the mass, $\kappa$ the Boltzmann constant, $\gamma$ the friction coefficient and $T$ the temperature. This equation is a special case of the general Fokker--Planck equation, or Kolmogorov forward equation, for continuous-time diffusion processes.
Setting
\[
W(t, x, v) =e^{-\f{m}{2\kappa T} (\f{|v|^2}{2} + V(x))} u(t, x,v),
\]
one sees that \eqref{kfp1} is equivalent to
\be \label{kfp2}
\frac{\partial u}{\partial t} + \left( v\cdot \nabla_x  - \nabla V(x)\cdot \nabla_v  + \gamma (
- \f{\kappa T }{m} \Delta_v + \f{m }{4 \kappa T }|v|^2 -\f n 2) \right)u =0.
\ee
We are interested in the resonances of the associated KFP operator in the overdamped regime $\gamma \to +\infty$.
Omitting the physical constants for the moment, we set
\be\label{operator}
P =  -\Delta_v+\f{1}{4}|v|^2-\f{n}{2} + v\cdot\na_x-\na V(x)\cdot\na_v, \quad (x,v) \in \bR^{2n}.
\ee
Corresponding to the two different viewpoints, scattering and kinetic, we shall use the following two decompositions of $P$:
\be
P = P_0 + W= H_0 + L,
\ee
where
\beas
P_0 &=& -\Delta_v+\f{1}{4}|v|^2-\f{n}{2} + v\cdot\na_x, \qquad
W   =   -\na V(x)\cdot\na_v, \\
H_0  & =&   -\Delta_v+\f{1}{4}|v|^2-\f{n}{2},\qquad
L  =  v\cdot\na_x-\na V(x)\cdot\na_v.
\eeas
By an abuse of notation, we still denote by $P$ (resp. $P_0$) its closure, with maximal domain, in $L^2(\bR^{2n}_{x,v})$.

In this work, the potential $V$ is a real-valued function of class at least $C^1$. To define resonances by analytic dilation, we assume in addition that $V(x)$ is  dilation-analytic: there exists $a\in ]0, \f \pi 4]$ such that $V_\theta(x) = V(e^{\theta}x)$, defined for $\theta$ real, extends holomorphically in $\theta$ to the disc $\{\theta \in \bC;\ |\theta|<a\}$, and that, for some positive constants $\rho$ and $C$,
\be \label{ass1}
|(\nabla_xV)(e^{\theta}x)| \le C \w{x}^{-\rho}, \qquad x\in \bR^n, |\theta| <a.
\ee
Under condition \eqref{ass1}, the resonances of the KFP operator $P$ are defined in \cite{w4} by analytic dilation in the $x$-variables. A definition by analytic distortion is also given in \cite{w4}, but it is only valid in dimension one, because global subellipticity fails for the distorted operator when $n \ge 2$. They can also be defined through appropriate phase-space deformation (see Section \ref{ssec-phase}). All these definitions are equivalent when they are applicable (\cite{hm,w4}).

For $\theta \in \bR$, let $U_\theta$ denote the unitary operator
\[
U_\theta f(x,v) = e^{\f{n\theta} 2} f(e^{\theta}x, v), \qquad f \in L^2(\bR^{2n}).
\]
Let $P(\theta) = U_\theta P U_\theta^{-1}$, and define $P_0(\theta)$, $W_\theta$ and $L_\theta$ in the same way. Then
\be
P(\theta) = P_0(\theta) + W_\theta  = H_0 + L_\theta.
\ee
Under assumption \eqref{ass1}, $P(\theta)$, defined for $\theta\in \bR$, extends to $\theta \in \bC$ with $|\Im\theta|<a$. Throughout the paper, we write $\varphi = \Im\theta$. The key ingredient in the definition of resonances by analytic dilation is the following global subelliptic estimate: for any $\phi_0\in \,]0, \f \pi 2[$, there exist positive constants $C$ and $s_0$ such that
\be \label{sub}
\| (H_0 +s)f\| +  \|\w{D_x}^{\f 2 3} f\|  + \| v\cdot \nabla_x f\| \le C \|(P_0( \theta) +s) f\|,
\ee
for $\theta \in \bC$ with $|\theta|\le \phi_0$, $s \ge s_0$ and $f \in \vS(\bR_{x,v}^{2n})$ (see \cite[Theorem 2.1]{w4}). It is also proved in \cite{w4} that for $n=1$ estimate \eqref{sub} still holds for the analytically distorted free KFP operator $\widetilde P_0(\theta)$, whereas for $n \ge 2$ and $\varphi \neq 0$ no estimate of the form
\be\label{nosubnge2}
\| H_0^{\f 3 4 +\ep }f\| \le C (\|\widetilde {P}_0( \theta) f\| +\|f\|), \qquad f \in \vS(\bR^{2n}),
\ee
holds, for any $\ep>0$ and $C>0$. Estimate \eqref{sub} is also used to prove a pseudospectral estimate, from which one obtains
\be
\sigma(P_0(\theta)) = \bigcup_{k\in\bN}\{k+ e^{-2i\varphi}r;\ r \ge 0\} .
\ee
 For the KFP operator with a decaying force, the eigenvalues $\vT$ of the diffusion operator play a special role and are called {\it thresholds}. With the present normalization of physical constants, $\vT = \sigma(H_0) =\bN$. Under condition \eqref{ass1}, $W_\theta$ is relatively compact with respect to $P_0(\theta)$. It follows that for $0 <\varphi <a$ the resolvent $R(z, \theta) =(P(\theta)  -z)^{-1}$, which is meromorphic in $0< \arg  z < \f{3\pi} 2$, admits a meromorphic continuation across the positive semi-axis into the region
\[
\vZ_\varphi = \{ z\in \bC;\ -2\varphi <\arg z \le 0\}\setminus \bigcup_{k\in\bN}\{k+ e^{-2i\varphi}r;\ r \ge 0\}
\]
of the lower half-plane. The set of resonances of the KFP operator $P$ is defined as
\be \label{resdefn}
{\rm Res}\,(P) = \bigcup_{ 0 <\varphi <a}\Big\{z_0 \in \vZ_\varphi;\ R(z, \theta) \mbox{ has a pole at } z=z_0 \Big\}.
\ee
Hence the resonances of $P$ in $\vZ_\varphi$ are discrete eigenvalues of $P(\theta)$ with $\varphi = \Im \theta >0$. The above definition includes the real poles of $R(z, \theta)$ in $\bR_+^* \setminus \vT$. It is known \cite{w4} that these poles are positive eigenvalues or real resonances of $P$ in a weighted $L^2$-space, satisfying an outgoing radiation condition discussed in \cite{w2}.
The domain in which resonances of the KFP operator are defined resembles that for the Schr\"odinger operator with constant magnetic field in dimension three \cite{bbr} where thresholds are the Landau levels.

As in quantum scattering, the resonant states of the KFP operator are metastable states. The real part of a resonance is related to the Kramers escape rate of a Brownian particle, and its inverse gives the lifetime of the particle (see \cite{mel}). In quantum
scattering, by contrast, it is the imaginary part of a resonance which is the inverse of the lifetime
\cite{dz}. For the mathematical analysis of the low-lying
eigenvalues of the KFP operator and of the Eyring--Kramers law, which are also related to the escape rates, see \cite{hln,hn,hss,hhs,hhs2,blm}.

The goal of this work is to study the resonances of $P$ in the large-friction limit.
Let $P(\gamma) = \gamma H_0 + L$, where the friction coefficient $\gamma$ is regarded as a large parameter, $P(\gamma, \theta) = \gamma H_0 + L_\theta$ and $P_0(\gamma, \theta ) = \gamma H_0 + e^{-\theta} v\cdot\na_x$. Then
\[
\sigma (P_0(\gamma, \theta )) = \bigcup_{k\in\bN}\{\gamma k+ e^{-2i\varphi}r;\ r \ge 0\}.
\]
 Under condition \eqref{ass1}, the resonances of $P(\gamma)$ are defined as in \eqref{resdefn}, the set of thresholds $\vT=\bN$ being replaced by $\vT_\gamma=\gamma \bN$.

The main result of this work describes the eigenvalues and resonances of $P(\gamma)$ in a
neighborhood of an arbitrary threshold $\gamma k$, $k \in \bN$, in terms of an operator
acting in the $x$-variables only, the Witten Laplacian at the level $k$. Let $\Pi_k$ be the eigenprojection of $H_0$ for the eigenvalue
$k$. The eigenspace $\Ran \Pi_k$ is spanned by the Hermite functions $\psi_\alpha$, $|\alpha| = k$,
defined in (\ref{hermite}). We use throughout the identification
\bea \label{ident}
\bC^{\mu_k} \ni u = (u_\alpha)_{|\alpha|=k} \ &\longleftrightarrow &\ \sum_{|\alpha|=k}u_\alpha
\psi_\alpha \in \Ran \Pi_k, \\
\mu_k = \mbox{\rm rank}\,\Pi_k &= & \# \{\alpha \in \bN^n; |\alpha|=k\} = \binom{n+k-1}{k} . \nn
\eea
The components of a vector of $\bC^{\mu_k}$ are thus numbered by the multi-indices $\alpha \in
\bN^n$ with $|\alpha| = k$. The corresponding matrices have entries $M_{\beta\alpha}$,
$|\alpha| = |\beta| = k$, where $M_{\beta\alpha}$ is the coefficient of $\psi_\beta$ in
$M\psi_\alpha$. With this choice of basis, for an $n\times n$ complex matrix $M = (M_{lm})$, the
operator $\dG_k(M)$ on $\bC^{\mu_k}$ is defined by
\be \label{dGkintro}
\left( \dG_k(M)\right)_{\beta\alpha} =
\left\{ \begin{array}{ll}
\sum_l M_{ll}\,\alpha_l , & \beta = \alpha, \\[2mm]
M_{lm}\sqrt{\alpha_m}\sqrt{\alpha_l+1}, & \beta = \alpha - e_m + e_l, \ l \neq m, \\[2mm]
0 & \mbox{ otherwise,}
\end{array}\right.
\ee
where $(e_1,\dots,e_n)$ is the canonical basis of $\bR^n$.
 Set
\be \label{Qkintro}
-\Delta^{(k)}_V(\theta) := -\Delta_V(\theta)\otimes 1_{\mu_k}
- \dG_k\left( (\mbox{\rm Hess}\,V)(e^\theta x)\right) \quad \mbox{ on } L^2(\bR^n_x;\bC^{\mu_k}),
\ee
which we call the \emph{Witten Laplacian at the level $k$} (Definition \ref{def-Qk} and
Proposition \ref{prop-Qk}). It is the Witten Laplacian dressed by $k$ quanta of the Hessian of $V$. Here
\be \label{Witten}
-\Delta_V(\theta) = - e^{-2\theta}\Delta + \f 1 4 (\nabla V)(e^\theta x) \cdot (\nabla V)(e^\theta x) - \f 1 2 (\Delta V)(e^\theta x).
\ee
Denote $-\Delta^{(k)}_V = -\Delta^{(k)}_V(0)$ and $-\Delta_V = -\Delta_V(0)$. In particular, $-\Delta^{(0)}_V(\theta) = -\Delta_V(\theta)$.
For $k=0$, one has $\mu_0=1$ and $\dG_0 = 0$, and we recover the Witten Laplacian
\be
-\Delta^{(0)}_V = -\Delta_V = -\Delta + \f 1 4 |\nabla V(x)|^2 - \f 1 2 \Delta V(x).
\ee
Thus $-\Delta^{(k)}_V$ is the higher-threshold substitute for the Witten Laplacian. For $n=1$, $\mu_k=1$
and $-\Delta^{(k)}_V = -\Delta_V - kV''(x)$. For $k=1$, $-\Delta^{(1)}_V = -\Delta_V\otimes 1_n - \mbox{\rm Hess}\,V$ on
$L^2(\bR^n;\bC^n)$. This is \emph{not} the Witten Laplacian on $1$-forms,
$-\Delta_V\otimes 1_n + \mbox{\rm Hess}\,V$ \cite{wit,cfks,hs4,hln}: the sign is opposite.

Under condition \eqref{ass1}, the resonances of $-\Delta^{(k)}_V$ can also be defined by analytic dilation \cite{ac}. They are eigenvalues of $-\Delta^{(k)}_V(\theta)$ in $\{-2\varphi<\arg\eta\le 0\}$, $0<\varphi<a$.
The correct spectral parameter relating the KFP operator to the Witten Laplacian at the level $k$, $-\Delta^{(k)}_V$, is given by
\be \label{etaintro}
\eta = \gamma (z-\gamma k) .
\ee
Let $\delta \in\,]0,\f\pi2[$ and
\[
\Lambda_\delta = \left\{ \eta \in \bC\setminus\{0\};\ \big|\arg\big(e^{2i\varphi}\eta\big)\big| \ge \delta\right\} .
\]
It is the complement of the open sector of half-angle $\delta$ around the half-line
$\sigma_{\rm ess}(-\Delta^{(k)}_V(\theta)) = e^{-2i\varphi}[0,\infty[$.

The main results of this work are proved under the following stronger version of \eqref{ass1}: $V$ extends to a holomorphic function $V(\zeta)$ on a domain of the form
\[
\{\zeta \in \bC^n, |\zeta|< C\} \cup \{\zeta \in \bC^n, |\Im \zeta| < \tan a \, |\Re \zeta|\}
\]
for some $a \in ]0, \f \pi 4] $ and $C>0$, and satisfies there
\be \label{ass1b}
|\nabla V(\zeta)| \le C\w{\zeta}^{-\rho}
\ee
for some $\rho>0$. By the Cauchy inequalities, \eqref{ass1b} implies that for every $\alpha \in \bN^n$
\be
|\partial_x^\alpha \nabla V (e^\theta x) | \le C_\alpha \w{x}^{-|\alpha|-\rho}, \qquad
x\in \bR^n,\ \theta\in \bC,\ |\theta| < a .
\ee
We denote by $\delta_V>0$ the supremum of the $\delta>0$ such that $V$ extends holomorphically to the union of the above domain and
of the strip $\{\zeta\in\bC^n;\ |\Im\zeta|<\delta\}$; see \eqref{deltaV}.

\begin{theorem}\label{th-mainintro} Assume condition \eqref{ass1b} and let $\theta$ be fixed with
$\varphi \in\,]0,a[$, $k \in \bN$ and $\delta\in\,]0,\f\pi2[$. There exist $c_1>0$, $C_0>0$
and $C>0$ such that the following holds.

\textup{(a)} There is $\gamma_0>0$ such that for $\gamma\ge\gamma_0$, $R(z,\gamma,\theta) = (P(\gamma,\theta)-z)^{-1}$ is holomorphic and
$\|R(z,\gamma,\theta)\| \le C\gamma$ in
\[
\vD^{(k)}_\gamma = \left\{ z\in\bC;\ \f{C_0}\gamma \le |z-\gamma k| \le c_1,\
\gamma(z-\gamma k) \in \Lambda_\delta \right\} .
\]
In particular $P(\gamma)$ has neither eigenvalues nor resonances in $\vD^{(k)}_\gamma$.

\textup{(b)} Let $\eta_0 \in \Lambda_\delta$ be an eigenvalue of
$-\Delta^{(k)}_V(\theta)$ of multiplicity\footnote{In this work the multiplicity of an eigenvalue
refers to its algebraic multiplicity.} $k_0$. Let $\ep_0>0$ be such that
$\{|\eta-\eta_0|\le2\ep_0\}\cap\sigma(-\Delta^{(k)}_V(\theta)) = \{\eta_0\}$ and
$\{|\eta-\eta_0|\le2\ep_0\}\subset\Lambda_{\delta/2}$. Then there is
$\gamma_1$ such that, for $\gamma\ge\gamma_1$, $P(\gamma,\theta)$ has exactly
$k_0$ eigenvalues, counted with  multiplicity, in the disc
$\{|z-\gamma k - \eta_0\gamma^{-1}| < \ep_0\gamma^{-1}\}$. These eigenvalues $z_l(\gamma)$, $1\le l\le k_0$, satisfy, for every
$N\in\bN^*$,
\be \label{expintro}
z_l(\gamma) = \gamma k + \f{\eta_0}\gamma + \sum_{m=1}^N c_{l,m}\gamma^{-1-\f m{k_0!}}
+ O\left( \gamma^{-1-\f{N+1}{k_0!}}\right), \qquad \gamma \to \infty .
\ee
\textup{(c)} Let $r>0$ and assume that $-\Delta^{(k)}_V(\theta)$ has no eigenvalue on the
boundary of $\Lambda_\delta\cap\{|\eta|\ge r\}$.
Then, for $\gamma \ge \gamma_2(r)$, in the truncated annulus
\[
\vA_\gamma(r) = \left\{ z;\ \f r\gamma \le |z-\gamma k|\le c_1,\
\gamma(z-\gamma k)\in\Lambda_\delta\right\} ,
\]
the eigenvalues of $P(\gamma, \theta)$, counted with multiplicity, are in bijection with those  of $-\Delta^{(k)}_V(\theta)$ in $\Lambda_\delta\cap\{r \le |\eta| \le c_1\gamma\}$, and they are located at $z = \gamma k + \eta_j\gamma^{-1} + O(\gamma^{-1-1/k_j})$, where $\eta_j$ is some eigenvalue of $-\Delta^{(k)}_V(\theta)$ and $k_j$ the
multiplicity of $\eta_j$.
\end{theorem}

 Under the conditions of Theorem \ref{th-mainintro}, the resonant states can be expanded asymptotically in powers of $\gamma^{-1}$. The denominator $k_0!$ in the exponents of \eqref{expintro} results from the finite-dimensional perturbation theory of eigenvalues in its general form \cite{k}. We also prove that under condition \eqref{ass1} there are neither eigenvalues nor resonances in
domains away from the thresholds $\vT_\gamma$ (Theorem \ref{th2.5}). Together, Theorem \ref{th-mainintro} and Theorem
\ref{th2.5} cover large regions around $\gamma k$ such as
\[
\begin{aligned}
\Big\{z;\ &\f{r}\gamma < |z-\gamma k| < (1-\delta) \gamma,\
-2\varphi +\delta < \arg (z-\gamma k) <\f{3\pi}2, \\
&{\rm dist}\Big(z, \bigcup_{j<k}\big(j\gamma +e^{-2i\varphi}\bR_+\big)\Big) >\delta \gamma \Big\}, \qquad r, \delta >0,
\end{aligned}
\]
in which all the poles, except those induced by the negative eigenvalues of $-\Delta^{(k)}_V$, $k\ge 1$,  lie in the lower half-plane, hence are resonances of $P(\gamma)$.

Part (a) of Theorem \ref{th-mainintro} is Corollary \ref{th-absence}, part (b) consists
of Theorems \ref{th-bij} and \ref{th3.10}, and part (c) is Corollary \ref{cor-sharp}.
They are deduced
from the stability estimate, Theorem \ref{th-red}, and, for part (b), from the reduced determinant of
Proposition \ref{prop-det}.
For $k\ge 1$, $-\Delta^{(k)}_V$ may have negative eigenvalues (Theorems \ref{th-neg} and \ref{th-largeg}).
By dilation analyticity, they are eigenvalues of $-\Delta^{(k)}_V(\theta)$ for any $\theta$ with $0<\varphi <a$.

\begin{theorem}\label{th-introneg} Assume \eqref{ass1b}, let $k\ge 1$ and let $\eta_0 <0$ be an eigenvalue of $-\Delta^{(k)}_V$ with multiplicity $k_0$.
Then the $k_0$ poles of $R(z, \gamma, \theta)$, counted with multiplicity, near $\gamma k + \f{\eta_0}\gamma $ satisfy
\[
z(\gamma) = \gamma k + \eta_0\gamma^{-1} + O(\gamma^{-2}).
\]
Their expansions are in integer powers of $\gamma^{-1}$, with real coefficients, see \eqref{intpow}. Moreover, for every $c<2\delta_V$, there is $C_c$ such that
$|\Im z(\gamma)|\le C_ce^{-c\gamma}$.
\end{theorem}

Theorem \ref{th-introneg} is  Corollary \ref{cor-intpow} for the expansion in powers of $\gamma^{-1}$ and Theorem \ref{th-real} for the exponential bound on the imaginary parts.

\medskip

The main difficulties of this work are related to the loss of $\f 13$ derivative of the KFP
operator. This results in the ill-posedness of the Grushin problem, and in domain problems. This Grushin problem
reduces the study of $P(\gamma,\theta)$ to that of an operator $\vE_{-+}$ acting in
the $x$-variables only. The operator $\vE_{-+}$ is the sum of a matrix-valued differential operator of order one
and of a remainder, which is a pseudodifferential operator of order $\f 43$. Hence, even at a
formal level, it is not clear why the first-order part determines the eigenvalues of
$P(\gamma, \theta)$.

Two results are basic to this work. The first one is the spectral
decomposition of the fibre resolvent, Lemma \ref{lem2.3}, proved by the time-dependent method already used
in \cite{w1} for $\theta=0$ and $\gamma=1$.
Lemma \ref{lem2.3} makes it possible to separate the levels and to
localize the frequencies to the region $|\xi| \le \gamma^\sigma$, $\f 34<\sigma<1$ (Corollary \ref{cor2.4} and
Lemma \ref{lem-loc}). The absence of eigenvalues and resonances away from the thresholds $\vT_\gamma$
(Theorem \ref{th2.5}) then follows easily.
The second one is a stability estimate, Theorem \ref{th-red}, saying that the eigenvalues of $P(\gamma,\theta)$ near $\gamma k$ can only
appear near the points $\gamma k + \eta_0\gamma^{-1}$, where $\eta_0$ is an eigenvalue of
$-\Delta^{(k)}_V(\theta)$. It is obtained by combining the Grushin reduction with the frequency
localization.
To show that poles of $R(z,\gamma,\theta)$ do appear near these points, we eliminate the levels $j\neq k$ by successive approximation. Each step is a Schur complement, and
the Witten Laplacian at the level $k$, $-\Delta^{(k)}_V(\theta)$, appears in the Schur complement at the level $k$.
The localization of resonances and their large-friction expansions are obtained by means of a Weinstein--Aronszajn determinant (see \cite{k}).
Finally, we show that the eigenvalues of $P(\gamma, \theta)$ produced by the negative eigenvalues of $-\Delta^{(k)}_V$ have exponentially small imaginary parts.
For this, we deform $P(\gamma)$ in the frequency variable $\xi$ at the frequencies $|\xi|\ge c\gamma$.
A flux identity, coming from the reality of $P(\gamma)$ and from its parity
symmetry $\bJ$ defined by (\ref{bJ}), then expresses their imaginary parts through the difference of two opposite deformations, which gives the exponential smallness.
\medskip

The organization of this paper is as follows. In Section \ref{sect2}, we establish a spectral
decomposition for the resolvent of
$\widehat P_0(\xi,\gamma,\theta) = \gamma H_0 + ie^{-\theta}v\cdot\xi$. We also prove estimates with frequency localization near each channel, which are used throughout. We then derive the absence
of eigenvalues and resonances away from the thresholds for $\gamma$ large. Section \ref{sect3}
contains the analysis of the ill-posed Grushin problem with the first $k+2$ eigenvalues of $H_0$, which gives rise to  the Witten Laplacians at the level $k$.  Section \ref{sect4} treats the poles near an arbitrary threshold $\gamma k$, $k\in\bN$.
  Section \ref{sec-poles} studies the poles produced by negative eigenvalues of $-\Delta_V^{(k)}$, $k \ge 1$, which lie
below the thresholds $\gamma k$. To overcome the difficulties due to high frequencies arising from analytic dilation, we
 use instead the method of phase-space deformation to define resonances.
Finally, Appendix \ref{app-Gr} contains the proof of Theorem \ref{th3.2} on the inverse of the Grushin problem, and Appendix \ref{app-free} a large-friction expansion of the dilated free resolvent, which is used to obtain the complete asymptotic expansion of the poles as $\gamma \to\infty$.

\medskip
\noindent{\bf Notation and conventions.} For $r\ge 0$ and $s \in \bR$, we introduce the weighted Sobolev space
$\vH^{r,s}$ adapted to the subellipticity of $P_0$:
\[
\vH^{r,s} = \{u \in \vS'(\bR^{2n}), (1 + \w{D_v}^2 + |v|^2 + \w{D_x}^{\f 2 3})^{\f r 2 } \w{x}^s u \in L^2\}.
\]
For $r <0$ and $s\in \bR$, $\vH^{r, s}$ is the dual space of $\vH^{-r, -s}$. We write
$\vH = \vH^0= L^2(\bR^{2n})$ and

\[
\vH^s =\vH^{0,s}= L^2(\bR^{2n}, \w{x}^{2s} dxdv).
\]

$\vL(r,s, r',s')$ is the space of bounded linear operators from $\vH^{r,s}$ to $\vH^{r',s'}$, and
$\vL(s,s') = \vL(0,s, 0,s')$. For example, $(P_0+1)^{-m}\in \vL(r, s, r+2m, s)$ for every
$m\in \bN^*$ and $r,s \in \bR$. For $s\in\bR$, $H^s = H^s(\bR^n)$, or $H^s(\bR^n;\bC^\mu)$, is the usual Sobolev space of order $s$, and
$H^\infty := \bigcap_sH^s$. The arguments $z$, $\gamma$, $\theta$, $\xi$ of operators such as $R(z,\gamma,\theta)$,
$R'(z,\gamma,\theta)$, $\vE_{-+}(z,\gamma,\theta)$ or $\widehat R_0(z,\xi,\gamma,\theta)$ are written in the
definitions and in the statements. They are omitted in the proofs when no confusion can arise. For example, in a proof, $R = R(z,\gamma,\theta)$ and $P = P(\gamma,\theta)$. They are also omitted when the corresponding parameter is fixed.
Throughout, $C$ denotes a positive constant which may change from
line to line. The scalar product $\w{\cdot, \cdot}$ is linear in the second argument.

\sect{Frequency localizations} \label{sect2}

Consider the KFP operator in the large-friction regime
\[
P(\gamma) = \gamma H_0 + L = \gamma \left(  -\Delta_v+\f{1}{4}|v|^2-\f{n}{2}\right ) + v\cdot\na_x-\na V(x)\cdot\na_v ,
\]
where the friction coefficient $\gamma \ge 1$ is regarded as a large parameter. We assume the
dilation-analyticity condition \eqref{ass1} for some $a>0$. For $\theta \in \bC$ with
$|\varphi| < a$, let
\[
P(\gamma, \theta) =\gamma H_0 + L_\theta = P_0(\gamma, \theta) + W_\theta ,
\]
where
\begin{align*}
 P_0(\gamma, \theta) &= \gamma H_0 + e^{-\theta} v\cdot\na_x , \\
 L_\theta &=  e^{-\theta} v\cdot\na_x-\na V(e^{\theta} x)\cdot\na_v , \\
W_\theta &= -\na V(e^{\theta} x)\cdot\na_v .
\end{align*}
After a Fourier transform in $x$, $P_0(\gamma,\theta)$ becomes the family
\[
\widehat{P}_0(\xi, \gamma, \theta) = \gamma H_0 + i e^{-\theta} v\cdot \xi, \qquad \xi\in\bR^n,
\]
and we set

\[
\widehat R_0(z, \xi, \gamma, \theta) = ( \widehat{P}_0(\xi, \gamma, \theta)-z)^{-1}
\]
 for
$z\not\in \sigma ( \widehat{P}_0(\xi, \gamma, \theta) )$, and $R_0(z,\gamma,\theta) := (P_0(\gamma,\theta)-z)^{-1}$.

\subsection{Spectral decomposition of the fibred free resolvent}

\begin{theorem} \label{th2.1} For any $\phi_0\in \,]0, \f \pi 2[$, there exist positive constants $C$ and $s_0$ such that
\be \label{e2.2b}
\| (s + \gamma H_0)f\| +  \gamma^{\f 1 3}\w{\gamma +|\xi|}^{\f 2 3 }\|f\|  + \| v\cdot \xi f\| \le C \|(\widehat{P}_0(\xi, \gamma, \theta) +s) f\|,
\ee
for $\theta \in \bC$ with $|\theta|\le \phi_0$, $s \ge \gamma s_0$, $\xi\in \bR^n$, $\gamma \ge 1$ and $f \in \vS(\bR_v^{n})$, where $\|\cdot\|$ is the norm of $L^2(\bR^n_v)$.
\end{theorem}
\pf
One has
\[
\widehat{P}_0(\xi, \gamma, \theta) +s =\gamma (\widehat{P}_0(\xi',  \theta) +s')
\]
 with $\xi' = \gamma^{-1}\xi$ and $s' = \gamma^{-1}s$, and \eqref{e2.2b} follows from \eqref{sub} applied to $\widehat{P}_0(\xi',  \theta) +s'$.
\ef

Taking $s = \gamma s_0$ in Theorem \ref{th2.1}, one sees that for any $M>0$ there exists $R>0$ such that
\be \label{e4.26}
\gamma \| (H_0+s_0)  \widehat R_0(z, \xi, \gamma, \theta)\| +  \gamma^{ \f 1 3} \w{\xi}^{\f 2 3} \|\widehat R_0(z, \xi, \gamma, \theta)\| \le C
\ee
uniformly in $\gamma \ge 1$, $\xi\in \bR^n$ with $|\xi| \ge R \gamma $ and $|z|\le M \gamma$. In particular, for
$z \not\in \sigma(P_0(\gamma, \theta))$,
\begin{align*}
 \|\widehat R_0(z, \xi, \gamma, \theta)\| &\le C \gamma^{-1},\\
 \w{\xi}^{\f 2 3} \|\widehat R_0(z, \xi, \gamma, \theta)\| &\le C \gamma^{-\f 1 3}, \\
 \|v  \widehat R_0(z, \xi, \gamma, \theta)\| +  \|  \nabla_v \widehat R_0(z, \xi, \gamma, \theta)\|  & \le C \gamma^{-\f 1 2} ,
\end{align*}
uniformly in $ \xi \in \bR^n$ with $|\xi| >R\gamma$, $|\theta| \le \phi_0$ and $\gamma \ge \f{|z|} M +1$.

Let $M>0$ and let $R=R_M$ be such that \eqref{e4.26} holds. Introduce a partition of unity
\be \label{pu}
\chi_0(\xi) +\chi_1(\xi) =1, \qquad  \xi\in \bR^n,
\ee
with $0 \le \chi_j \le 1$, $\chi_0\in C_0^\infty(\bR^n)$ and $\chi_0(\xi ) =1$ for $|\xi| \le R$. Let $D = i^{-1} \partial_x$.
Theorem \ref{th2.1} gives
\be \label{highfrequencies}
\| R_0(z, \gamma, \theta) \chi_1(\gamma^{-1}D)\| \le C\gamma^{-1}
\ee
uniformly for $|z|\le M \gamma$.
Using the properties of the Hermite functions, one checks as in \cite{w1} that
\[
\sigma( \widehat{P}_0(\xi, \gamma, \theta) ) =\left\{E_k^\theta(\xi, \gamma);\ k \in \bN\right \},
\qquad
E_k^\theta(\xi, \gamma) = \gamma k + \gamma^{-1} e^{-2\theta} \xi^2 .
\]
The functions $E_k^\theta(\xi, \gamma)$ are called the \emph{channels}.

Let $\zeta = \gamma^{-1} e^{-\theta} \xi$, and write $\zeta = \zeta_1 + i \zeta_2$ with $\zeta_j \in \bR^n$. The Riesz projection $\Pi^\zeta_k$ associated with the eigenvalue $E_k^\theta(\xi, \gamma)$ can be computed in terms of Hermite functions. To this end, let
\be \label{creation}
a_l = \p_{v_l} + \f{v_l}2, \qquad a^*_l = -\p_{v_l} + \f{v_l}2, \qquad 1\le l \le n,
\ee
be the annihilation and creation operators associated with $H_0$. They satisfy
$[a_l,a_m] = [a^*_l,a^*_m]= 0$, $[a_l,a^*_m] = \delta_{lm}$. Moreover,
\be \label{aa}
H_0 = \sum_{l=1}^n a^*_la_l, \qquad v_l = a_l + a^*_l, \qquad
\p_{v_l} = \f 12 \left( a_l - a^*_l\right).
\ee
Let $\psi_0(v) = (2\pi)^{-\f n4}e^{-\f{|v|^2}4}$ be the normalized ground state, $a_l\psi_0 = 0$ for
all $l$, and, for $\alpha \in \bN^n$,
\be \label{hermite}
\psi_\alpha = \f{(a^*)^\alpha}{\sqrt{\alpha!}}\,\psi_0
= \f{(a^*_1)^{\alpha_1}\cdots (a^*_n)^{\alpha_n}}{\sqrt{\alpha_1!\cdots\alpha_n!}}\,\psi_0 .
\ee
Then $\{\psi_\alpha; \alpha \in \bN^n\}$ is an orthonormal basis of $L^2(\bR^n_v)$,
$H_0\psi_\alpha = |\alpha|\psi_\alpha$ and
\be \label{aaction}
a^*_l\psi_\alpha = \sqrt{\alpha_l+1}\,\psi_{\alpha+e_l}, \qquad
a_l\psi_\alpha = \sqrt{\alpha_l}\,\psi_{\alpha-e_l}, \qquad 1 \le l \le n .
\ee
In particular, $\Ran \Pi_k = {\rm span}\{\psi_\alpha; |\alpha| = k\}$ has dimension
$\mu_k = \binom{n+k-1}{k}$. By \eqref{aa}, both $v_l$ and $\p_{v_l}$ map $\Ran \Pi_j$ into
$\Ran \Pi_{j-1}\oplus \Ran\Pi_{j+1}$.
With this notation, the Riesz projection $\Pi^\zeta_k$ associated with the eigenvalue
$E^\theta_k(\xi,\gamma)$ is given by
\begin{equation}\label{Riesz}
\Pi_k^\zeta =\sum_{|\alpha| =k} \w{\overline{\psi_\alpha(\cdot + 2i \zeta)}, \ \cdot\ } \psi_{\alpha}(\cdot + 2i\zeta),
\end{equation}
and $\Pi_k = \Pi_k^0$ is the eigenprojection of $H_0$ associated with the eigenvalue $k$.

The norms of the projections $\Pi^\zeta_k$ enter every estimate of this section, and
they grow with $|\zeta|$. We collect here the bounds which are used. They follow from the explicit
computation of $\|\psi_\alpha(\cdot+2i\xi)\|$ made in \cite[Appendix A]{w1}. Recall that

\[
\mu_k = \#\{\alpha \in \bN^n; |\alpha| = k\} \le C_n(1+k)^{n-1}.
\]

\begin{lemma}\label{lem-proj} \textup{(i)} For $\zeta = \zeta_1 + i \zeta_2 \in \bC^n$ with
$\zeta_1, \zeta_2 \in \bR^n$, one has $\|\Pi^\zeta_k\| = \|\Pi^{\zeta_1}_k\|$.

\textup{(ii)} For $\xi \in \bR^n$,
\be \label{projnorm}
\left\| \Pi^\xi_k\right\| \ \le\ \varpi_k(\xi) := \sum_{|\alpha| = k} e^{2|\xi|^2}\prod_{m=1}^n
\left( \sum_{j=0}^{\alpha_m}\binom{\alpha_m}{j}\f{(4\xi_m^2)^j}{j!}\right)
\ \le\ \mu_k\, e^{2|\xi|^2 + 4|\xi|\sqrt k} .
\ee

\textup{(iii)} For $t>0$ and $\xi \in \bR^n$,
\be \label{sumproj}
\sum_{k=0}^\infty e^{-t(k+|\xi|^2)}\,\varpi_k(\xi) \ \le\
\f{e^{-|\xi|^2\left( t - 2 - \f 4{e^t-1}\right)}}{(1-e^{-t})^n} \, .
\ee

\textup{(iv)} For every $\ell \in \bN$, $R>0$ and $M>0$, there exists $C>0$ such
that
\be \label{tailproj}
\sum_{k=\ell+1}^\infty e^{-t(k+|\xi|^2)}\,\varpi_k(\xi) \ \le\ C\,e^{-(\ell+1)t},
\qquad t \ge M, \quad |\xi| \le R .
\ee
\end{lemma}
\pf The first inequality in \eqref{projnorm} and \eqref{sumproj} are proved in \cite[Appendix A]{w1}.
Part (i) holds because $\Pi^\zeta_k$ is conjugate to $\Pi^{\zeta_1}_k$ by the unitary translation
$v\mapsto v-2\zeta_2$. The second inequality in \eqref{projnorm} follows from
$\binom{p}{j}\le p^j/j!$ and the Cauchy--Schwarz inequality, and \eqref{tailproj} follows from
\eqref{projnorm} and $\mu_k\le C_n(1+k)^{n-1}$. \ef

The following result is a spectral decomposition of the fibred resolvent
$\widehat R_0(z, \xi, \gamma, \theta)$;
see \cite{w1} for the case $\theta =0$, $\gamma =1$.

\begin{lemma}\label{lem2.3} Let $\ell \in \bN$, $b\in\,]\ell, \ell +1[$ and $\phi_0 \in \,]0,  \f \pi 2[$.
Then there exists $\chi_0 \in C_0^\infty(\bR^n_\xi)$ such that
\be \label{R0}
\widehat R_0(z, \xi, \gamma, \theta) = \sum_{k=0}^\ell \chi_0(\gamma^{-1} \xi) \frac{\Pi_k^\zeta}{E_k^\theta(\xi, \gamma)-z} + r_\ell(z,\xi, \gamma, \theta),
\ee
for all $\xi \in \bR^n$ and $z\not\in \sigma ( \widehat{P}_0(\xi, \gamma, \theta) )$ with $|z| < b \cos \varphi\, \gamma $.
Here $r_\ell(z,\xi, \gamma, \theta)$ is holomorphic in $z$ in the disc $\{z;\ |z| <   b \cos \varphi \, \gamma\}$ and satisfies
\be \label{reste}
\|r_\ell(z, \xi, \gamma, \theta)\| \le C\gamma^{-1}
\ee
and
\be \label{reste1}
 \|\widehat{P}_0(\xi, \gamma, \theta)r_\ell(z, \xi, \gamma, \theta)\| \le C,
\ee
uniformly in $z\in \bC$ with $ |z| <   b \cos \varphi \, \gamma$, $|\theta|\le \phi_0$, $ \xi \in \bR^n$ and $ \gamma \ge 1$.
\end{lemma}
\pf To simplify the notation, we only give the proof for $\theta = i \varphi$ with $ \varphi>0$. The case of general $\theta$ is similar.
Let $\chi_0(\xi) + \chi_1(\xi) =1$ be the partition of unity \eqref{pu} with $M>\ell+1$,
so that $\chi_0$ is supported in $\{|\xi| \le R\}$ for some $R = R_M$ and $\chi_1(\gamma^{-1}\xi)$
is supported in $\{|\xi| \ge R\gamma\}$.
By \eqref{e4.26},
\[
\widetilde r_{\infty} = \chi_1(\gamma^{-1} \xi) \widehat R_0
\]
 satisfies the remainder estimate \eqref{reste}.
We write $\xi' = \gamma^{-1}\xi$ and $t' = \gamma\cos\varphi\,t$, so that
$\zeta = \gamma^{-1}e^{-i\varphi}\xi$ has real part $\zeta_1 = \cos\varphi\,\xi'$ and
\be \label{ReEk}
\Re\left( e^\theta E^\theta_k(\xi,\gamma)\right)
= \cos\varphi\left( \gamma k + \gamma^{-1}|\xi|^2\right), \qquad
t\,\Re\left( e^\theta E^\theta_k(\xi,\gamma)\right) = t'\left( k + |\xi'|^2\right).
\ee
Since $v\cdot\xi$ is symmetric,
\[
\Re \w{f, e^\theta \widehat P_0f}
= \gamma\cos\varphi\,\w{f,H_0f} \ge 0,
\]
so that $-e^\theta\widehat P_0$ generates a
contraction semigroup, and for $\Re (e^\theta z) <0$
\[
\chi_0(\gamma^{-1} \xi) \widehat R_0 =e^{\theta} \chi_0(\gamma^{-1} \xi)  \int_0^\infty e^{-te^\theta (\widehat P_0  -z)} dt.
\]
Let
\[
\widetilde{r}_0 = e^{\theta} \chi_0(\gamma^{-1} \xi) \int_{0}^{M\gamma^{-1}} e^{-te^\theta (\widehat P_0  -z)} dt .
\]
Since $\Re \left(e^\theta \widehat P_0 \right) \ge 0$, the
integrand has norm at most $e^{t\Re(e^\theta z)}$, and since $\Re (e^\theta z) \le |z| \le M \gamma$
\[
\left\| \widetilde{r}_0\right\| \ \le \ \int_0^{M\gamma^{-1}}e^{t\Re(e^\theta z)}\,dt \ \le \
\f M\gamma\, e^{M^2}
\]
uniformly in $\xi\in \bR^n$ and $\gamma \ge 1$. Therefore $\widetilde{r}_0 $ also satisfies \eqref{reste}.

By Lemma \ref{lem-proj}(i)--(ii) and
\[
|(\zeta_1)_m| = \cos\varphi\,|\xi'_m| \le |\xi'_m|,
\]

\be \label{Pizeta}
\left\| \Pi^\zeta_k\right\| = \left\| \Pi^{\zeta_1}_k\right\| \le \varpi_k(\zeta_1) \le \varpi_k(\xi') ,
\ee
so that, by \eqref{ReEk} and Lemma \ref{lem-proj}(iii),
\begin{align} \label{sumPi}
\sum_{k=0}^\infty e^{-t\Re (e^\theta E_k^\theta(\xi, \gamma) ) }\|\Pi_k^\zeta\|
&
\ \le \ \f{e^{-|\xi'|^2\left(t'-2 - \f{4}{e^{t'}-1}\right) }}{(1-e^{-t'})^n}, \qquad t'>0.
\end{align}
Since the right-hand side is finite for every $t>0$, the spectral decomposition
\be \label{semigroupdec}
e^{-t e^\theta \widehat P_0 } = \sum_{k=0}^\infty e^{-t e^{\theta} E_k^\theta(\xi, \gamma)}\Pi_k^\zeta
\ee
is norm-convergent for $t>0$. Indeed, completing the square,
\be \label{square}
\widehat P_0 = \gamma\left( H_0 + i v\cdot \zeta\right), \qquad
E^\theta_k(\xi,\gamma) = \gamma\left( k + \zeta^2\right), \qquad \zeta^2 = \sum_m \zeta_m^2 ,
\ee
so that \eqref{semigroupdec} is the statement of \cite[Proposition A.3]{w1} for the complex
parameter $\zeta$. Its proof applies verbatim. Moreover, \eqref{semigroupdec} may be
integrated term by term on $[M\gamma^{-1},\infty[$ against $e^{te^\theta z}\,dt$. Indeed, on that interval,
$t' \ge M\cos\varphi$. Hence the right-hand side of \eqref{sumPi} is bounded for $\xi \in \supp\chi_0(\gamma^{-1}\cdot)$, while $\Re(e^\theta z)<0$.
One has
\[\chi_0(\gamma^{-1} \xi) e^{\theta} \int_{M\gamma^{-1}}^\infty e^{t e^\theta z} \left(\sum_{k=0}^\ell e^{-t e^\theta E_k^\theta(\xi, \gamma) } \Pi_k^\zeta\right) dt
=  \sum_{k=0}^\ell  \frac{\Pi_k^\zeta \chi_0(\gamma^{-1} \xi)}{E_k^\theta(\xi, \gamma)-z} - r_{\ell,0} ,
\]
where
\[
r_{\ell,0}= \chi_0(\gamma^{-1} \xi) e^{\theta} \int_0^{M\gamma^{-1}} e^{t e^\theta z} \left(\sum_{k=0}^\ell e^{-t e^\theta E_k^\theta(\xi, \gamma) } \Pi_k^\zeta\right) \, dt .
\]
On the support of $\chi_0(\gamma^{-1}\cdot)$ one has $|\xi'| \le R$, so that

\[
\|\Pi^\zeta_k\| \le \varpi_k(\xi') \le C
\]
 for $k \le \ell$ by \eqref{Pizeta} and
\eqref{projnorm}. The exponentials are bounded by $e^{t\Re(e^\theta z)} \le e^{M^2}$ on
$[0,M\gamma^{-1}]$. Hence $\|r_{\ell,0}\| \le C\gamma^{-1}$, that is, $r_{\ell,0}$
satisfies \eqref{reste}. Set
\[
r_{\ell, 1} =\chi_0(\gamma^{-1} \xi) e^{\theta} \int_{M\gamma^{-1}}^\infty e^{t e^\theta z} \left(\sum_{k=\ell + 1}^\infty e^{-t e^\theta E_k^\theta(\xi, \gamma) } \Pi_k^\zeta\right) \, dt.
\]
By \eqref{Pizeta}, \eqref{ReEk} and Lemma \ref{lem-proj}(iv) applied with $M$ replaced by
$M\cos\varphi$,
\begin{align*}
\left\|\sum_{k=\ell + 1}^\infty e^{-t e^\theta E_k^\theta(\xi, \gamma) } \Pi_k^\zeta \right \|
& \le \sum_{k=\ell +1}^\infty e^{-t'(k +|\xi'|^2)} \left \|\Pi_k^{\zeta_1}\right \| \le C_\ell e^{- (\ell + 1) t'},
\end{align*}
for
\[
t' = \gamma \cos \varphi \, t \ge M\cos\varphi
\]
 and $\xi \in \supp \chi_0(\gamma^{-1}\cdot)$. For $\Re (e^{i\varphi} z)  < b \cos \varphi \, \gamma$ with $b <\ell +1$, one has
$e^{t\Re(e^{i\varphi}z)} \le e^{bt'}$, and therefore
\[
\|r_{\ell, 1}\| \le C_\ell \int_{0}^\infty e^{-(\ell + 1 -b) t'} \, dt = \f{C_\ell}{\gamma \cos \varphi \, (\ell + 1 -b)}.
\]
Thus $r_{\ell,1}$ also satisfies \eqref{reste}. This proves \eqref{R0} with
\be \label{reste-decomposition}
r_\ell = \widetilde r_0 +  \widetilde r_\infty - r_{\ell, 0} + r_{\ell, 1}
\ee
satisfying \eqref{reste}.
All the constants above are uniform in $|\theta|\le\phi_0$, and so is \eqref{reste}.
The identity \eqref{R0} has been established for $\Re(e^\theta z)<0$. From the proof one sees that $r_\ell$ is holomorphic for $z$ in the disc
 $\{|z| < b\cos\varphi\,\gamma\}$. By analytic continuation, \eqref{R0} holds
in $\{|z| < b\cos\varphi\,\gamma\}\setminus \sigma(\widehat P_0)$.

Finally, $\widehat P_0\Pi_k^\zeta = E_k^\theta(\xi, \gamma)\Pi_k^\zeta$, since $\Pi_k^\zeta$ is the Riesz projection of
$\widehat P_0$ with eigenvalue $E_k^\theta(\xi, \gamma)$. Applying $\widehat P_0$ to \eqref{R0} one gets
\[
\widehat P_0r_\ell = 1 - \sum_{k=0}^\ell\chi_0(\gamma^{-1}\xi)\Pi_k^\zeta + z\, r_\ell ,
\]
and \eqref{reste1} follows from \eqref{reste}, $|z|\le b\gamma$, and the bound \eqref{Pizeta} of $\Pi_k^\zeta$ on the support of $\chi_0(\gamma^{-1}\cdot)$.
\ef

\begin{rmk}
The identity \eqref{R0} may be differentiated in $\theta$. Indeed
$\widehat P_0(\xi,\gamma,\theta)$ is a holomorphic family of type \textup{(A)} in $\theta$.
Moreover $\Pi^\zeta_k$ and $E^\theta_k(\xi,\gamma)$ are holomorphic in $\theta$ by \eqref{Riesz}, so that
$r_\ell(z,\xi,\gamma,\theta)$ is holomorphic in $\{|\theta|<\phi_0\}$.
Since \eqref{reste} is
uniformly bounded, Cauchy's inequality gives
\be \label{restetheta}
\left\| \p_\theta r_\ell(z,\xi,\gamma,\theta)\right\| \le C\gamma^{-1} ,
\ee
uniformly in $|z| < b'\cos\varphi\,\gamma$, $\xi\in\bR^n$ and $\gamma\ge1$, for every $b'<b$.
\end{rmk}

For $k\in\bN$ let $\Upsilon^{(k)} := \nabla_v\Pi_k$, that is, $\Upsilon^{(k)}_j = \p_{v_j}\Pi_k$, $1\le j\le n$. Since $\p_{v_j}$ maps
$\Ran\Pi_k$ into $\Ran\Pi_{k-1}\oplus\Ran\Pi_{k+1}$, $\Upsilon^{(k)}$ is nilpotent in the sense that
\[
\Upsilon^{(k)}_i \Upsilon^{(k)}_j =0, \qquad 1\le  i, j \le n .
\]
On the support of $\chi_0(\gamma^{-1} \cdot)$, $|\zeta|$ is bounded, and
\be \label{Y}
\|\nabla_v \Pi_k^\zeta - \Upsilon^{(k)}\| \le C |\zeta| .
\ee
Let $\theta \in \bC$ be fixed with $\Im \theta = \varphi \in \,]0, \f \pi 4 [$. For $\delta >0$ set
\begin{align*}
\Omega_0  &= \left\{ z\in \bC;\ 0<|z| < 1,\  -2\varphi +\delta < \arg z <  \f {3\pi }2 \right \},\\[2mm]
\Omega_1(\gamma) &= \left\{z \in \bC;\  \delta < | z | < (1- \delta)\gamma,\   - 2 \varphi +\delta < \arg z< \f{3\pi} 2\right \},
\end{align*}
and, for $k\ge1$,
\be \label{Omega1k}
\Omega_{1,k}(\gamma) = \Big\{ z= k\gamma + w;\ w\in \Omega_1(\gamma),\
{\rm dist}\Big(z, \bigcup_{j=0}^{k-1}\big(j\gamma +e^{-2i\varphi}\bR_+\big)\Big) >\delta \gamma\Big\} ,
\ee
with $\Omega_{1,0}(\gamma) := \Omega_1(\gamma)$.

\begin{cor}[frequency localizations]\label{cor2.4} Let $k\in\bN$. There are $\gamma_0$ and $C$ such that the following holds.

\textup{(a)} For $z\in \Omega_{1,k}(\gamma)$, $\gamma \ge \gamma_0$ and $\xi\in \bR^n$,
\be \label{Ups1}
\left \| \nabla_v \widehat R_0(z, \xi, \gamma, \theta) -    \f{ \chi_0(\gamma^{-1} \xi)} {E_k^\theta(\xi, \gamma)-z}\Upsilon^{(k)} \right \| \le C \gamma^{-\f 1 2} .
\ee

\textup{(b)} Let $ \f 12 < \sigma < 1$. Then there exists $R>1$ such that \eqref{Ups1} holds for $z\in\bC$ with $|z-\gamma k|<1$, $\gamma \ge \gamma_0$ and $\xi\in \bR^n$
with $|\xi| \ge R \gamma^{\sigma}$.
\end{cor}
\pf We apply \eqref{R0} with $\ell$ so large that $\ell \cos \varphi > k+\delta^{-1}$. Using \eqref{reste1}, one sees that
 the remainder $r_{\ell}(z, \xi, \gamma, \theta)$ satisfies
\[
\|\nabla_v r_\ell\|\le C \gamma^{-1}
\]
for $\xi \in \bR^n$ and $|z|<(k+\delta^{-1})\gamma$. Indeed, as in the proof of Lemma \ref{lem2.3},
$\gamma\cos\varphi\,\|H_0^{\f12}r_\ell f\|^2=\Re\w{r_\ell f,e^\theta\widehat P_0r_\ell f}\le C\gamma^{-1}\|f\|^2$ by \eqref{reste} and \eqref{reste1}, so that $\|\nabla_vr_\ell\|\le C\|(1+H_0)^{\f12}r_\ell\|\le C\gamma^{-1}$. For $j\neq k$, one has $|E^\theta_j(\xi,\gamma)-z|\ge c\gamma$. In case (a), this follows from the definition of
$\Omega_{1,k}(\gamma)$ for $j<k$, and from $\Re E^\theta_j(\xi,\gamma)\ge\gamma j$ for $j>k$. In case (b), write $z = \gamma k+w$ with $|w|<1$. Then
$E^\theta_j(\xi,\gamma)-z = \gamma(j-k)+\gamma^{-1}e^{-2\theta}\xi^2-w$ with $\arg(e^{-2\theta}\xi^2) = -2\varphi$, so that
$|E^\theta_j(\xi,\gamma)-z|\ge\gamma\sin2\varphi-1$. Since $\nabla_v\Pi^\zeta_j$ is bounded on the support of $\chi_0(\gamma^{-1}\cdot)$, it follows that
\[
\left \|\frac{\chi_0(\gamma^{-1} \xi)}{E_j^\theta(\xi, \gamma)-z} \nabla_v   \Pi_j^\zeta \right \| \le C \gamma^{-1}, \qquad j\neq k,\ j\le\ell .
\]
For the term $j=k$ in \eqref{R0}, we use \eqref{Y} and, with $\eta = \gamma(z-\gamma k)$, the estimate
\[
 \left|\f{|\zeta|}{E_k^\theta(\xi, \gamma)-z}\right | \le \f{C |\xi|}{|e^{-2\theta}\xi^2-\eta|} \le \f{ C'|\xi|}{|\xi|^2 + \gamma}  \le C'' \gamma^{-\f 1 2}
\]
uniformly in $z\in \Omega_{1,k}(\gamma)$, $\xi \in \bR^n$ and $\gamma \ge 1$. This proves (a). Similarly,
\[
 \left |\f{|\zeta|}{E_k^\theta(\xi, \gamma)-z}\right | \le \f{C |\xi|}{e^{-2\Re\theta}|\xi|^2-\gamma } \le C' \gamma^{-\sigma}
\]
uniformly in $|z-\gamma k|<1$, $|\xi| \ge R\gamma^{\sigma}$ and $\gamma \ge \gamma_0$, which proves (b) since $\sigma>\f12$.
\ef

Let $\chi_0(\xi)+\chi_1(\xi) = 1$ be the partition of unity \eqref{pu} and set, for $\sigma \in\,]\f 34,1[$,
\be \label{chisig}
\chi^\sigma_j := \chi_j\big( \gamma^{-\sigma}D_x\big), \qquad j = 0,1 .
\ee

\begin{lemma}\label{lem-loc} Assume \eqref{ass1} and let $k\in\bN$ and $\f 34 < \sigma <1$. Then
\be \label{chi1small}
\left\| \chi^\sigma_1 R_0(z,\gamma,\theta)W_\theta\right\| \le C\gamma^{-\f 12}, \qquad |z-\gamma k|<1, \quad \gamma \ge \gamma_0 .
\ee
\end{lemma}
\pf One has $W_\theta = -\na_v\cdot(\na V)(e^\theta x)$. The proof of Corollary \ref{cor2.4}(b) applies verbatim to
$\widehat R_0\nabla_v$, with $\Upsilon^{(k)}$ replaced by $\Pi_k\nabla_v = -(\Upsilon^{(k)})^*$. It shows that, up to $O(\gamma^{-\f12})$,
$\chi^\sigma_1R_0W_\theta$ is equal to
\[
\vY^\sigma_k(z,\gamma) = \gamma\,\chi^\sigma_1\chi_0(\gamma^{-1}D_x)
\left( -e^{-2\theta}\Delta_x - \eta\right)^{-1}\left( \Upsilon^{(k)}\right)^*\!\cdot(\na V)(e^\theta x) , \qquad \eta = \gamma(z-\gamma k) .
\]
As $|e^{-2\theta}\xi^2-\eta| \ge c|\xi|^2$ for $|\xi| \ge R\gamma^\sigma$ and $|\eta|<\gamma$,
\[
\|\vY^\sigma_k(z,\gamma)\| \le C\gamma^{1-2\sigma}\le C\gamma^{-\f12},
\]

since $\sigma>\f34$. \ef

Estimate \eqref{chi1small} is the only consequence of the frequency localizations used below. On the
range of $\chi^\sigma_1$, the potential is negligible. Hence the analysis near the threshold $\gamma k$ may be performed
on the range of $\chi^\sigma_0$, where $|\xi| \le C\gamma^\sigma$ with $\sigma<1$.

\subsection{Eigenvalue- and resonance-free domains}

As an immediate consequence of frequency localizations, we obtain the following result on eigenvalue- and resonance-free domains away from thresholds.

\begin{theorem} \label{th2.5} Assume condition (\ref{ass1}). Let $\theta \in \bC$ with $\Im \theta = \varphi \in ]0, a[$. Then, there exists $\gamma_0>0$
such that for $\gamma \ge \gamma_0$ the resolvent $R(z, \gamma, \theta)$ is holomorphic in
$\Omega_1(\gamma)$ and, for every $k\ge 1$ and $\gamma\ge\gamma_0(k)$, in the domain $\Omega_{1,k}(\gamma)$ of \eqref{Omega1k}.
\end{theorem}
\pf Let $\Upsilon = \Upsilon^{(0)}$ and $Y= - \Upsilon\cdot (\nabla V)(e^{\theta}x)$. Then $Y^2=0$.
By Corollary \ref{cor2.4}(a) with $k=0$,
\be \label{e4.17}
1+ W_\theta R_0  = 1 + \vY(z,\gamma) + O(\gamma^{-\f 1 2})
\ee
for $z\in \Omega_1(\gamma)$, where
\be\label{e4.18}
\vY(z,\gamma)  = Y  \chi_0(\gamma^{-1} D_x) (-e^{-2\theta}\gamma^{-1} \Delta_x-z)^{-1}.
\ee
Since $\vY(z, \gamma) ^2 =0$,  $1+ \vY(z, \gamma) $ is invertible with the inverse given by
\[
(1+  \vY(z, \gamma)  )^{-1} = 1-  \vY(z, \gamma) .
\]
Since $(-e^{-2\theta} \gamma^{-1} \Delta -z)^{-1}$ is uniformly bounded for $z\in \Omega_1(\gamma)$, there exists $\gamma_0$ such that
$1 + W_\theta R_0$ is invertible for $z\in \Omega_1(\gamma)$ and $\gamma \ge \gamma_0$ with the inverse given by
\[
(1 + W_\theta R_0 )^{-1} = 1 - \vY(z, \gamma) +  O(\gamma^{-\f 1 2}).
\]
Consequently,
\[
R =  R_0 (1+ W_\theta R_0)^{-1}
\]
is holomorphic for $z\in \Omega_1(\gamma)$.
The proof for $k \ge 1$ is the same. Corollary \ref{cor2.4}(a) is used with $\Upsilon^{(k)}$, which is nilpotent, and
$(-e^{-2\theta}\gamma^{-1}\Delta_x+\gamma k-z)^{-1}$ is uniformly bounded for $z\in\Omega_{1,k}(\gamma)$.
\ef

In particular, for $\gamma$ large, $P(\gamma)$ has neither eigenvalues nor resonances in $\Omega_1(\gamma)$, nor, for $k\ge1$, in $\Omega_{1,k}(\gamma)$. The rest of this work is devoted to the eigenvalues and resonances near the thresholds $\vT_\gamma$.

\sect{A class of Witten Laplacians}\label{sect3}

To study the eigenvalues and resonances near the threshold $\gamma k$, $k \in \bN$, we use a
Grushin problem \cite{sz} on the first $k+2$ eigenvalues of $H_0$. The number $k+2$ is the smallest choice for
which the level $k$ is separated from the part of the spectrum which is eliminated.
The blocks of the Grushin problem produce the Witten Laplacian at the level $k$, $-\Delta^{(k)}_V(\theta)$, introduced in Subsection \ref{ssec-grushin}. For $k=0$, it is the dilated Witten Laplacian $-\Delta_V(\theta)$. In this section, $k\in\bN$ is fixed.

\subsection{The Grushin problem}\label{ssec-high}

Let $\Pi_j$ denote the orthogonal projection associated with the eigenvalue $j$ of $H_0$. For $J\subset\bN$ put $\Pi_J = \sum_{j\in J}\Pi_j$,
\[
P^J(\gamma,\theta) = \Pi_JP(\gamma,\theta)\Pi_J, \qquad P^J_0(\gamma,\theta) = \Pi_JP_0(\gamma,\theta)\Pi_J.
\]
Set
$\Pi = \sum_{j=0}^{k+1} \Pi_j$, $\Pi' = 1-\Pi = \Pi_{J'}$, $P' = P^{J'}$, $P'_0 = P^{J'}_0$, with $J' = \{k+2,k+3,\dots\}$,
and, when these inverses exist,
\be \label{Rprimedef}
R^J_0(z,\gamma,\theta) = (P^J_0(\gamma,\theta)-z)^{-1}\Pi_J, \quad R^J(z,\gamma,\theta) = (P^J(\gamma,\theta)-z)^{-1}\Pi_J,
\ee
$R'_0 = R^{J'}_0$ and $R' = R^{J'}$. As stated in the paragraph \emph{Notation}, the arguments
$z,\gamma,\theta$ are written in the statements and are often omitted in the proofs.

\begin{figure}[ht]
\centering
\includegraphics[width=0.4\textwidth]{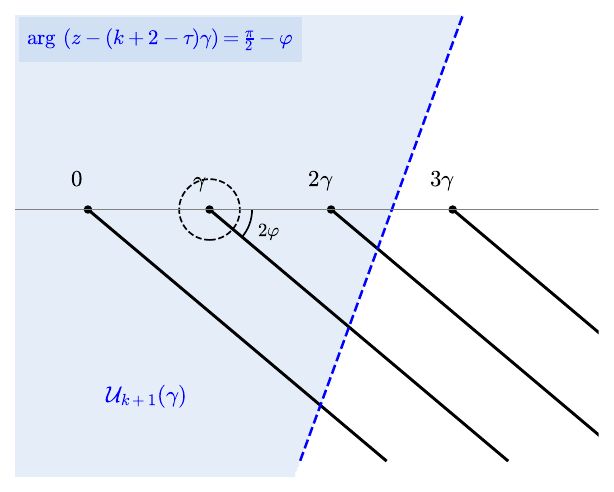}
\caption{\small The spectrum of $P_0(\gamma,\theta)$, $\Im\theta = \varphi>0$ (half-lines), and the half-plane $\vU_{k+1}(\gamma)$ of \eqref{sector}, drawn for $k = 1$; the
dashed circle is $\{|z-\gamma k|\le c_1\}$.}\label{fig-U}
\end{figure}

We consider sets $J$ lying on one side of the level $k$, that is, $J\subset\{m,m+1,\dots\}$ with
$m\ge k+1$, or $J\subset\{0,\dots,M\}$ with $M\le k-1$. For $\theta\in\bC$ with $|\varphi|<a \le \f \pi 4$ and a
fixed $\tau\in\,]0,1[$, let $\vU_J(\gamma)$ be the half-plane
\be \label{sector}
\vU_J(\gamma) = \left\{ \begin{array}{ll}
\left\{z\in\bC;\ \Re(e^{i\varphi}z) < (m-\tau)\gamma\cos\varphi\right\}, & J\subset\{m,m+1,\dots\},\\[2mm]
\left\{z\in\bC;\ \Re(e^{i\varphi}z) > (M+\tau)\gamma\cos\varphi\right\}, & J\subset\{0,\dots,M\},
\end{array}\right.
\ee and write $\vU_{k+1}(\gamma) := \vU_{J'}(\gamma)$. Each $\vU_J(\gamma)$ contains the
disc $\{|z-\gamma k|\le c_1\}$ for any fixed $c_1$ and $\gamma$ large. Figure \ref{fig-U} shows
\[
\sigma(P_0(\gamma,\theta)) = \bigcup_{j\in\bN}\big(\gamma j+e^{-2i\varphi}[0,\infty[\big)
\]
 and the half-plane $\vU_{k+1}(\gamma)$.

Let $f\in\Ran\Pi_J\cap\vS$. For $\theta = i \varphi$ with $\varphi$ real, one has
\begin{align}
|\w{f, ( P^J_0(\gamma, \theta)-z)f}| &\ge  \left|\Re \w{f, e^{i\varphi} ( P^J_0(\gamma, \theta)-z)f}\right| \nonumber\\
& = \left|\w{f, (\gamma  \cos \varphi \, H_0 - \Re (e^{i\varphi} z)) f}\right|, \label{e3.1}
\end{align}
and $H_0\ge m$, resp. $H_0\le M$, on $\Ran\Pi_J$. Hence, for $z\in\vU_J(\gamma)$,
$\|(P^J_0(\gamma,\theta)-z)f\|\ge\tau\gamma\cos\varphi\,\|f\|$, and the same holds for the adjoint.
Therefore $R^J_0(z,\gamma,\theta)$ exists on $\Ran\Pi_J$. The same
holds for every $\theta\in\bC$ with imaginary part $\varphi$, since $P^J_0(\gamma,\theta)$ and
$P^J_0(\gamma,i\varphi)$ are unitarily equivalent. Moreover,
\be \label{e4.19}
\| R^J_0(z, \gamma, \theta)\| \le \f{1}{\tau\gamma\cos\varphi} , \qquad z\in\vU_J(\gamma),\ \gamma\ge1 .
\ee
If $J\subset\{0,\dots,M\}$, $H_0$ is bounded on $\Ran\Pi_J$. If $J\subset\{m,m+1,\dots\}$, applying
\eqref{e3.1} to $f = R^J_0g$ gives
$$\gamma\cos\varphi\,\w{f,H_0f} \le \|f\|\,\|g\| + \Re(e^{i\varphi}z)\|f\|^2.$$
If
$\Re(e^{i\varphi}z)\le0$, this gives the bound below directly. Otherwise,
$\Re(e^{i\varphi}z)\le m\gamma\cos\varphi$, and \eqref{e4.19} gives $\w{f,H_0f}\le C\gamma^{-2}\|g\|^2$. In
both cases,
\be \label{e4.19b}
\| (1+H_0)^{\f 12} R^J_0(z, \gamma, \theta)\| \le C\gamma^{-1}, \qquad z\in\vU_J(\gamma),\ \gamma\ge1 .
\ee

\begin{lemma}\label{lem3.1} Assume \eqref{ass1} and let $J$ be as above. There exists $\gamma_0\ge 1$
such that for $z \in \vU_J(\gamma)$ and $\gamma \ge \gamma_0$, $P^J(\gamma, \theta)-z$ is invertible on
$\Ran \Pi_J$, and $R^J(z, \gamma, \theta)$ satisfies
\be \label{e4.20}
\| R^J(z, \gamma, \theta)\| + \|(1+H_0)^{\f 12}  R^J(z, \gamma, \theta)\| \le C \gamma^{-1}
\ee
uniformly in $z \in \vU_J(\gamma)$ and $\gamma \ge \gamma_0$. In particular, $R'(z,\gamma,\theta)$ satisfies \eqref{e4.20} for $z\in\vU_{k+1}(\gamma)$. Moreover, $R^J(z,\gamma,\theta)$ is holomorphic in $(z,\theta)$, for $|\theta|<a$ and $z\in\vU_J(\gamma)$.
\end{lemma}
\pf With $W^J_\theta = \Pi_JW_\theta\Pi_J$, $W_\theta = -(\nabla V)(e^\theta x)\cdot \nabla_v$, one has

\[
P^J-z = (1+W^J_\theta R^J_0)(P^J_0-z)
\]
 on $\Ran\Pi_J$, and,
by \eqref{e4.19b} and \eqref{ass1},
\[
\|W^J_\theta R^J_0\| \le C\|\na_v R^J_0\| \le C'\gamma^{-1} \le \f12
\]
for $\gamma$ large. Hence $R^J = R^J_0(1+W^J_\theta R^J_0)^{-1}$, and \eqref{e4.20} follows from
\eqref{e4.19} and \eqref{e4.19b}. The holomorphy in $z$ is clear. The holomorphy in $\theta$ follows from the fact that $\{P(\gamma,\theta)\}_{|\theta|<a}$ is a holomorphic family of type (A), whose domain is independent of $\theta$ by \eqref{sub}. Note that the projections $\Pi_J$ do not depend on $\theta$.
\ef

The following result is a subelliptic estimate for $P'(\gamma, \theta)$, analogous to Theorem \ref{th2.1}, but without the term $\|v\cdot \nabla_x f\|$.

\begin{prop}\label{prop-subRJ} Assume condition \eqref{ass1}. There is $\gamma_0$ such that for every $M>0$ there is $C_M$ with
\be \label{subRJ}
\|\gamma H_0 R'(z, \gamma, \theta)\| + \gamma^{\f 1 3}\| (\gamma +\w{D_x})^{\f 2 3} R'(z, \gamma, \theta)\| \le C_M
\ee
for $z \in \vU_{k+1}(\gamma)$ with $|z|\le M\gamma$, and $\gamma \ge \gamma_0$.
\end{prop}
\pf The method is the commutator argument which gives the global subelliptic estimates of the KFP operator, see \cite{hn,hln} and \cite[Theorem 2.1]{w4}. We adapt it to the present situation. Recall that $R'=R^J$ with $J=J'=\{m,m+1,\dots\}$, $m=k+2$. Consider first $V=0$. A dilation in $x$ shows that $P^J_0$ is unitarily equivalent to $P^J_0(\gamma,i\varphi)$ and changes $\w{D_x}$ by a bounded factor. Hence we may assume $\theta=i\varphi$, $|\varphi|<\f\pi4$. After the Fourier transform in $x$ and the scaling $\xi\mapsto\gamma^{-1}\xi$, $z\mapsto\gamma^{-1}z$, which factors out $\gamma$ as in the proof of Theorem \ref{th2.1}, we have to estimate on $\Ran\Pi_J$ the operator $H_0+\kappa B-z$, where $\kappa=ie^{-i\varphi}$ and $B=\Pi_J(v\cdot\xi)\Pi_J=b+b^*$, with $b=\Pi_Ja_\xi\Pi_J$ and $a_\xi=\xi\cdot a$.
The operator $a_\xi^*$ raises the level and $a_\xi$ lowers it. Hence $b^*=a^*_\xi$ on $\Ran\Pi_J$, while $b=a_\xi$ on $\Ran\Pi_{m+1}\oplus\Ran\Pi_{m+2}\oplus\cdots$ and $b=0$ on $\Ran\Pi_m$. For $f\in\Ran\Pi_J\cap\vS$, one has
\be \label{bb}
\w{f,[b,b^*]f} = |\xi|^2\|f\|^2+\|a_\xi\Pi_mf\|^2\ \ge\ |\xi|^2\|f\|^2, \qquad \|bf\|+\|b^*f\|\le 2|\xi|\,\|(1+H_0)^{\f12}f\|.
\ee
Let $s\ge 1$, $D=H_0+s$, $g=(D+\kappa B)f$, $S=i[H_0,B]=i(b^*-b)$, $G=\|g\|$ and $X=\|f\|$. The real and imaginary parts of $\w{g,f}$ give $\w{Df,f}\le (1+|\tan\varphi|)GX\le 2GX$, hence $\|(1+H_0)^{\f12}f\|^2\le 2GX$. Next, $[D,B]=-iS$ and $|\kappa|=1$ give
\[
G^2 = \|Df\|^2+\|Bf\|^2+2\sin\varphi\,\Re\w{Df,Bf}+\cos\varphi\,\w{f,Sf},
\]
and, since $|\sin\varphi|<2^{-\f12}$ and $\|Sf\|\le 2|\xi|\,\|(1+H_0)^{\f12}f\|$,
\be \label{sub1}
\|Df\|^2+\|Bf\|^2\le C\left(G^2+|\xi|\,G^{\f12}X^{\f32}\right).
\ee
Finally, $[B,S]=2i[b,b^*]$ shows that $\Im\w{Bf,Sf}=\w{f,[b,b^*]f}$, and the real part of $\w{g,Sf}=\w{Df,Sf}+\bar\kappa\w{Bf,Sf}$ gives, by \eqref{bb} and \eqref{sub1},
\[
|\xi|^2X^2\le C\left(G+\|Df\|+\|Bf\|\right)|\xi|\,(GX)^{\f12}
\le C\left(|\xi|\,G^{\f32}X^{\f12}+|\xi|^{\f32}G^{\f34}X^{\f54}\right).
\]
For $|\xi|\ge1$, writing $X=\lambda G|\xi|^{-\f23}$ one obtains $\lambda^2\le C(\lambda^{\f12}+\lambda^{\f54})$, hence $X\le CG|\xi|^{-\f23}$. For $|\xi|\le1$, $X\le 2G$ because $\w{Df,f}\ge X^2$. Then \eqref{sub1} gives $\|Df\|+\|Bf\|\le CG$, and we have proved
\be \label{subJs}
\|(H_0+s)f\| + \|Bf\| + \w{\xi}^{\f23}\|f\| \le C\|(H_0+\kappa B+s)f\|, \qquad s\ge 1 .
\ee
Now let $z\in\vU_J(\gamma)$, $|z|\le M\gamma$, and $g=(H_0+\kappa B-\gamma^{-1}z)f$. If $|\xi|\ge R$ with $R$ large, we write
$$H_0+\kappa B-\gamma^{-1}z=(H_0+\kappa B+1)-(\gamma^{-1}z+1).$$ Since $|\gamma^{-1}z+1|\le M+1$ is absorbed by $\w{\xi}^{\f23}\|f\|$ in \eqref{subJs}, the three terms of \eqref{subJs} are bounded by $C\|g\|$. If $|\xi|\le R$, then $\|f\|+\|(1+H_0)^{\f12}f\|\le C\|g\|$ by \eqref{e4.19} and \eqref{e4.19b}, hence
$$\|Bf\|\le R\|(1+H_0)^{\f12}f\|\le C\|g\| \qquad \mbox{and} \qquad \|H_0f\|\le\|g\|+M\|f\|+\|Bf\|\le C\|g\|.$$ This proves \eqref{subRJ} for $V=0$.

In general, $R^J=R^J_0(1+W^J_\theta R^J_0)^{-1}$ with $\|W^J_\theta R^J_0\|\le C\gamma^{-1}\le\f12$ by the proof of Lemma \ref{lem3.1}. The operators $\gamma H_0$ and $\gamma^{\f13}(\gamma+\w{D_x})^{\f23}$ are placed on the left of $R^J_0$, so \eqref{subRJ} follows from the case $V=0$.
\ef

\begin{lemma} \label{lem3.4}
Assume condition \eqref{ass1b}, let $J$ be as in Lemma \ref{lem3.1}, $R^J =R^J(z,\gamma, \theta)$ and $r \in \bR$. Then
\be \label{elm4.7}
\| \w{x}^{r} R^J\w{x}^{-r}\| + \| \w{D_x}^{r} R^J\w{D_x}^{-r}\|
 \le C \gamma^{-1}
\ee
uniformly for $z\in \vU_J(\gamma)$ and $\gamma \ge \gamma_0$. Moreover, for $M>0$,
\be \label{elm4.7b}
\|\w{v}^{r} R'\w{v}^{-r}\| + \| \w{D_v}^{r} R'\w{D_v}^{-r}\|
\le C \gamma^{-1}
\ee
uniformly for $z\in \vU_{k+1}(\gamma)$ with $|z|\le M\gamma$, and $\gamma \ge \gamma_0$. In particular
\be \label{RS}
R'(z,\gamma,\theta)\,\vS(\bR^{2n}) \subset \vS(\bR^{2n}), \qquad z\in\vU_{k+1}(\gamma),\ |z|\le M\gamma, \; \gamma\ge\gamma_0 .
\ee
\end{lemma}
\pf Estimate (\ref{elm4.7}) for $\w{x}^{r}$ follows from
\[
\w{x}^{r} R^J\w{x}^{-r} = R^J + R^J O(v \w{x}^{r-1})R^J \w{x}^{-r} = O(\gamma^{-1}) + O(\gamma^{-1})\, \w{x}^{r} R^J\w{x}^{-r}
\]
by Lemma \ref{lem3.1}, for $\gamma$ large. The same argument applies to $\w{D_x}^r$.
For $\w{v}^r$ in \eqref{elm4.7b}, $[\Pi', \w{v}^r](1+ H_0)^\ell$ is bounded  for any $\ell \ge 0$, and
\[
[P',\w{v}^r] = O(\gamma\w{v}^{r-1}\na_v) + O(\gamma\w{v}^{r-2}) + O(\w{v}^{r-1}),
\]
 so that $R'[P',\w{v}^r]\w{v}^{-r+1} =O(1)$ by Lemma \ref{lem3.1}.
This gives \eqref{elm4.7b} for $\w{v}^r$ with $0\le r \le 1$. The general case $r\in \bR$ follows from successive commutators and duality.
For $\w{D_v}^r$ in (\ref{elm4.7b}), the commutator of $R'$ with $D_v$ gives rise to $\p_x$, which can be estimated as
\[
[D_v, R'] = R' (O(\p_x) + O(\gamma v)) R'= O(\gamma^{-1}).
\]
Indeed, the $O(\gamma v)$ term can be estimated by Lemma \ref{lem3.1}. Moreover, Proposition \ref{prop-subRJ} gives $\|\w{D_x}^{\f23}R'\|\le C\gamma^{-\f13}$, hence $\|\w{D_x}^{\f13}R'\|\le C\gamma^{-\f23}$ by \eqref{e4.20} and interpolation.
This proves (\ref{elm4.7b}) for $r=1$. The general case $r\in \bR$ follows from successive commutators and complex interpolation. The details are omitted. \ef

Let $\mu_j = \mbox{\rm rank}\,\Pi_j$ be the multiplicity of the eigenvalue $j$ of $H_0$ and set
\be \label{Gj}
\vG_j := L^2(\bR^n_x, \bC^{\mu_j}), \qquad \vG := \bigoplus_{j=0}^{k+1}\vG_j
\simeq L^2(\bR^n_x,\bC^\mu), \qquad \mu := \sum_{j=0}^{k+1} \mu_j .
\ee
By \eqref{ident}, the components of an element of $\vG_j$ are numbered by the multi-indices
$\alpha \in \bN^n$ with $|\alpha| = j$. Define, for $j\in\bN$,
\be \label{Epm}
\begin{aligned}
&E^{(j)}_+ : \vG_j \to L^2(\bR^{2n}), \quad u \mapsto \sum_{|\alpha|=j} u_\alpha \otimes \psi_\alpha, \\
&E^{(j)}_- : L^2(\bR^{2n}) \to \vG_j, \quad f \mapsto \left( \w{\psi_\alpha, f}_v\right)_{|\alpha|=j},
\end{aligned}
\ee
where $\w{\cdot,\cdot}_v$ is the scalar product in the $v$-variables, and let
$E_\pm = \bigoplus_{j=0}^{k+1} E^{(j)}_\pm$. Then
\be \label{Erules}
\begin{aligned}
&E^{(i)}_-E^{(j)}_+ = \delta_{ij}1_{\mu_j}, \qquad E^{(j)}_+E^{(j)}_- = \Pi_j, \\
&\left(E^{(j)}_\pm\right)^* = E^{(j)}_\mp, \qquad E_+E_- = \Pi, \qquad E_-E_+=1 .
\end{aligned}
\ee
In the rest of this section, $\gamma$ and $\theta$ are fixed, and we omit them from the notation. For example, we write $P = P(\gamma,\theta)$, $R(z) = R(z,\gamma,\theta)$, $R'(z) = R'(z,\gamma,\theta)$ and $\vE_{-+}(z) = \vE_{-+}(z,\gamma,\theta)$. Consider the Grushin problem
\be \label{Gr1}
\vP(z) = \left(\begin{array}{cc}
P-z & E_+ \\
E_- & 0
\end{array} \right): D(P)\times  \vG  \to L^2(\bR_{x,v}^{2n})\times  \vG .
\ee

This Grushin problem is \emph{not} well
posed in the sense of \cite{sz}, because
$\vP(z)$, as an unbounded operator on $L^2(\bR_{x,v}^{2n})\times  \vG $, is closed and
injective, but not surjective (see Proposition \ref{prop-Gr} in Appendix \ref{app-Gr}).

By Lemma \ref{lem3.1}, $R'(z)$ exists for $z \in \vU_{k+1}(\gamma)$ and
$\gamma\ge\gamma_0$. Under \eqref{ass1b}, it maps $\vS(\bR^{2n})$ into $\vS(\bR^{2n})$ by \eqref{RS}.
The following computations are valid on
$\vS(\bR^{2n}) \times \vS(\bR^n;\bC^\mu)$. One has:
\[
\vP(z) \left(\begin{array}{cc}
R'(z) & E_+ \\
E_- &  E_- (z-P )E_+
\end{array} \right) =  \left(\begin{array}{cc}
1 + \Pi L_\theta R'(z) & \Pi'L_\theta E_+ \\
0 &  1
\end{array} \right).
\]
Since $R'\Pi = 0$, the operator $\vK := \Pi L_\theta R'(z)$ satisfies $\vK^2=0$, so
that $(1 + \vK)^{-1} = 1 - \vK$ and
\[
\left(\begin{array}{cc}
1 + \vK &  \Pi'L_\theta E_+ \\
0 &  1
\end{array} \right)^{-1}
= \left(\begin{array}{cc}
1-\vK    &  - (1-\vK) \Pi'L_\theta E_+ \\
0 &  1
\end{array} \right).
\]
The same computation on the other side gives a left inverse. Thus, on
$\vS(\bR^{2n}) \times \vS(\bR^n;\bC^\mu)$,
\[
\vP(z)^{-1}
= \left(\begin{array}{cc}
R'(z) & E_+ \\
E_- &  E_- (z-P )E_+
\end{array} \right) \left(\begin{array}{cc}
1-\vK   &  -(1-\vK) \Pi'L_\theta E_+ \\
0 &  1
\end{array} \right).
\]
In particular $\vP(z)$ is a \emph{bijection} of
$\vS(\bR^{2n})\times\vS(\bR^n;\bC^\mu)$ onto itself. The
Grushin problem is thus exactly solvable in $\vS$, where \eqref{r1}, \eqref{r2} and \eqref{r4} below are
identities. What fails is only the continuity of $\vP(z)^{-1}$ for the $L^2$-norms.
This is the content of Proposition \ref{prop-Gr}.

We write this operator in the form
\be
\vP(z)^{-1} = \left(\begin{array}{cc}
\vE(z) & \vE_+(z) \\
\vE_-(z) &  \vE_{-+}(z)
\end{array} \right)
\ee
where
\beas
\vE(z)  & = & R'(z) \\
\vE_+(z) & = & (1 - R'(z)L_\theta) E_+ \\
\vE_-(z) & = & E_- (1- L_\theta  R'(z)) \\
\vE_{-+}(z) & = & E_- (z-P )E_+ +  E_- L_\theta  R'(z)L_\theta E_+.
\eeas
The identities $\vP \vP^{-1} =  \vP^{-1} \vP= 1$ give, on $\vS$,
\bea
(P-z)\vE(z) + E_+ \vE_-(z) & =& 1, \label{r1} \\
(P-z) \vE_+(z) + E_+\vE_{-+}(z) & =&0, \label{r2} \\
\vE_-(z) (P-z) + \vE_{-+}(z) E_- &=& 0. \label{r4}
\eea
Let $\bJ$ be the parity operator
\be \label{bJ}
(\bJ u)(x,v) = u(x,-v).
\ee
Since
$\psi_\alpha(-v) = (-1)^{|\alpha|}\psi_\alpha(v)$, one has $\bJ = \bJ^*=\bJ^{-1}$,
$\bJ\Pi_j = \Pi_j\bJ$, $\bJ|_{\Ran\Pi_j} = (-1)^j$, and $\bJ$ preserves $\vS(\bR^{2n})$. As $H_0$
is even in $v$ and $L_\theta$ is odd in $v$,
\be \label{tRprime}
{}^tP = \bJ P\bJ, \qquad
{}^tR'(z) = \bJ R'(z)\bJ .
\ee
Each factor of $\vP(z)^{-1}$ therefore extends by transposition to a continuous operator on $\vS'(\bR^{2n})$. By duality, the above formulas also hold in $\vS'$.

\subsection{The blocks of $\vE_{-+}$}\label{ssec-blocks}
The operator $\vE_{-+}(z)$ acts on $\vG = \bigoplus_{j=0}^{k+1}\vG_j$. We denote by
$\vE_{-+}^{(ij)} : \vG_j \to \vG_i$ its blocks with respect to this decomposition. They are
computed from the following lemma.
Set, for $1 \le l \le n$,
\be \label{Xl}
X^\theta_l = e^{-\theta}\p_{x_l} + \f 12 (\p_{x_l}V)(e^\theta x), \qquad
{}^tX^\theta_l = -e^{-\theta}\p_{x_l} + \f 12 (\p_{x_l}V)(e^\theta x).
\ee
Here ${}^t X^\theta_l$ is the transpose of $X^\theta_l$ with respect to the bilinear form
$\int f g\,dx$. With $a_l$, $a^*_l$ as in \eqref{creation}, set
\be \label{AAdag}
A := \sum_{l=1}^n a_l\,{}^tX^\theta_l , \qquad A^\dag := \sum_{l=1}^n a^*_l X^\theta_l .
\ee
Thus $A^\dag$ raises the Hermite level by one and $A$ lowers it by one. Set finally, for
$j \in \bN$,
\be \label{Xblocks}
\bX^\theta_j := E^{(j+1)}_- L_\theta E^{(j)}_+ \, : \, \vG_j \to \vG_{j+1}, \qquad
{}^t\bX^\theta_j := - E^{(j)}_-L_\theta E^{(j+1)}_+ \, :\, \vG_{j+1}\to \vG_j .
\ee

\begin{lemma}\label{lem-blocks} One has
\be \label{LAB}
L_\theta = A^\dag - A .
\ee
Consequently $\w{\psi_\beta, L_\theta\psi_\alpha}_v = 0$ unless $\big| |\beta|-|\alpha|\big| = 1$,
\be \label{entries}
\w{\psi_{\alpha+e_l}, L_\theta \psi_\alpha}_v = \sqrt{\alpha_l+1}\; X^\theta_l, \qquad
\w{\psi_{\alpha}, L_\theta \psi_{\alpha+e_l}}_v = - \sqrt{\alpha_l+1}\; {}^tX^\theta_l ,
\ee
and $\bX^\theta_j$ is the $\mu_{j+1}\times \mu_j$ matrix with entries
$(\bX^\theta_j)_{\beta\alpha} = \sqrt{\alpha_l+1}\,X^\theta_l$ if $\beta = \alpha + e_l$ and $0$
otherwise. Moreover, ${}^t\bX^\theta_j$ is its transpose, entry by entry.
\end{lemma}
\pf By \eqref{aa} and the expression of $L_\theta$ given in Section \ref{sect2}, with
$g_l = (\p_{x_l}V)(e^\theta x)$,
\beas
L_\theta &=& \sum_l \left[ e^{-\theta}(a_l+a^*_l)\p_{x_l} - \f 12 g_l (a_l - a^*_l)\right] \\
&=& \sum_l\left[ a_l \left( e^{-\theta}\p_{x_l} - \f 12 g_l\right)
+ a^*_l\left( e^{-\theta}\p_{x_l}+\f 12 g_l\right)\right]
= \sum_l \left[ -a_l\, {}^tX^\theta_l + a^*_l\, X^\theta_l\right],
\eeas
which is \eqref{LAB}. Then \eqref{entries} follows from \eqref{aaction}, and the matrix form of
$\bX^\theta_j$ from \eqref{Xblocks} and \eqref{entries}.
\ef

\begin{rmk}\label{rmk-LR} The operator $L_\theta R'(z)$ is not bounded on
$L^2(\bR^{2n})$. Indeed,
\[
\Pi L_\theta R'(z) = \Pi_{k+1}L_\theta\Pi_{k+2}R'(z)
= -E^{(k+1)}_+\,{}^t\bX^\theta_{k+1}E^{(k+2)}_-R'(z), \]
where ${}^t\bX^\theta_{k+1}$ is a
first-order differential operator in $x$, while $R'$ gains only $\f 23$ derivative in $x$ (compare
Theorem \ref{th2.1}).
This is the source of the ill-posedness of the Grushin problem and of the domain problems below.
\end{rmk}

For $n=1$ one has $\mu_j = 1$, $\bX^\theta_j = \sqrt{j+1}\,X^\theta$, and \eqref{entries} reduces to
the recurrence relations of the Hermite functions. For $j=0$, $\mu_0 = 1$, $\mu_1 = n$ and
$\bX^\theta_0 = {}^t(X^\theta_1,\dots,X^\theta_n)$,
${}^t\bX^\theta_0 = ({}^tX^\theta_1,\dots,{}^tX^\theta_n)$. The Witten Laplacian
\be \label{witten}
-\Delta_V(\theta) = \sum_{l=1}^n {}^tX^\theta_lX^\theta_l = {}^t\bX^\theta_0\bX^\theta_0
\ee
thus arises from the interaction between the drift operator $L_\theta$ and the diffusion operator
$H_0$. This is the mechanism already observed in \cite{hln}.

\begin{prop}\label{prop-tri} With respect to $\vG = \bigoplus_{j=0}^{k+1} \vG_j$, the operator
$\vE_{-+}(z)$ is block tridiagonal:
\be \label{Tentries}
\vE_{-+}^{(jj)} = (z - \gamma j) 1_{\mu_j} - \delta_{j,k+1}N_{k+1}, \qquad
\vE_{-+}^{(j+1,j)} = -\bX^\theta_j, \qquad \vE_{-+}^{(j,j+1)} = {}^t\bX^\theta_j ,
\ee
all the other blocks being zero. Here
\be \label{N}
N_{k+1} = N_{k+1}(z) = {}^t\bX^\theta_{k+1}\, R'_{k+2}\, \bX^\theta_{k+1}, \qquad
R'_{k+2} := E^{(k+2)}_- R'(z) E^{(k+2)}_+ .
\ee
Equivalently, writing $u_j = E^{(j)}_-u$ and $\phi_j = E^{(j)}_-\phi$, the equation
$-\vE_{-+}(z)E_-u = \phi$ is the block tridiagonal system
\be \label{system}
\bX^\theta_{j-1}u_{j-1} + (\gamma j - z)u_j - {}^t\bX^\theta_j u_{j+1} = \phi_j,
\qquad 0 \le j \le k+1, \quad (u_{-1}:=0),
\ee
where, at the top level $j=k+1$, the term $-{}^t\bX^\theta_{k+1}u_{k+2}$ has been replaced by $N_{k+1}u_{k+1}$.
\end{prop}
\pf We compute directly on the equation $(P-z)u = f$. This fixes every sign. Put
$u_j = E^{(j)}_-u$, $f_j = E^{(j)}_-f$ for $j \in \bN$. Since $\gamma H_0 - z$ acts on
$\Ran \Pi_j$ as the scalar $\gamma j - z$ and, by Lemma \ref{lem-blocks},
$E^{(j)}_-L_\theta = \bX^\theta_{j-1}E^{(j-1)}_- - {}^t\bX^\theta_j E^{(j+1)}_-$, the equation is
equivalent to the infinite block tridiagonal system in $L^2(\bR^n_x,\bC^{\mu_j})$:
\be \label{systeminf}
\bX^\theta_{j-1}u_{j-1} + (\gamma j - z)u_j - {}^t\bX^\theta_j u_{j+1} = f_j, \qquad j \in \bN,
\quad (u_{-1}:=0).
\ee
Eliminate the levels $\ge k+2$: applying $\Pi'$ to $(P-z)u=f$ gives
\[
u' = R'f - R'E^{(k+2)}_+\bX^\theta_{k+1}u_{k+1}, \qquad\mbox{hence}\qquad
u_{k+2} = E^{(k+2)}_-R'f - R'_{k+2}\bX^\theta_{k+1} u_{k+1},
\]
and the equation \eqref{systeminf} with
$j=k+1$ becomes
\[
\bX^\theta_{k}u_{k} + \left[ (\gamma(k+1) - z) + N_{k+1}\right] u_{k+1}
= f_{k+1} + {}^t\bX^\theta_{k+1}E^{(k+2)}_-R'f ,
\]
the equations $j \le k$ being unchanged. The left-hand sides are the rows of \eqref{system}, and
the right-hand side is the $(k+1)$-th component of $\vE_-f$ by the formula for $\vE_-$ above. Thus
$-\vE_{-+}E_-u = \vE_-f$, which is \eqref{r4}, and the blocks are given by \eqref{Tentries}.
\ef

Therefore the matrix $-\vE_{-+}$ is the block tridiagonal matrix with diagonal blocks
\[
(\gamma j -z)1_{\mu_j} + \delta_{j,k+1}N_{k+1} ,
\]
subdiagonal blocks $\bX^\theta_j$ and superdiagonal blocks
$-{}^t\bX^\theta_j$. Explicitly, in the decomposition
$\vG = \bigoplus_{j=0}^{k+1}\vG_j$, with the row index $j$ (resp. column index $j'$)
running from $0$ to $k+1$,
\be \label{tableE}
-\vE_{-+}(z) =
\begin{pmatrix}
-z & -{}^t\bX^\theta_0 & 0 & \cdots & 0 \\[2pt]
\bX^\theta_0 & \gamma-z & -{}^t\bX^\theta_1 & \ddots & \vdots \\[2pt]
0 & \bX^\theta_1 & 2\gamma-z & \ddots & 0 \\[2pt]
\vdots & \ddots & \ddots & \ddots & -{}^t\bX^\theta_k \\[2pt]
0 & \cdots & 0 & \bX^\theta_k & \gamma(k+1)-z+N_{k+1}
\end{pmatrix},
\ee
that is, the entry $(j,j')$ is $\bX^\theta_{j'}$ if $j=j'+1$, $-{}^t\bX^\theta_{j}$ if $j'=j+1$,
$(\gamma j-z)1_{\mu_j}+\delta_{j,k+1}N_{k+1}$ if $j=j'$, and $0$ otherwise. This is
\eqref{system} read row by row.

Let $h^\theta = (\mbox{\rm Hess}\,V)(e^\theta x)$, a complex symmetric matrix with entries
\[
h^\theta_{lm}(x) = (\p_{x_l}\p_{x_m}V)(e^\theta x) ,
\]
and let
\be \label{dGdef}
\dG(M) := \sum_{l,m=1}^n M_{lm}\,a^*_la_m
\ee
be the second quantization of an $n\times n$ complex matrix. $\dG(h^\theta) $ acts on $\bigoplus_j L^2(\bR^n_x;\allowbreak\Ran\Pi_j)$ and leaves each
summand stable.

\begin{lemma}[the key commutator]\label{lem-comm} One has
\be \label{comm}
[X^\theta_l, X^\theta_m] = [{}^tX^\theta_l, {}^tX^\theta_m] = 0, \qquad
[{}^tX^\theta_l, X^\theta_m] = -h^\theta_{lm}, \qquad
\sum_l {}^tX^\theta_lX^\theta_l = -\Delta_V(\theta) ,
\ee
and, with $A$ and $A^\dag$ as in \eqref{AAdag},
\be \label{AB}
[A,A^\dag] = -\Delta_V(\theta) - \dG(h^\theta) .
\ee
\end{lemma}
\pf With $g_l = (\p_{x_l}V)(e^\theta x)$ one has $\p_{x_l}g_m = e^\theta h^\theta_{lm}$, hence
\[
[{}^tX^\theta_l, X^\theta_m] = \left[-e^{-\theta}\p_l, \f 12 g_m\right] + \left[\f 12 g_l, e^{-\theta}\p_m\right]
= -\f 12 e^{-\theta}(\p_lg_m) - \f 12 e^{-\theta}(\p_mg_l) = -h^\theta_{lm}
\]
by symmetry of the Hessian. The same computation with two $+$ or two $-$ signs gives the first two
identities, and the third is \eqref{witten}. Since the $a$'s commute with the $X^\theta$'s,
\[
[A,A^\dag] = \sum_{l,m}\left( a_la^*_m\,{}^tX^\theta_lX^\theta_m - a^*_ma_l\, X^\theta_m{}^tX^\theta_l\right)
= \sum_{l,m}a^*_ma_l\,[{}^tX^\theta_l,X^\theta_m] + \sum_l {}^tX^\theta_lX^\theta_l ,
\]
where we used $a_la^*_m = a^*_ma_l + \delta_{lm}$. Then \eqref{AB} follows from \eqref{comm} and the
symmetry of $h^\theta$.
\ef

\subsection{The Witten Laplacians at the level \texorpdfstring{$k$}{k}}\label{ssec-grushin}
The blocks of $\vE_{-+}$ at the level $k$ and at its two neighboring levels produce the following operator on $\vG_k$, which plays for the threshold $\gamma k$ the role played by the Witten Laplacian $-\Delta_V$ for the threshold $0$.

\begin{defn}\label{def-Qk} Assume (\ref{ass1}).  For $k \in \bN$ set
\be \label{Qk}
-\Delta^{(k)}_V(\theta) := {}^t\bX^\theta_k \bX^\theta_k - \bX^\theta_{k-1}{}^t\bX^\theta_{k-1}
\qquad \mbox{ on } \ \vG_k = L^2(\bR^n_x,\bC^{\mu_k}),
\ee
with the convention $\bX^\theta_{-1}=0$. We write $-\Delta^{(k)}_V = -\Delta^{(k)}_V(0)$, and call $-\Delta^{(k)}_V$ the \emph{Witten Laplacian at the level $k$}.
\end{defn}
Since $X^\theta_l$ and ${}^tX^\theta_l$ only involve $(\na V)(e^\theta x)$, $-\Delta^{(k)}_V(\theta)$ is well defined
under \eqref{ass1}, as an operator from $H^1$ to $H^{-1}$. The expression \eqref{QkE} of $-\Delta^{(k)}_V(\theta)$
below, with a bounded potential and domain $H^2$, uses in addition the boundedness of the second
derivatives $(\p_{x_l}\p_{x_m}V)(e^\theta x)$, which holds under \eqref{ass1b}.

By \eqref{aaction}, the second quantization \eqref{dGdef} commutes with $H_0$, and we set
\be \label{dGkdef}
\dG_k(M) := E^{(k)}_-\dG(M)E^{(k)}_+ .
\ee
By \eqref{aaction},
\[
a^*_la_m\psi_\alpha = \sqrt{\alpha_m}\sqrt{\alpha_l+1-\delta_{lm}}\,
\psi_{\alpha-e_m+e_l},
\]
so that the matrix of $\dG_k(M)$ in the basis $\{\psi_\alpha\}_{|\alpha|=k}$
is
\[
\left( \dG_k(M)\right)_{\beta\alpha} =
\left\{ \begin{array}{ll}
\sum_l M_{ll}\,\alpha_l , & \beta = \alpha, \\[2mm]
M_{lm}\sqrt{\alpha_m}\sqrt{\alpha_l+1}, & \beta = \alpha - e_m + e_l, \ l \neq m, \\[2mm]
0 & \mbox{ otherwise,}
\end{array}\right.
\]
which is \eqref{dGkintro}.

\begin{lemma}\label{lem-dG} Let $k \in \bN$ and let $M = (M_{lm})$ be an $n\times n$ complex matrix.

\textup{(i)} $\dG_k(M)^* = \dG_k(M^*)$.

\textup{(ii)} If $M$ is real symmetric, with eigenvalues $\nu_1,\dots,\nu_n$, then
\[
\sigma(\dG_k(M)) = \Big\{\sum_l\alpha_l\nu_l;\ |\alpha| = k\Big\} ,
\]
and

\[
k\min_l\nu_l \le \dG_k(M) \le k\max_l\nu_l.
\]

\end{lemma}
\pf (i) follows from $(a^*_la_m)^* = a^*_ma_l$.

(ii) Let $O$ be a real orthogonal matrix such that $OM\,{}^tO = {\rm diag}(\nu_1,\dots,\nu_n)$, and
$(Uf)(v) = f({}^tOv)$. Then $U$ is unitary on $L^2(\bR^n_v)$ and commutes with $H_0$, so that it
leaves $\Ran\Pi_k$ invariant. Moreover, $Ua_lU^{-1} = \sum_pO_{pl}a_p$ and
$Ua^*_lU^{-1} = \sum_pO_{pl}a^*_p$. Hence

\[
U\dG(M)U^{-1} = \dG(OM\,{}^tO) = \sum_p\nu_p\,a^*_pa_p,
\]
 which, by \eqref{aaction}, has the
eigenvector $\psi_\alpha$ with eigenvalue $\sum_p\alpha_p\nu_p$. Since $\sum_p\alpha_p = k$, this
proves (ii).
\ef

\begin{prop}\label{prop-Qk} Assume \eqref{ass1b}. Then
\be \label{QkE}
-\Delta^{(k)}_V(\theta) = E^{(k)}_-[A,A^\dag]E^{(k)}_+ = -\Delta_V(\theta)\otimes 1_{\mu_k} - \dG_k(h^\theta) .
\ee
\end{prop}
\pf By \eqref{Xblocks} and \eqref{Erules},
\[
{}^t\bX^\theta_k\bX^\theta_k = -E^{(k)}_-L_\theta \Pi_{k+1}L_\theta E^{(k)}_+, \qquad
\bX^\theta_{k-1}{}^t\bX^\theta_{k-1} = -E^{(k)}_-L_\theta\Pi_{k-1}L_\theta E^{(k)}_+ .
\]
Since $L_\theta = A^\dag-A$ with $A^\dag$ raising and $A$ lowering the level,
\[
\Pi_{k+1}L_\theta E^{(k)}_+ = \Pi_{k+1}A^\dag E^{(k)}_+, \qquad
E^{(k)}_-L_\theta \Pi_{k+1} = -E^{(k)}_-A\Pi_{k+1},
\]
so that ${}^t\bX^\theta_k\bX^\theta_k = E^{(k)}_-AA^\dag E^{(k)}_+$. In the same way,
$\bX^\theta_{k-1}{}^t\bX^\theta_{k-1} = E^{(k)}_-A^\dag AE^{(k)}_+$. Hence
$-\Delta^{(k)}_V(\theta) = E^{(k)}_-[A,A^\dag]E^{(k)}_+$, and \eqref{QkE} follows from Lemma \ref{lem-comm}.
\ef

By \eqref{QkE} and the expression of $-\Delta_V(\theta)$, we may write
\be \label{Wk}
-\Delta^{(k)}_V(\theta) = -e^{-2\theta}\Delta_x + \vW^\theta_k, \qquad
\vW^\theta_k := \f 14 \sum_{l=1}^n\left( (\p_{x_l}V)(e^\theta \cdot)\right)^2
- \f 12 (\Delta V)(e^\theta \cdot) - \dG_k(h^\theta) .
\ee
Let $\Lambda_\delta$, $\delta\in\,]0,\f\pi2[$, be as in Section \ref{sect1}.
The operator $-e^{-2\theta}\Delta_x$ is normal with spectrum $e^{-2i\varphi}[0,\infty[$. Hence
$\|(-e^{-2\theta}\Delta_x-\eta)^{-1}\| \le (|\eta|\sin\delta)^{-1}$ for $\eta\in\Lambda_\delta$. Moreover, $\vW^\theta_k$ is a bounded matrix-valued function with bounded derivatives by \eqref{ass1b}. A Neumann series and elliptic regularity thus give a constant $C_0$ such that, for $|\theta|<a$ and $\eta \in \Lambda_\delta$ with $|\eta| \ge C_0$,
\be \label{highen}
\left\| (-\Delta^{(k)}_V(\theta) - \eta)^{-1}\right\| \le C_\delta |\eta|^{-1}, \qquad
\left\| \w{D_x}^{s+2}(-\Delta^{(k)}_V(\theta)-\eta)^{-1}\w{D_x}^{-s}\right\| \le C_{\delta, s}, \quad s\in \bR .
\ee

\subsection{The inverse of the Grushin problem}\label{ssec-inv}

The inverse of $\vP(z)$ on $L^2(\bR^{2n})\times\vG$ is more subtle to study, because of domain problems.

By Proposition \ref{prop-Gr} in Appendix \ref{app-Gr}, for $\gamma\ge\gamma_0$ and $z\in\vU_{k+1}(\gamma)$, $\vP(z) : D(P)\times\vG \to
L^2(\bR^{2n})\times\vG$ is closed and injective, and its range is dense but is not the whole space.
Hence $\vP(z)^{-1}$, with domain $\Ran\,\vP(z)$, is a closed, densely defined
and unbounded operator on $L^2(\bR^{2n})\times\vG$, with range $D(P)\times\vG$. We define its entries on their
natural domains: let
\bea
\vD_-(z) &=& \{f\in L^2(\bR^{2n});\ (f,0)\in\Ran\vP(z)\}, \nn \\
\vD_+(z) &=& \{w\in\vG;\ (0,w)\in\Ran\vP(z)\}, \label{Dpmdef}
\eea
and
\bea
\vP(z)^{-1}(f,0) &=& (\vE(z)f,\ \vE_-(z)f), \qquad f\in\vD_-(z), \nn \\
\vP(z)^{-1}(0,w) &=& (\vE_+(z)w,\ \vE_{-+}(z)w), \qquad w\in\vD_+(z), \label{Pinvnat}
\eea
and $\ker\vE_{-+}(z) := \{w\in\vD_+(z);\ \vE_{-+}(z)w = 0\}$. We write $\vD_\pm$ for $\vD_\pm(z)$.
One has (see Appendix \ref{app-Gr})
\be \label{Dm}
\vD_- = \{f\in L^2(\bR^{2n});\ R'f\in D(P)\}, \qquad \vE f = R'f ,
\ee
and, under \eqref{ass1b},
\be \label{RanP}
\Ran\,\vP(z) = \left\{ (f,g) \in L^2(\bR^{2n})\times\vG ;\
\vE_+g \in L^2(\bR^{2n}) \ \mbox{ and } \ \vE_-f+\vE_{-+}g \in \vG \right\} .
\ee
The entries \eqref{Pinvnat} are the restrictions to $\vD_\pm$ of the operators $\vE$, $\vE_\pm$, $\vE_{-+}$ defined above through $\vS'$.

\begin{theorem} \label{th3.2} Assume \eqref{ass1}, and let $\gamma \ge \gamma_0$ and $z \in \vU_{k+1}(\gamma)$.

\textup{(i)} If $z \not\in\sigma(P)$, then $\vE_{-+}(z)$, with domain $\vD_+$, is a closed,
densely defined operator on $\vG$. It is invertible, with inverse
\be \label{Empinv}
\vE_{-+}(z)^{-1} = - E_- R(z) E_+ \in \vL(\vG),
\ee
that is, $\vE_{-+}\vE_{-+}^{-1} = 1$ on $\vG$ and $\vE_{-+}^{-1}\vE_{-+} = 1$ on $\vD_+$.
Moreover,
\bea
\vE_{-+}(z)^{-1}\vE_-(z) &=& -E_-R(z), \label{EEmp}\\
R(z) &=& \vE(z) - \vE_+(z)\vE_{-+}(z)^{-1}\vE_-(z)
\label{Grushindilated2}
\eea
on $\vD_-$.

\textup{(ii)} $E_-$ is a bijection from $\ker(P-z)$ onto $\ker\vE_{-+}(z)$, with inverse
$\vE_+(z)$. In particular, if $z\notin\sigma(P_0)$, then $z\notin\sigma(P)$
if and only if $\ker\vE_{-+}(z) = \{0\}$.
\end{theorem}

\noindent The proof of Theorem \ref{th3.2} is given in Appendix \ref{app-Gr}.

\begin{cor}\label{cor-Grass1b} Assume \eqref{ass1b}, and let $\gamma\ge\gamma_0$ and $z\in\vU_{k+1}(\gamma)$.
Then $\vS(\bR^{2n})\subset\vD_-$ and $H^2(\bR^n;\bC^\mu)\subset\vD_+$. Consequently, if
$z\notin\sigma(P)$, the operators $\vE_{-+}^{-1}\vE_-$ and $\vE-\vE_+\vE_{-+}^{-1}\vE_-$, defined on
the dense subspace $\vD_-$, extend respectively to $-E_-R(z)\in\vL(L^2(\bR^{2n}),\vG)$ and to
$R(z)\in\vL(L^2(\bR^{2n}))$. In particular
\[
R(z) = \vE(z)-\vE_+(z)\vE_{-+}(z)^{-1}\vE_-(z)\qquad \mbox{ on $\vS(\bR^{2n})$}.
\]
\end{cor}
\pf If $f\in\vS(\bR^{2n})$, then $R'f\in\vS(\bR^{2n})\subset D(P)$ by \eqref{RS}, so that
$f\in\vD_-$ by \eqref{Dm}. If $w\in H^2(\bR^n;\bC^\mu)$, then $L_\theta E_+w\in L^2$ and
$\bX^\theta_{k+1}E^{(k+1)}_-w\in H^1$, so that $\vE_+w\in L^2(\bR^{2n})$. Moreover,
$N_{k+1}E^{(k+1)}_-w\in L^2$ by \eqref{elm4.7}. Hence $\vE_{-+}w\in\vG$ by \eqref{Tentries}, and $w\in\vD_+$ by
\eqref{RanP}. The last assertion follows from \eqref{EEmp} and \eqref{Grushindilated2}.
\ef

\sect{Poles near the thresholds}\label{sect4}

Throughout this section, $\theta$ is fixed with $\varphi\in\,]0,a[$, $k\in\bN$ is fixed, and
\be \label{eta}
z = \gamma k + z', \qquad \eta = \gamma z' = \gamma (z-\gamma k) .
\ee
We keep the notation of Sections \ref{sect2} and
\ref{sect3}. In particular, we use the Grushin data $E_\pm$, $\vE$, $\vE_\pm$, $\vE_{-+}$, the blocks
$\bX^\theta_j$, ${}^t\bX^\theta_j$ of \eqref{Xblocks}, with the convention $\bX^\theta_{-1}=0$, and the block structure \eqref{Tentries}. We also use the Witten Laplacian at the level $k$, $-\Delta^{(k)}_V(\theta)$, of Definition \ref{def-Qk}.
Condition \eqref{ass1b} is assumed throughout this section.

\subsection{Successive approximation}\label{ssec-cf}

By Proposition \ref{prop-tri}, solving $(P(\gamma,\theta)-z)u=f$ amounts to solving the block
tridiagonal system \eqref{system}. We do this by successive elimination of the unknowns. Each
step is a Schur complement (Feshbach map, see \cite{sz}). The resulting operators are given by
the successive approximation \eqref{cf} below, which has the form of a matrix continued fraction (see \cite{wall}). For the
KFP operator, Schur complements with respect to the kernel of the collision operator have been used in
the study of hypocoercivity \cite{bfls}. The elimination down to the lowest level $j=0$, that is, the
case $k=0$, is the
classical derivation of the Smoluchowski equation in the overdamped limit \cite{ne,ls}. Here the
elimination is exact. The levels $\ge k+2$ have already been eliminated in Proposition
\ref{prop-tri}, at the cost of the single block $N_{k+1}$.

The blocks of $P(\gamma,\theta)-z$ away from the level $k$ are controlled by Lemmas \ref{lem3.1} and
\ref{lem3.4}, since $z = \gamma k+z'\in\vU_J(\gamma)$ for $|z'|\le c_1$ and $\gamma$ large. The operators
$\bX^\theta_{k+1}$ and ${}^t\bX^\theta_{k+1}$ are of first order in $x$, with coefficients bounded
together with all their derivatives by \eqref{ass1b}. Moreover, $\|\w{D_x}^sR'\w{D_x}^{-s}\| \le C_s\gamma^{-1}$
by Lemma \ref{lem3.4}. Hence the block $N_{k+1}$ of \eqref{N} satisfies
\be \label{Nest}
\left\| N_{k+1}g\right\|_{H^s} \le C_s\gamma^{-1}\|g\|_{H^{s+2}}, \qquad s \in \bR .
\ee
Here $H^s$ denotes the Sobolev space $H^s(\bR^n_x;\bC^{\mu})$, for the relevant $\mu$.

\subsubsection{The successive approximation and the Schur complement}
The elimination described above produces the following operators. Denote
\be \label{cf}
\begin{aligned}
&T^{(0)} := -z \ (k\ge1), \qquad
T^{(j)} := (\gamma j - z) + \bX^\theta_{j-1}\left(T^{(j-1)}\right)^{-1}{}^t\bX^\theta_{j-1}
\quad (1\le j \le k-1), \\
&T^{(k+1)} := (\gamma(k+1) -z) + N_{k+1} ,
\end{aligned}
\ee
Each $T^{(j)}$ is an operator on $\vG_j$, defined recursively provided the previous one is
invertible. We also put
\be \label{Bk}
\vB_{k-1} := \bX^\theta_{k-2}\left(T^{(k-2)}\right)^{-1}{}^t\bX^\theta_{k-2}
= T^{(k-1)} - \big(\gamma(k-1)-z\big) , \qquad \vB_0 := 0 \ \mbox{ when } k=1 .
\ee
For $0 \le j \le k-1$ let
$\vM^{[0,j]} := -\big( \vE^{(il)}_{-+}\big)_{0\le i,l\le j}$ be the matrix of the first $j+1$
equations of \eqref{system}. It is the matrix of $P^J-z$, $J = \{0,\dots,j\}$, in the bases \eqref{ident}. Hence, by Lemma \ref{lem3.1}, it is invertible for $z\in\vU_J(\gamma)$ and $\gamma\ge\gamma_0$, with inverse the matrix of $R^J(z,\gamma,\theta)$, and it is holomorphic in $(z,\theta)$. This holds in particular for $|z-\gamma k|\le c_1$ and $\gamma$ large, since $\vU_J(\gamma)\supset\{|z-\gamma k|\le c_1\}$. Then $T^{(j)}$ is the Schur complement of $\vM^{[0,j]}$ at the level
$j$, that is, the operator obtained by eliminating the unknowns $u_0,\dots,u_{j-1}$ from those
equations. Similarly, let $J = \{k+1, k+2, \dots\}$ and $P^J = \Pi_JP(\gamma,\theta)\Pi_J$ as in
\S\ref{ssec-high}. Then $T^{(k+1)}$ is the Schur complement of $\vM^{[k+1,\infty[} := P^J-z$ at the level $k+1$. The levels
$\ge k+2$ are eliminated by means of $R'$, which produces the block $N_{k+1}$.
Let $\Theta_a := \{\theta\in\bC;\ |\theta|<a\}$.

\begin{lemma}\label{lem-cf} Assume \eqref{ass1b}. Let $z = \gamma k+z'$ with
$|z'|\le c_1$. There is $\gamma_0$ such that, for $\gamma\ge\gamma_0$ and $\theta\in\Theta_a$, all
the operators $T^{(j)}$ in \eqref{cf} are invertible, and their inverses depend holomorphically on $\theta\in\Theta_a$. One has
\be \label{cfSchur}
\left( T^{(j)}\right)^{-1} = \left[ \left( \vM^{[0,j]}\right)^{-1}\right]_{jj}, \qquad
\left( T^{(k+1)}\right)^{-1} = \left[ \left( \vM^{[k+1,\infty[}\right)^{-1}\right]_{k+1,k+1} ,
\ee
for $0\le j \le k-1$ and, for every $s$,
\be \label{Sigmaest}
\begin{aligned}
\|(T^{(j)})^{-1}\| + \|\w{D_x}^s(T^{(j)})^{-1}\w{D_x}^{-s}\| &\le \f{C_s}{\gamma},
\qquad 0 \le j \le k-1, \\
\|(T^{(k+1)})^{-1}\| + \|\w{D_x}^s(T^{(k+1)})^{-1}\w{D_x}^{-s}\| &\le \f{C_s}{\gamma} .
\end{aligned}
\ee
\end{lemma}
\pf For $\theta\in\Theta_a$ and $z = \gamma k+z'$, $|z'|\le c_1$, one has
$\Re(e^{i\varphi}z) = \gamma k\cos\varphi+\Re(e^{i\varphi}z')$, so that $z\in\vU_J(\gamma)$ for every $J$
lying on one side of the level $k$ as soon as $\gamma(1-\tau)\cos a>c_1$. This condition is
independent of $\theta$. Hence Lemmas \ref{lem3.1} and \ref{lem3.4}, which hold for $|\varphi|<a$, apply to
all the blocks $P^J-z$ below. Moreover, the resolvents $R^J$ are holomorphic in $\theta\in\Theta_a$ by Lemma \ref{lem3.1}. Hence so are all the operators constructed from them below.
We prove the invertibility and \eqref{cfSchur} by induction on $j$, for $k\ge1$ (for $k=0$ only $T^{(1)}$ occurs). For $j=0$,
$T^{(0)} = -z$ is the scalar $-(\gamma k+z')$, $|z| \ge \gamma k - c_1$, and $\vM^{[0,0]} = -z$. Assume $T^{(j-1)}$
invertible and \eqref{cfSchur} at rank $j-1$. Both $\vM^{[0,j]}$ and $\vM^{[0,j-1]}$ are the matrices of $P^{J}-z$ for $J = \{0,\dots,j\}$ and
$J=\{0,\dots,j-1\}$ in the bases \eqref{ident}. They are
invertible by Lemma \ref{lem3.1}, since $z\in\vU_J(\gamma)$ for $\gamma$ large. Write $\vM^{[0,j]}$
as a $2\times 2$ block matrix with respect to the levels $\{0,\dots,j-1\}$ and $\{j\}$. If a
$2\times2$ block matrix and its $(1,1)$ block are invertible, then the Schur complement of the
$(1,1)$ block is invertible and its inverse is the $(2,2)$ entry of the inverse. Here, by
\eqref{Tentries} and the induction hypothesis, that Schur complement is
$(\gamma j - z) + \bX^\theta_{j-1}(T^{(j-1)})^{-1}{}^t\bX^\theta_{j-1} = T^{(j)}$. The same
argument applies to $T^{(k+1)}$. The complementary block of $\vM^{[k+1,\infty[}$ is the matrix
of $P'-z$, which is invertible by Lemma \ref{lem3.1}. The Schur complement is then
$(\gamma(k+1) - z)+N_{k+1}$ by \eqref{N}. Estimates \eqref{Sigmaest} follow from Lemmas \ref{lem3.1} and \ref{lem3.4}, since an entry of a matrix has norm at
most that of the matrix.
\ef

By Lemma \ref{lem-cf}, $T^{(k-1)}$ and $T^{(k+1)}$ are invertible, and we define the operators $\ep_\pm$ by
\be \label{eps}
\gamma \left(T^{(k-1)}\right)^{-1} =: -1 + \ep_-, \qquad
\gamma \left(T^{(k+1)}\right)^{-1} =: 1 + \ep_+ .
\ee
Since, by \eqref{cf} and \eqref{Bk}, $\gamma-T^{(k+1)} = z'-N_{k+1}$ and $\gamma+T^{(k-1)} = \vB_{k-1}-z'$,
\be \label{epsid}
\ep_+ = \left(T^{(k+1)}\right)^{-1}(z' - N_{k+1}), \qquad
\ep_- = \left(T^{(k-1)}\right)^{-1}(\vB_{k-1}-z').
\ee
Moreover, by \eqref{Bk} and \eqref{Sigmaest} for $T^{(k-2)}$ ($\vB_0=0$),
\be \label{Best}
\|\vB_{k-1}g\|_{H^s}\le C_s\gamma^{-1}\|g\|_{H^{s+2}} .
\ee

\begin{cor}\label{cor-eps} Under the hypotheses of Lemma \ref{lem-cf}, for every $s \in \bR$,
\be \label{epsest}
\left\| \ep_\pm g\right\|_{H^s} \le
C_s\left( |z'|\gamma^{-1}\|g\|_{H^s} + \gamma^{-2}\|g\|_{H^{s+2}}\right) .
\ee
\end{cor}
\pf It follows from \eqref{epsid}, \eqref{Sigmaest}, \eqref{Nest} and \eqref{Best}. \ef

Eliminating in the same way all the unknowns $u_j$, $j \neq k$, from the system \eqref{system}
leaves the single operator
\be \label{Sk}
S_k := (\gamma k - z) + \bX^\theta_{k-1}\left(T^{(k-1)}\right)^{-1}{}^t\bX^\theta_{k-1}
+ {}^t\bX^\theta_k\left(T^{(k+1)}\right)^{-1}\bX^\theta_k
\ee
on $\vG_k$. We split $\vG = \vG_k\oplus\vG'$, $\vG' = \bigoplus_{j\neq k}\vG_j$.

\begin{prop}[Schur complement at the level $k$]\label{prop-schur} Assume \eqref{ass1b} and the hypotheses of Lemma
\ref{lem-cf}. Then, as an operator on $\vS(\bR^n;\bC^{\mu_k})$, $S_k$ is the Schur complement of
$\vG'$ in $-\vE_{-+}(z,\allowbreak\gamma,\theta)$.
If moreover $P(\gamma,\theta)-z$ is invertible, then all
the entries of $\vE_{-+}^{-1}$ are bounded, $\big[\vE_{-+}^{-1}\big]_{kk} = -S_k^{-1}$, and
\be \label{SkR}
S_k^{-1} = E^{(k)}_-R(z,\gamma,\theta)E^{(k)}_+.
\ee
\end{prop}
\pf Since the levels $k-1$ and $k+1$ are not adjacent, the restriction of $-\vE_{-+}$ to $\vG'$ is
block diagonal. It is equal to $\vM^{[0,k-1]}\oplus T^{(k+1)}$ in the notation of Lemma \ref{lem-cf}, and
it is invertible by Lemmas \ref{lem3.1} and \ref{lem-cf}. The Schur complement of $\vG'$ in
$-\vE_{-+}$ is \eqref{Sk} by \eqref{Tentries} and \eqref{cfSchur}. The inverses involved are
bounded on every $H^s$ by \eqref{Sigmaest}. If $P-z$ is invertible, then
$\vE_{-+}^{-1} = -E_-RE_+$ by \eqref{Empinv}, which is bounded. The $(k,k)$ entry of
the inverse of a $2\times2$ block matrix is the inverse of the Schur complement. This gives
$\big[\vE_{-+}^{-1}\big]_{kk} = -S_k^{-1}$ and \eqref{SkR}.
\ef

\subsubsection{The block of the level $k$}
Let
\[
{\rm Rem} := {}^t\bX^\theta_k\,\ep_+\,\bX^\theta_k
+ \bX^\theta_{k-1}\,\ep_-\,{}^t\bX^\theta_{k-1}, \qquad {\rm Rem} = {\rm Rem}_{z'}+{\rm Rem}_N ,
\]
\[
{\rm Rem}_{z'} := z'\big({}^t\bX^\theta_k(T^{(k+1)})^{-1}\bX^\theta_k
- \bX^\theta_{k-1}(T^{(k-1)})^{-1}\,{}^t\bX^\theta_{k-1}\big) .
\]
Here $\ep_\pm$ are the operators of \eqref{eps}.

\begin{theorem}\label{th-main-red} Assume \eqref{ass1b} and the hypotheses of Lemma \ref{lem-cf}
with $c_1 = \f 12$. Then, as an identity between operators on $\vS(\bR^n;\bC^{\mu_k})$,
\be \label{main}
\gamma S_k = -\Delta^{(k)}_V(\theta)-\eta + {\rm Rem} ,
\ee
and, for every $s\in\bR$,
\be \label{Remest}
\|{\rm Rem}_{z'}g\|_{H^s} \le C_s|\eta|\gamma^{-2}\|g\|_{H^{s+2}}, \qquad
\|{\rm Rem}_Ng\|_{H^s} \le C_s\gamma^{-2}\|g\|_{H^{s+4}} .
\ee
\end{theorem}
\pf By \eqref{Sk} and \eqref{eps},
\begin{align*}
\gamma S_k &= -\gamma z' + \bX^\theta_{k-1}\left(\gamma\left(T^{(k-1)}\right)^{-1}\right){}^t\bX^\theta_{k-1}
+ {}^t\bX^\theta_k\left(\gamma\left(T^{(k+1)}\right)^{-1}\right)\bX^\theta_k \\
&= -\eta - \bX^\theta_{k-1}{}^t\bX^\theta_{k-1} + {}^t\bX^\theta_k\bX^\theta_k + {\rm Rem},
\end{align*}
which is \eqref{main} by Definition \ref{def-Qk} and Proposition \ref{prop-Qk}. By \eqref{epsid},
\[
{\rm Rem}_N = -{}^t\bX^\theta_k(T^{(k+1)})^{-1}N_{k+1}\bX^\theta_k
+ \bX^\theta_{k-1}(T^{(k-1)})^{-1}\vB_{k-1}{}^t\bX^\theta_{k-1}, \]
and \eqref{Remest} follows from
\eqref{Sigmaest}, \eqref{Nest}, \eqref{Best} and $|z'| = |\eta|\gamma^{-1}$.
\ef

By Lemma \ref{lem-cf}, the operators $T^{(j)}$, $N_{k+1}$, $S_k$ and ${\rm Rem}$ are defined and
holomorphic in $\theta\in\Theta_a$, as operators from $H^{s+2}$ to $H^s$, for $z = \gamma k+z'$, $|z'|\le c_1$.
For $\theta$ real, they are conjugated by $U_\theta$ to their values at $\theta = 0$. The full asymptotic expansions of the poles as $\gamma \to \infty$ rely on the following lemma.

\begin{lemma}[expansion of the remainder]\label{lem-Remexp} Let $r_0>0$. There are
differential operators $M_p(\eta,\theta)$, $p\ge2$,
on $\vG_k$, independent of $\gamma$, such that, for every $N\ge2$
and $s\in\bR$, there is $C$ with
\be \label{Remexp}
\Big\| \Big({\rm Rem}-\sum_{p=2}^N\gamma^{-p}M_p(\eta,\theta)\Big)g\Big\|_{H^s}
\le C\gamma^{-N-1}\|g\|_{H^{s+2N+4}}, \qquad g\in H^\infty,\ |\eta|\le r_0,\ \gamma\ge\gamma_0 .
\ee
The operators $M_p(\eta,\theta)$ are of order at most $p+2$, with coefficients bounded together with all their derivatives, polynomial in $\eta$ and
holomorphic in $\theta\in\Theta_a$.
For $\theta$ real, $M_p(\eta,\theta) = U_\theta M_p(\eta,0)U_\theta^{-1}$.
For $\theta = 0$ and $\eta$ real, the $M_p(\eta,0)$ are symmetric on $H^\infty$.
\end{lemma}
\pf By \eqref{eps} and \eqref{epsid},
$\gamma(T^{(k+1)})^{-1} = (1+Y_+)^{-1}$ and $\gamma(T^{(k-1)})^{-1} = -(1+Y_-)^{-1}$, with
$Y_+ = \gamma^{-1}(N_{k+1}-z')$, $Y_- = \gamma^{-1}(z'-\vB_{k-1})$ and $z' = \gamma^{-1}\eta$. By
\eqref{Nest} and \eqref{Best}, $\|Y_\pm g\|_{H^s}\le C\gamma^{-2}\|g\|_{H^{s+2}}$. Moreover, $(1+Y_\pm)^{-1}$
are bounded on every $H^s$, uniformly in $\gamma$, by \eqref{Sigmaest}. The identity
$(1+Y)^{-1} = \sum_{r=0}^{N}(-Y)^r + (1+Y)^{-1}(-Y)^{N+1}$ reduces the expansion of $\ep_\pm$ to
that of $N_{k+1}$ and $\vB_{k-1}$. By Proposition \ref{prop-Rprime}(ii),
$R'_{k+2} = E^{(k+2)}_-R'E^{(k+2)}_+$, hence $N_{k+1}$, has an expansion in powers of $\gamma^{-1}$ of
the form \eqref{Remexp}, with coefficients which are differential operators in $x$, polynomial in
$\eta$ and holomorphic in $\theta\in\Theta_a$. For $\vB_{k-1}$ one argues by induction on $j$ with
$\gamma(T^{(j)})^{-1} = (j-k)^{-1}(1+Y_j)^{-1}$,
\[
Y_j = \gamma^{-1}(j-k)^{-1}\big(-z'+\bX^\theta_{j-1}(T^{(j-1)})^{-1}\,{}^t\bX^\theta_{j-1}\big),
\]
using
\eqref{Sigmaest} in the same way. Collecting the powers of $\gamma^{-1}$ gives \eqref{Remexp}. There is
no term $p=1$ since $\|Y_\pm g\|_{H^s} = O(\gamma^{-2})\|g\|_{H^{s+2}}$.

For $\theta$ real, $U_\theta$ acts in $x$ only and commutes with the $E^{(j)}_\pm$ and the $\Pi_j$.
Hence $T^{(j)}$, $S_k$ and ${\rm Rem}$ at $\theta$ are conjugated by $U_\theta$ to those at $\theta=0$.
Since the coefficients of an expansion \eqref{Remexp} are unique,
$M_p(\eta,\theta) = U_\theta M_p(\eta,0)U_\theta^{-1}$. Let $\theta = 0$ and $\eta$, hence $z$, be
real. Then $X^0_l = \p_{x_l}+\f12\p_{x_l}V$ has real coefficients, so that
${}^t\bX^0_j = (\bX^0_j)^*$. By induction on $j$ in \eqref{cf}, $T^{(j)}$ is symmetric for
$0\le j\le k-1$. Since $\bJ$ commutes with $\Pi'$ and equals $(-1)^{k+2}$ on $\Ran\Pi_{k+2}$,
\eqref{tRprime} and the reality of the coefficients give $R'_{k+2}(z)^* = R'_{k+2}(z)$, so that $N_{k+1}$ and $T^{(k+1)}$ are symmetric.
Hence ${\rm Rem}$ is symmetric on $H^\infty$, and so are the $M_p(\eta,0)$, again by uniqueness of the
coefficients.

Finally, by \eqref{wN}, the coefficient of $\gamma^{-q}$ in the expansion of $R'$ contains at most $q-1$
factors $L'$, hence contains at most $q-1$ derivatives in $x$. The remainder of Proposition
\ref{prop-Rprime}(ii) at the order $N$ loses $N$ derivatives. Since $\bX^\theta_j$ and ${}^t\bX^\theta_j$
are of order $1$, the coefficient of $\gamma^{-p}$ in $Y_\pm$, $Y_j$ contains at most $p$ derivatives in $x$, and that of
${\rm Rem}$ at most $p+2$. Counting the derivatives lost by the remainders in the same way gives
\eqref{Remexp} with the loss $2N+4$.
\ef

\subsection{Main results}\label{ssec-main}
Throughout the rest of this section, \eqref{ass1b} holds, $\delta\in\,]0,\f\pi2[$,
$\sigma\in\,]\f34,1[$, $\Lambda_\delta$ is as in Section \ref{sect1}, and $\chi^\sigma_j$ is as in \eqref{chisig}.
The main step is the stability estimate of Theorem \ref{th-red} below. In its proof, the high frequencies $|\xi|\ge\gamma^\sigma$ are
controlled by the free resolvent (Lemma \ref{lem-loc}).
At low frequencies, the Grushin problem of Section \ref{sect3} reduces the
equation to the level $k$, where the Schur complement $\gamma S_k$ is a perturbation of $-\Delta^{(k)}_V(\theta)-\eta$.
The following lemma then recovers $u$ from $E_-u$ and its high-frequency part.

\begin{lemma}\label{lem-uE} Let $\sigma\in\,]\f34,1[$, $z\in\bC$ with $|z-\gamma k|<1$, and $u \in D(P(\gamma,\theta))$,
and set $f = (P(\gamma,\theta)-z)u$. Then, for $\gamma$ large,
\be \label{uE}
\|u\| \le 2\|E_-u\| + C\left( \|\chi^\sigma_1u\| + \gamma^{-1}\|f\|\right) .
\ee
\end{lemma}
\pf By density, it suffices to prove \eqref{uE} for $u$ in the core $\vS(\bR^{2n})$. Since $z\in\vU_{k+1}(\gamma)$ for $\gamma$ large,
Lemma \ref{lem3.1} applies. Applying $\Pi'$ to
$(P-z)u=f$ gives $\Pi'u = R'\Pi'f - \vK u$ with $\vK := R'\Pi'L_\theta\Pi$, so that,
by Lemma \ref{lem3.1},
\be \label{Kbound}
\|\vK u\| \le \|u\| + C\gamma^{-1}\|f\| .
\ee
Since $[\chi^\sigma_0,R'] = R'[W_\theta,\chi^\sigma_0]R'$ and
\[
[\chi^\sigma_0,L_\theta] = -[\chi^\sigma_0,(\na V)(e^\theta x)]\cdot\na_v = O(\gamma^{-\sigma})\na_v,
\]
one has
\[
\chi^\sigma_0\vK = \vK\chi^\sigma_0 + R'[W_\theta,\chi^\sigma_0]\vK + R'\Pi'[\chi^\sigma_0,L_\theta]\Pi .
\]
Here
\[
\|\vK\chi^\sigma_0\| \le \|R'\|\,\|L_\theta\Pi\w{D_x}^{-1}\|\,\|\w{D_x}\chi^\sigma_0\|
\le C\gamma^{\sigma-1},
\]
and the last two terms are $O(\gamma^{-1-\sigma})$ by Lemma \ref{lem3.1}. With
\eqref{Kbound}, this gives
\[
\|\chi^\sigma_0\Pi'u\| \le C\gamma^{-1}\|f\| + \|\chi^\sigma_0\vK u\|
\le C\left( \gamma^{\sigma-1}\|u\| + \gamma^{-1}\|f\|\right) .
\]
As
\[
\|\chi^\sigma_0\Pi u\| \le \|\Pi u\| = \|E_-u\|,
\]
 one gets
\[
\|u\| \le \|\chi^\sigma_0u\| + \|\chi^\sigma_1u\| \le \|E_-u\| + C\gamma^{\sigma-1}\|u\| +
\|\chi^\sigma_1u\| + C\gamma^{-1}\|f\|,
\]
and \eqref{uE} follows since $\sigma<1$.
\ef

Let $K$ be an operator of finite rank on $\vG_k$ of the
form $K = \sum_{i}\w{g_i,\cdot}f_i$ with $f_i,g_i\in H^4$, possibly $K=0$, and put
\[
P^K(\gamma,\theta) := P(\gamma,\theta) + \gamma^{-1}E^{(k)}_+KE^{(k)}_-, \qquad
-\Delta^{(k)}_{V,K} := -\Delta^{(k)}_V(\theta)+K .
\]
For $\eta \notin \sigma(-\Delta^{(k)}_{V,K})$ let
\[
M' := \sup_{-4\le s\le 0}\left\| \w{D_x}^{s+2}(-\Delta^{(k)}_{V,K}-\eta)^{-1}\w{D_x}^{-s}\right\| .
\]

The finite-rank term $K$ in the following theorem is introduced to treat $\eta$ near an eigenvalue
$\eta_0$ of $-\Delta^{(k)}_V(\theta)$. With $K = c\,\gU$, $\gU$ the Riesz projection at $\eta_0$, $-\Delta^{(k)}_{V,K}-\eta$ is
invertible near $\eta_0$ and the theorem gives the invertibility of $P^K(\gamma,\theta)-z$ there.
The poles of $R(z,\gamma,\theta)$ are then the zeros of a determinant of size $k_0$ (Proposition
\ref{prop-det}).
Here $\eta_0\in\Lambda_\delta$, and $\eta$ is taken in $\Lambda_{\delta/2}$. Theorem \ref{th-red} is then applied with $\delta/2$ in place of $\delta$.
Since $\eta = \gamma(z-\gamma k)$, $\arg\eta = \arg(z-\gamma k)$, and this condition reads $|\arg(e^{2i\varphi}(z-\gamma k))|\ge\f\delta2$.
That is, $z$ stays away from the half-line $\gamma k+e^{-2i\varphi}\bR_+\subset\sigma_{\rm ess}(P(\gamma,\theta))$.
For such $z$, one has $\eta\notin\sigma_{\rm ess}(-\Delta^{(k)}_V(\theta))$.

The following theorem is basic for establishing a bijection between the poles of the KFP operator near the threshold $\gamma k$ and those of the Witten Laplacian at the level $k$.

\begin{theorem}[stability estimate]\label{th-red} There exist $c_1\in\,]0,\f12]$ and $C>0$ and, for every $M_0>0$, there exists $\gamma_0$ such that the following holds. Let $\gamma \ge \gamma_0$ and let $\eta \in \Lambda_\delta$ with $|\eta| \le c_1\gamma$,
$\eta \notin \sigma(-\Delta^{(k)}_{V,K})$ and $M'\le M_0$. Then $z = \gamma k + \eta/\gamma$
belongs to the resolvent set of $P^K(\gamma,\theta)$ and
\be \label{resest}
\|(P^K(\gamma,\theta)-z)^{-1}\| \le C\gamma(1+M') .
\ee
\end{theorem}
\pf Let $u\in\vS(\bR^{2n})$, $f = (P^K-z)u$, $u_j = E^{(j)}_-u$ and
$\chi_j = \chi^\sigma_j$. It suffices to prove
\be \label{apriori}
\|u\| \le C\gamma(1+M')\|f\| .
\ee
Indeed, by density of $\vS(\bR^{2n})$ in $D(P)$ for the graph norm, \eqref{apriori} then holds
on $D(P)$. Hence $P^K-z$ is injective. It is invertible since $z$
does not belong to
\[
\sigma_{\rm ess}(P^K) = \sigma_{\rm ess}(P_0) .
\]

The perturbation $B := \gamma^{-1}E^{(k)}_+KE^{(k)}_-$ satisfies $\Pi'B = B\Pi' = 0$. Hence Lemma
\ref{lem-uE} holds for $P^K$, and so do \eqref{r1}, \eqref{r2} and \eqref{r4}, Proposition \ref{prop-schur} and
Theorem \ref{th-main-red}.
The $(k,k)$ block of $\vE_{-+}$ is replaced by
$\vE^{(kk)}_{-+}-\gamma^{-1}K$, $S_k$ by $S_k+\gamma^{-1}K$ and $-\Delta^{(k)}_V(\theta)$ by $-\Delta^{(k)}_{V,K}$, while
${\rm Rem}$ is unchanged.

For the high frequency part, one has
\[
\chi_1u = \chi_1R_0(f-W_\theta u-Bu) .
\]
By
\eqref{R0}, $\|\chi_1R_0\| \le C\gamma^{1-2\sigma}$, since, as in the proof of Corollary
\ref{cor2.4}(b),
\[
|E^\theta_k(\xi,\gamma)-z| \ge c\gamma^{-1}|\xi|^2 \ \mbox{ for } |\xi|\ge R\gamma^\sigma, \qquad
|E^\theta_j(\xi,\gamma)-z|\ge c\gamma \ \mbox{ for } j\ne k .
\]
With \eqref{chi1small} and $\|B\|\le C\gamma^{-1}$,
\be \label{hf}
\|\chi_1u\| \le C\gamma^{-\f12}(\|u\|+\|f\|) .
\ee

 By \eqref{r4}, $-\vE_{-+}E_-u = \vE_-f =: \phi$, with
$\phi_j = E^{(j)}_-f$ for $j\ne k+1$ and
$\phi_{k+1} = E^{(k+1)}_-f+{}^t\bX^\theta_{k+1}E^{(k+2)}_-R'f$. With respect to
$\vG = \vG_k\oplus\vG'$, the block of $-\vE_{-+}$ on $\vG'$ is $\bT := \vM^{[0,k-1]}\oplus T^{(k+1)}$.
Its off-diagonal blocks $\vX_+:\vG'\to\vG_k$, $\vX_-:\vG_k\to\vG'$ have the entries $\bX^\theta_{k-1}$,
$-{}^t\bX^\theta_k$ and $-{}^t\bX^\theta_{k-1}$, $\bX^\theta_k$ by \eqref{Tentries}. Writing
$u' = (u_j)_{j\ne k}$, one gets
\[
u' = \bT^{-1}(\phi'-\vX_-u_k), \qquad (\gamma S_k+K)u_k = \psi := \gamma\big(\phi_k-\vX_+\bT^{-1}\phi'\big).
\]
By Theorem \ref{th-main-red}, $\gamma S_k+K = (-\Delta^{(k)}_{V,K}-\eta)(1+\vR_K{\rm Rem})$ on $\vS$, with
$\vR_K := (-\Delta^{(k)}_{V,K}-\eta)^{-1}$. Hence
\be \label{uk}
u_k = -\vR_K{\rm Rem}\,u_k + \vR_K\psi .
\ee

Making use of \eqref{Remest} and $\|\chi_0\w{D_x}^2\| \le C\gamma^{2\sigma}$, one has
\begin{align*}
\|\chi_0\vR_K{\rm Rem}\| &\le \|\vR_K\w{D_x}^2\|\,\|\w{D_x}^{-2}{\rm Rem}_{z'}\| \\
&\quad+ \|\chi_0\w{D_x}^2\|\,\|\w{D_x}^{-2}\vR_K\w{D_x}^4\|\,\|\w{D_x}^{-4}{\rm Rem}_N\| \\
&\le CM_0\big(c_1\gamma^{-1}+\gamma^{2\sigma-2}\big) \le \f12
\end{align*}
for $\gamma$ large, since $\sigma<1$. The operator $\bT^{-1}$ is bounded on every $H^s$ with norm $\le C\gamma^{-1}$ (Lemmas
\ref{lem3.1} and \ref{lem-cf}). The operators $\vX_\pm$ and ${}^t\bX^\theta_{k+1}$ are of first order, and
$\|\w{D_x}^sR'\w{D_x}^{-s}\|\le C\gamma^{-1}$. With $\|\vR_K\|+\|\vR_K\w{D_x}^2\|\le 2M'$, this gives
$\|\vR_K\psi\| \le CM'\gamma\|f\|$. Applying $\chi_0$ to \eqref{uk} and using
$\|u_k\|\le\|\chi_0u_k\|+\|\chi_1u\|$, we obtain $\|u_k\| \le CM'\gamma\|f\| + 2\|\chi_1u\|$. In the
same way,
\[
\|\chi_0\bT^{-1}\vX_-\| \le \|\chi_0\w{D_x}\|\,\|\w{D_x}^{-1}\bT^{-1}\vX_-\| \le C\gamma^{\sigma-1}
\]
and $\|\chi_0\bT^{-1}\phi'\| \le C\gamma^{-1}\|f\|$, so that
$$\|u'\| \le C\gamma^{\sigma-1}\|u_k\| + C\gamma^{-1}\|f\| + \|\chi_1u\|.$$
Therefore
$\|E_-u\| \le C\gamma(1+M')\|f\| + C\|\chi_1u\|$, and Lemma \ref{lem-uE} and \eqref{hf} give
$$\|u\| \le C\gamma(1+M')\|f\| + C\gamma^{-\f12}\|u\|,$$
whence \eqref{apriori} for $\gamma$ large.
\ef

With $C_0$ as in \eqref{highen} and $c_1>0$, let
\be \label{Dk}
\vD^{(k)}_\gamma := \left\{ z \in \bC; \ \f{C_0}\gamma \le |z-\gamma k| \le c_1, \
\gamma(z-\gamma k) \in \Lambda_\delta\right\} .
\ee

\begin{cor}[pole-free domains near thresholds]\label{th-absence} Let $\Im\theta = \varphi\in\,]0,a[$. There are $c_1>0$ and $\gamma_0$ such that for $\gamma\ge\gamma_0$ the
resolvent $R(z,\gamma,\theta)$ is holomorphic, with $\|R(z,\gamma,\theta)\|\le C\gamma$, in $\vD^{(k)}_\gamma$.
In particular $P(\gamma)$ has neither eigenvalues nor resonances in $\vD^{(k)}_\gamma$.
\end{cor}
\pf For $z \in \vD^{(k)}_\gamma$, $\eta = \gamma(z-\gamma k)\in\Lambda_\delta$ with
$C_0 \le |\eta| \le c_1\gamma$, and \eqref{highen} gives $M'\le C_\delta$,
independent of $\gamma$ and $\eta$. Theorem \ref{th-red}, with $K=0$ and $M_0 = C_\delta$, then applies.
\ef

Let $\Im\theta = \varphi\in\,]0,a[$, let $\eta_0\in\Lambda_\delta$ be an eigenvalue of $-\Delta^{(k)}_V(\theta)$ of multiplicity $k_0$, and let $\ep_0>0$ be
such that $\{|\eta-\eta_0|\le 2\ep_0\}\cap\sigma(-\Delta^{(k)}_V(\theta)) = \{\eta_0\}$ and
$\{|\eta-\eta_0|\le 2\ep_0\}\subset\Lambda_{\delta/2}$. Let $\gU = \gU^\theta_k$ be the Riesz
projection of $-\Delta^{(k)}_V(\theta)$ at $\eta_0$, $c := 4\ep_0$, and
$\vQ_k(\eta) := -\Delta^{(k)}_V(\theta)+c\,\gU-\eta$.
 Since
\[
\sigma\big(-\Delta^{(k)}_V(\theta)+c\gU\big) =
\Big(\sigma\big(-\Delta^{(k)}_V(\theta)\big)\setminus\{\eta_0\}\Big)\cup\{\eta_0+c\},
\]
$\vQ_k(\eta)$ is invertible for
$|\eta-\eta_0|\le2\ep_0$. For $N\ge0$ let
\be \label{WN}
W_N(\eta,\gamma,\theta) := \sum (-1)^r\gamma^{-(p_1+\cdots+p_r)}\,
\gU \vQ_k^{-1}M_{p_1}\vQ_k^{-1}M_{p_2}\cdots \vQ_k^{-1}M_{p_r}\vQ_k^{-1}\gU .
\ee
Here the sum is over $r\ge0$ and $p_1,\dots,p_r\ge2$ with $p_1+\cdots+p_r\le N$, the $M_p$ are those
of Lemma \ref{lem-Remexp}, and $\vQ_k = \vQ_k(\eta)$.

Let $\gamma$ be large and $|\eta-\eta_0|\le2\ep_0$. The operator $\gamma S_k+c\gU$ is invertible on $\vG_k$ by Proposition \ref{prop-det} below, and we put
\[
d(\eta,\gamma,\theta) := \det\Big(1-c\,\gU(\gamma S_k+c\gU)^{-1}\gU\Big|_{\Ran\gU}\Big), \qquad
d_N(\eta,\gamma,\theta) := \det\big(1-c\,W_N(\eta,\gamma,\theta)|_{\Ran\gU}\big) .
\]

\begin{prop}[reduced determinant]\label{prop-det} Let $\varphi\in\,]0,a[$. For $\gamma$ large and
$|\eta-\eta_0|\le2\ep_0$, $\gamma S_k + c\gU$ is invertible on $\vG_k$ with
$\|(\gamma S_k+c\gU)^{-1}\|\le C$, and the following holds.

\textup{(i)} $z = \gamma k+\eta\gamma^{-1}$, which satisfies $|\arg(e^{2i\varphi}(z-\gamma k))|\ge\f\delta2$, is a pole of $R(z,\gamma,\theta)$ of multiplicity $m$ if and
only if $\eta$ is a zero of order $m$ of $d(\cdot,\gamma,\theta)$.

\textup{(ii)} For every $N$, $d = d_N + O(\gamma^{-N-1})$ uniformly in $|\eta-\eta_0|\le2\ep_0$.

\textup{(iii)} If $\eta_0<0$, then $\eta_0$ is an isolated eigenvalue of $-\Delta^{(k)}_V(\theta)$ for every
$\theta\in\Theta_a$. Moreover, $\ep_0$ may be chosen independent of $\theta$, and $W_N$ is defined and holomorphic
in $(\eta,\theta)$ for $|\eta-\eta_0|\le2\ep_0$, $\theta\in\Theta_a$. For $\theta$ real,
$W_N(\eta,\gamma,\theta) = U_\theta W_N(\eta,\gamma,0)U_\theta^{-1}$. For $\theta = 0$ and $\eta$
real, $W_N(\eta,\gamma,0)$ is Hermitian.
\end{prop}
\pf Apply Theorem \ref{th-red} with $\delta/2$ and $K = c\gU$, which is of finite rank with range and
co-range in $H^\infty$. Since $M'$ is bounded on the compact set
$\{|\eta-\eta_0|\le2\ep_0\}$, it gives the invertibility of $P^K-z$ with
$\|(P^K-z)^{-1}\|\le C\gamma$. By Proposition \ref{prop-schur} for $P^K$,
$(S_k+\gamma^{-1}c\gU)^{-1} = E^{(k)}_-(P^K-z)^{-1}E^{(k)}_+$, which gives the bound on
$(\gamma S_k+c\gU)^{-1}$.

(i) $P = P^K-B$ with $B = \gamma^{-1}cE^{(k)}_+\gU E^{(k)}_-$ of finite rank,
and $P^K-z$ is invertible. By the Weinstein--Aronszajn formula \cite[IV.6]{k}, the
multiplicity of $z$ as an eigenvalue of $P$ is the order of the zero at $z$ of
\[
\det\big(1-B(P^K-z)^{-1}\big) = \det\big(1-\gamma^{-1}c\,\gU E^{(k)}_-(P^K-z)^{-1}
E^{(k)}_+\gU|_{\Ran\gU}\big).
\]
The right-hand side is $d(\eta,\gamma,\theta)$. Since the map $\eta\mapsto z$ is affine,
orders of zeros are the same in $z$ and in $\eta$.

(ii) The space $\Ran\gU$ is finite dimensional and contained in $H^\infty$, and $\vQ_k(\eta)^{-1}$ is
bounded on every $H^s$. For $g\in\Ran\gU$ let $h_N := \sum_{r=0}^N(-\vQ_k^{-1}{\rm Rem})^r\vQ_k^{-1}g$. By
\eqref{main}, $\gamma S_k+c\gU = \vQ_k+{\rm Rem}$, and a telescoping computation gives
$(\gamma S_k+c\gU)h_N = g-(-{\rm Rem}\,\vQ_k^{-1})^{N+1}g$. The inverse $(\gamma S_k+c\gU)^{-1}$ is a left
inverse on $H^\infty$. Indeed, for $h\in H^\infty$, complete $h$ by $u' = -\bT^{-1}\vX_-h$ as in the proof of
Theorem \ref{th-red}. This gives an element of $H^\infty(\bR^n;\bC^\mu)$ to which Theorem \ref{th3.2} (i) and Corollary \ref{cor-Grass1b} apply
($\vE_{-+}^{-1}\vE_{-+} = 1$ on $\vD_+\supset H^2$). Hence $(\gamma S_k+c\gU)^{-1}g = h_N + O(\gamma^{-2N-2})\|g\|$ by
\eqref{Remest}. Expanding each ${\rm Rem}$ in $h_N$ by Lemma \ref{lem-Remexp} and collecting the powers
of $\gamma^{-1}$ gives $\gU(\gamma S_k+c\gU)^{-1}\gU = W_N+O(\gamma^{-N-1})$ on $\Ran\gU$, hence (ii).

(iii) The number $\eta_0<0$ lies at positive distance from $e^{-2i\varphi}[0,\infty[$ for $|\varphi|<a$.
The eigenvalues of $-\Delta^{(k)}_V(\theta)$ outside the closed sector between $\bR_+$ and $e^{-2i\varphi}\bR_+$
are those of $-\Delta^{(k)}_V$ \cite{ac}. Hence the first assertion holds.
Moreover, $\gU^\theta_k$ and $\vQ_k(\eta)^{-1}$ are
holomorphic in $\theta\in\Theta_a$, since $\{-\Delta^{(k)}_V(\theta)\}$ is a holomorphic family of type (A).
With
Lemma \ref{lem-Remexp}, this gives the holomorphy of $W_N$ and, for $\theta$ real, the conjugation by
$U_\theta$.
For $\theta = 0$ and $\eta$ real, $\gU$ is an orthogonal projection, $c$ is real, $\vQ_k(\eta)$
is self-adjoint and the $M_p(\eta,0)$ are symmetric.
The adjoint of each term of \eqref{WN} is the
term with $(p_1,\dots,p_r)$ reversed, so that $W_N$ is Hermitian.
\ef

The following result gives the locations of the poles near $\gamma k$ of $R(z, \gamma, \theta)$ in terms of the eigenvalues of the Witten Laplacian at the level $k$, $-\Delta_V^{(k)}(\theta)$.

\begin{theorem}[bijection]\label{th-bij} Let $\Im\theta = \varphi\in\,]0,a[$, let $\eta_0$ be an eigenvalue of $-\Delta^{(k)}_V(\theta)$ of
multiplicity $k_0$ with $\eta_0 \in \Lambda_\delta$, and let
$\ep_0>0$ be such that $\{|\eta-\eta_0|\le 2\ep_0\}\cap\sigma(-\Delta^{(k)}_V(\theta)) = \{\eta_0\}$ and
$\{|\eta-\eta_0|\le2\ep_0\}\subset \Lambda_{\delta/2}$. Then there is $\gamma_0$ such that for
$\gamma\ge\gamma_0$, $R(z,\gamma,\theta)$ has exactly $k_0$ poles, counted with
multiplicity, in the disc $\{|z-\gamma k - \eta_0\gamma^{-1}| < \ep_0\gamma^{-1}\}$.
\end{theorem}
\pf Let $d$, $d_N$ be as in Proposition \ref{prop-det}. By Proposition
\ref{prop-det}(i), the poles of $R$ in the disc, with multiplicity, are the points
$z = \gamma k+\eta\gamma^{-1}$, $\eta$ running over the zeros of $d(\cdot,\gamma,\theta)$ in
$\{|\eta-\eta_0|<\ep_0\}$, with multiplicity. On $\Ran\gU$, $-\Delta^{(k)}_V(\theta)-\eta_0$ is nilpotent and
$W_0 = (-\Delta^{(k)}_V(\theta)+c-\eta)^{-1}\gU$, so that
$d_0(\eta) = \big((\eta_0-\eta)/(\eta_0+c-\eta)\big)^{k_0}$. As $d-d_0 = O(\gamma^{-1})$, Rouch\'e's
theorem gives exactly $k_0$ zeros of $d(\cdot,\gamma,\theta)$ in $\{|\eta-\eta_0|<\ep_0\}$ for $\gamma$
large.
\ef

For $r>0$ let
\be \label{annulus}
\vA_\gamma(r) := \left\{ z;\ \f r\gamma \le |z-\gamma k| \le c_1, \
\gamma(z-\gamma k)\in\Lambda_\delta\right\}
\ee
be the truncated \emph{annulus}, where $c_1$ is as in Corollary \ref{th-absence}.

\begin{cor}\label{cor-sharp} Let $\Im\theta = \varphi\in\,]0,a[$, $\delta \in\,]0,\f\pi2[$, let $r>0$, and assume that $-\Delta^{(k)}_V(\theta)$ has no eigenvalue on the boundary of
$\Lambda_\delta\cap\{|\eta|\ge r\}$. Then,
for $\gamma \ge \gamma_2(r)$, $P(\gamma)$ has in $\vA_\gamma(r)$
exactly as many eigenvalues and resonances, with multiplicity, as $-\Delta^{(k)}_V$ has eigenvalues and
resonances in $\Lambda_\delta\cap\{r\le|\eta|\le c_1\gamma\}$. Each of them is located at
$z = \gamma k + \eta_j\gamma^{-1}+O(\gamma^{-1-1/k_j})$, where $\eta_j$ is an eigenvalue of $-\Delta^{(k)}_V(\theta)$ and $k_j$ is its multiplicity.
\end{cor}
\pf By \eqref{highen}, $\sigma(-\Delta^{(k)}_V(\theta))\cap\Lambda_\delta \subset \{|\eta|<C_0\}$. The set $\Lambda_\delta\cap\{r\le|\eta|\le C_0\}$ is compact and
lies at distance $\ge r\sin\delta$ from $\sigma_{\rm ess}(-\Delta^{(k)}_V(\theta)) = e^{-2i\varphi}[0,\infty[$.
Hence $-\Delta^{(k)}_V(\theta)$ has there finitely many eigenvalues $\eta_1,\dots,\eta_p$, of
multiplicities $k_1,\dots,k_p$. None of them depends on $\gamma$, and they exhaust
$\sigma(-\Delta^{(k)}_V(\theta))\cap\Lambda_\delta\cap\{r\le|\eta|\le c_1\gamma\}$. By assumption, no $\eta_j$
lies on the boundary of $\Lambda_\delta\cap\{r\le|\eta|\le c_1\gamma\}$. Choose $\ep_0>0$ so
small that the discs $\{|\eta-\eta_j|\le2\ep_0\}$ are disjoint and contained in the interior of this set, and let $\vA$ be
$\Lambda_\delta\cap\{r\le|\eta|\le c_1\gamma\}$ with these discs removed. On $\vA$, the resolvent of
$-\Delta^{(k)}_V(\theta)$ exists and
\[
M' := \sup_{\eta\in \vA}\ \sup_{-4\le s\le0}\left\| \w{D_x}^{s+2}(-\Delta^{(k)}_V(\theta)-\eta)^{-1}\w{D_x}^{-s}\right\|
<\infty .
\]
This bound is independent of $\gamma$. For $|\eta|\ge C_0$, it follows from \eqref{highen}. For
$r\le|\eta|\le C_0$, it follows from the continuity of $\eta \mapsto (-\Delta^{(k)}_V(\theta)-\eta)^{-1}$ in
$\vL(H^s,H^{s+2})$ on a compact subset of the resolvent set of the $\gamma$-independent operator
$-\Delta^{(k)}_V(\theta)$. Since $|\eta|\le c_1\gamma$ on $\vA$, Theorem \ref{th-red}, with $K=0$ and $M_0=M'$, applies for $\gamma$ large. Hence $R$ is holomorphic, with
$\|R\|\le C\gamma(1+M')$, at every $z$ with $\gamma(z-\gamma k) \in \vA$. Theorem \ref{th-bij} then
counts $k_j$ poles inside each of the discs $\{|z-\gamma k - \eta_j\gamma^{-1}|<\ep_0\gamma^{-1}\}$. With $d$ and $d_0$ as in its proof,
$d-d_0 = O(\gamma^{-1})$, while $|d_0(\eta)|\ge c|\eta-\eta_j|^{k_j}$ for $|\eta-\eta_j|\le2\ep_0$.
Rouch\'e's theorem shows that $d$ has no zero in
$\{\ep_0\le|\eta-\eta_j|\le2\ep_0\}$, and that its $k_j$ zeros lie in $\{|\eta-\eta_j|\le C\gamma^{-1/k_j}\}$.
This gives the location of the poles.
\ef

\begin{rmk} The operator $-\Delta_V $ has no eigenvalues (except possibly $0$), but it may have
resonances. Replace $V$ by $g V$, where $g\ge 1$ is a large coupling parameter. The associated
Witten Laplacian becomes
\[
\begin{aligned}
-\Delta_{gV} &= -\Delta + \f{g^2}4|\nabla V(x)|^2 - \f g2 \Delta V(x) \\
&= h^{-2} \Big(-h^2 \Delta + \f14|\nabla V(x)|^2 - \f h2 \Delta V(x)\Big),
\qquad h = g^{-1}.
\end{aligned}
\]
There is a large literature on the semiclassical analysis of resonances of Schr\"odinger operators
(see, for example, \cite{hs2}).
For a fixed
large coupling constant, $-\Delta_{gV}$ has resonances due to the tunneling effect. By
Corollary \ref{cor-sharp} with $k=0$, they give rise to low-lying resonances of $P(\gamma)$ for $\gamma$ large. The real parts of these resonances determine the Kramers escape rates.
\end{rmk}

The reduced determinant of Proposition \ref{prop-det} gives complete asymptotic expansions of the
poles as $\gamma\to\infty$. Let $\Im\theta = \varphi\in\,]0,a[$, $k\in\bN$, let
$\eta_0\in\Lambda_\delta$ be an eigenvalue of $-\Delta^{(k)}_V(\theta)$ of multiplicity $k_0$ (for $k=0$, a resonance
of $-\Delta_V$ with $-2\varphi<\arg\eta_0<0$), let $\ep_0$ be as in Theorem \ref{th-bij}, and let
$z_1(\gamma),\dots,z_{k_0}(\gamma)$ be the poles of $R(z,\gamma,\theta)$ in
$\{|z-\gamma k-\eta_0\gamma^{-1}|<\ep_0\gamma^{-1}\}$, counted with multiplicity (Theorem
\ref{th-bij}). One has the following result.

\begin{theorem}[asymptotic expansions]\label{th3.10} There exist $c_{l,m}\in\bC$ such that, for every $N\in\bN^*$,
\be \label{expansion}
z_l(\gamma) = \gamma k + \eta_0\gamma^{-1} + \sum_{m=1}^{N}c_{l,m}\gamma^{-1-\f m{k_0!}}
+ O\left( \gamma^{-1-\f{N+1}{k_0!}}\right), \qquad \gamma \to \infty, \quad 1\le l\le k_0 .
\ee
\end{theorem}
\pf By Proposition \ref{prop-det}(i), $z_l(\gamma) = \gamma k+\gamma^{-1}\eta_l(\gamma)$, where the
$\eta_l(\gamma)$ are the zeros of $d(\cdot,\gamma,\theta)$ in $\{|\eta-\eta_0|<\ep_0\}$, with multiplicity.
Fix $N$. The function $d_N$ is a polynomial in $\gamma^{-1}$ with coefficients holomorphic in $\eta$,
and its value at $\gamma^{-1} = 0$ is $d_0 = ((\eta_0-\eta)/(\eta_0+c-\eta))^{k_0}$ (see the proof of Theorem
\ref{th-bij}). By the Weierstrass preparation theorem, in the variables $(\eta,\gamma^{-1})$ near $(\eta_0,0)$,
$d_N = e\,p_N$ \cite[Theorem 6.1.1]{ho}.
Here $e$ does not vanish, $p_N$ is a monic polynomial of degree $k_0$ in $\eta$
with coefficients holomorphic in $\gamma^{-1}$, and $p_N = (\eta-\eta_0)^{k_0}$ at $\gamma^{-1} = 0$.
By the
Newton--Puiseux theorem \cite[II.1.2]{k}, its roots are convergent series in powers of $\gamma^{-1/p}$, where $p\le k_0$ is the length of the
corresponding cycle. Hence they are series in powers of $\gamma^{-1/k_0!}$:
\[
\zeta_l(\gamma) = \eta_0+\sum_{m\ge1}a_{l,m}\gamma^{-m/k_0!}, \quad 1\le l\le k_0. \]
Let
$\rho = C_1\gamma^{-(N+1)/k_0}$ and let $\vO_\rho$ be the $\rho$-neighborhood of
$\{\zeta_1(\gamma),\dots,\zeta_{k_0}(\gamma)\}$. On $\partial\vO_\rho$, $|d_N|\ge c\rho^{k_0} > |d-d_N|$
if $C_1$ is large, by Proposition \ref{prop-det}(ii). Rouch\'e's theorem, applied to each connected
component of $\vO_\rho$, whose diameter is at most $2k_0\rho$, shows that the zeros $\eta_l(\gamma)$ of
$d$ can be labeled so that $|\eta_l(\gamma)-\zeta_l(\gamma)|\le 2k_0\rho$. The same argument applied to
$d_N$ and $d_{N+1}$ shows that the $a_{l,m}$, $m\le N$, do not depend on $N$, since $N/k_0!<(N+1)/k_0$, and \eqref{expansion}
holds with $c_{l,m} = a_{l,m}$.
\ef

\begin{cor}\label{cor-intpow} Let $k\ge1$ and let $\eta_0<0$ be an eigenvalue of $-\Delta^{(k)}_V$ of multiplicity $k_0$. Then there exist $b_{l,m}\in\bR$
such that, for every $N\in\bN^*$,
\be \label{intpow}
z_l(\gamma) = \gamma k+\eta_0\gamma^{-1}+\sum_{m=1}^Nb_{l,m}\gamma^{-1-m}+O\big(\gamma^{-2-N}\big), \qquad \gamma\to\infty, \quad 1\le l\le k_0 .
\ee
\end{cor}

\pf Let $d$, $d_N$ and $W_N$ be as in Proposition \ref{prop-det}, with $\ep_0$ independent of $\theta\in\Theta_a$ (Proposition
\ref{prop-det}(iii)), and let $h := \gamma^{-1}$. If $\theta_1\in\Theta_a$ and $\phi_1,\dots,\phi_{k_0}$ is a basis of $\Ran\gU^{\theta_1}_k$, then
$\gU^\theta_k\phi_1,\dots,\gU^\theta_k\phi_{k_0}$ is a basis of $\Ran\gU^\theta_k$ for $\theta$ near $\theta_1$, in which the matrix of $W_N|_{\Ran\gU}$ is holomorphic in
$\theta$. Hence $d_N(\eta,\gamma,\theta)$ is holomorphic in $\theta\in\Theta_a$. For $\theta$ real, $\Ran\gU^\theta_k = U_\theta\Ran\gU^0_k$ and
$W_N(\eta,\gamma,\theta) = U_\theta W_N(\eta,\gamma,0)U_\theta^{-1}$, so that $d_N(\eta,\gamma,\theta) = d_N(\eta,\gamma,0)$. By the identity theorem, this holds for
$\theta\in\Theta_a$.

At $\theta = 0$, $\gU^0_k$ is the orthogonal projection onto $\ker(-\Delta^{(k)}_V-\eta_0)$ and $\vQ_k(\eta)^{-1}\gU^0_k = (\eta_0+c-\eta)^{-1}\gU^0_k$. Since the
$M_p$ start at $p = 2$, \eqref{WN} gives
\[
\begin{aligned}
&1-cW_N(\eta,\gamma,0)|_{\Ran\gU^0_k} = (\eta_0+c-\eta)^{-1}G(\eta,h), \\
&G(\eta,h) := \eta_0-\eta-c(\eta_0+c-\eta)\sum_{p=2}^Nh^pF_p(\eta).
\end{aligned}
\]
Here the $F_p$ are holomorphic in $\{|\eta-\eta_0|\le2\ep_0\}$, with values in the operators on $\Ran\gU^0_k$. As $W_N(\eta,\gamma,0)$ is Hermitian for
$\eta$ real and $\gamma$ large (Proposition \ref{prop-det}(iii)), the $F_p(\eta)$ are Hermitian for $\eta$ real. Hence $G(\eta,h)^* = G(\bar\eta,h)$ for
$h\in\bR$. Moreover, $G(\eta,h)-G(\bar\eta,h) = (\eta-\bar\eta)(-1+O(h^2))$ for $|\eta-\eta_0|<\ep_0$, by Cauchy's inequality. If $h$ is real and small,
$|\eta-\eta_0|<\ep_0$ and $G(\eta,h)u = 0$ with $u\neq0$, then
\[
0 = \w{u,G(\eta,h)u}-\overline{\w{u,G(\eta,h)u}} = \w{u,\big(G(\eta,h)-G(\bar\eta,h)\big)u} = (\eta-\bar\eta)\big(-1+O(h^2)\big)\|u\|^2 ,
\]
so that $\eta$ is real. Now $d_N$ is a polynomial in $h$, and we write it $d_N(\eta,h,\theta)$.
For every small $h\in\bR$, of either sign, the zeros
of $d_N(\cdot,h,\theta) = d_N(\cdot,h,0)$ in $\{|\eta-\eta_0|<\ep_0\}$, that is, the roots of the polynomial $p_N$ of the proof of Theorem \ref{th3.10}, are real. This excludes cycles of
length $p\ge2$, as in the proof of Rellich's theorem \cite[II.6.1]{k}. Indeed, let $\zeta(h) = \eta_0+\sum_{m\ge1}a_mh^{m/p}$. The reality of the $p$ branches for
$h>0$ gives $a_m\in\bR$ and $e^{2i\pi m/p}a_m\in\bR$, hence $a_m = 0$ unless $2m/p\in\bN$.
Thus $\zeta(h) = A(h)+h^{1/2}B(h)$, where $A$ and $B$ are
power series with real coefficients.
Moreover, $B = 0$ since $\zeta(h)$ is real for $h<0$. Therefore the roots $\zeta_l$ of $p_N$ are power series in
$\gamma^{-1}$ with real coefficients, and the proof of Theorem \ref{th3.10}, with $k_0!$ replaced by $1$ and $N$ by $k_0(N+1)$, gives
\eqref{intpow}. \ef

\sect{Poles below the threshold \texorpdfstring{$\gamma k$}{gamma k}, \texorpdfstring{$k \ge 1$}{k >= 1}}\label{sec-poles}

In this section, we assume $k \ge 1$ and study the poles produced by the negative eigenvalues of $-\Delta^{(k)}_V$, which exist in some cases.

\subsection{Negative eigenvalues of \texorpdfstring{$-\Delta^{(k)}_V$}{-Delta_V^(k)}}\label{ssec-neg}
 The results of this subsection may be known. Since we could not find appropriate references, we give a succinct proof. The closest results we know concern the semiclassical Witten Laplacians on forms \cite{hs4}, which are, however, non-negative. The particularity here is that $-\Delta_V$ is non-negative and $-\dG_k(\mbox{\rm Hess}\,V)$ is a $\mu_k \times \mu_k$ matrix which can produce negative eigenvalues.

Throughout this subsection, $V\in C^\infty(\bR^n)$ is non-constant. Instead of \eqref{ass1b},
we only assume that, for some $\rho>0$,
\be \label{ass1d}
|\p^\alpha_xV(x)| \le C_\alpha\w x^{-\rho-|\alpha|+1}, \qquad x\in\bR^n,\ \alpha\in\bN^n,\ |\alpha|\ge1.
\ee
 $-\Delta^{(k)}_V$ is then defined by
\eqref{QkE} with $\theta = 0$.
When $n \le 2$ and $\rho>1$ in \eqref{ass1d}, the potential $V$ is bounded and $|\na V|^2 \le C\w x^{-2\rho}$ is integrable, so
that
\be \label{Q}
\vQ := \int_{\bR^n}e^{-V(x)}\,\na V(x)\otimes \na V(x)\,dx
\ee
is a well defined real symmetric non-negative $n\times n$ matrix.

\begin{theorem}\label{th-neg} Assume \eqref{ass1d} and let $V$ be non-constant.

\textup{(i)} For every $n$ and $k$, $\inf\sigma(-\Delta^{(k)}_V)<0$ if and only if there is
$g \in C^\infty_0(\bR^n;\bC^{\mu_k})$ such that
\be \label{crit}
\int_{\bR^n}e^{-V}|\na g|^2\,dx \ <\ \int_{\bR^n}e^{-V}
\w{\dG_k(\mbox{\rm Hess}\,V(x))g(x),g(x)}\,dx .
\ee
In this case, $\inf\sigma(-\Delta^{(k)}_V)$ is a negative eigenvalue of $-\Delta^{(k)}_V$ of finite multiplicity.

\textup{(ii)} Let $n\le 2$ and $\rho>1$ in \eqref{ass1d}. Then $-\Delta^{(k)}_V$ has at least one negative
eigenvalue for every $k\ge1$.
\end{theorem}
\pf \emph{(i)} Let $h = \mbox{Hess}\,V$. For $g \in C^\infty_0(\bR^n;\bC^{\mu_k})$ put
$u = e^{-V/2}g$. Since $X_l = \p_l + \f12\p_lV$ at $\theta=0$, $X_lu = e^{-V/2}\p_lg$, and by
\eqref{QkE}
\be \label{quadform}
\begin{aligned}
\w{(-\Delta^{(k)}_V)u,u} &= \sum_l \|X_lu\|^2 - \int e^{-V}\w{\dG_k(h)g,g}\,dx \\
&= \int e^{-V}|\na g|^2\,dx - \int e^{-V}\w{\dG_k(h)g,g}\,dx .
\end{aligned}
\ee
The map $g\mapsto e^{-V/2}g$ is a bijection of $C^\infty_0(\bR^n;\bC^{\mu_k})$ onto itself, and this space is a core
for $-\Delta^{(k)}_V$. Hence $-\Delta^{(k)}_V$ has a negative quadratic form value if and only if \eqref{crit} holds
for some $g$. Since $\sigma_{\rm ess}(-\Delta^{(k)}_V)=[0,\infty[$, this is equivalent to
$\inf\sigma(-\Delta^{(k)}_V)<0$, and $\inf\sigma(-\Delta^{(k)}_V)$ is then a discrete eigenvalue.

\emph{(ii)} Since $\rho>1$, $V$ is bounded and $c \le e^{-V}\le C$. Take
$g = g_R\otimes s$ with $s\in\bC^{\mu_k}$, $|s|=1$, and $g_R\in C^\infty_0(\bR^n)$ real-valued, $0\le g_R\le1$, depending on a parameter $R$ and such that
\be \label{gR}
g_R\to1 \ \mbox{ pointwise}, \qquad \int|\na g_R|^2\,dx\to0 \qquad (R\to\infty).
\ee
Such $g_R$ are given below for $n=1,2$. Integrating by parts and using
$\w{x}^{n-1}e^{-V}|\p_mV| \to 0$ at infinity,
\[
\int e^{-V}\p_l\p_mV\,dx = \int e^{-V}\p_lV\p_mV\,dx = \vQ_{lm}
\]
by \eqref{Q}. The integrals are absolutely convergent because $|{\rm Hess}\,V|\le C\w{x}^{-1-\rho}$ and
$|\na V|^2 \le C\w{x}^{-2\rho}$ are integrable for $n\le2$, $\rho>1$. By linearity,
$\int e^{-V}\dG_k(h(x))\,dx = \dG_k(\vQ)$. Moreover
\[
{\rm tr}\,\vQ = \int e^{-V}|\na V|^2\,dx>0
\]
since $V$ is non-constant, so that the largest eigenvalue $\nu_{\max}(\vQ)$ of $\vQ$ is $>0$. By \eqref{gR}, $\int e^{-V}|\na g|^2\,dx\to0$, while, by dominated convergence, the right-hand side of \eqref{crit} tends to $\w{\dG_k(\vQ)s,s}$.

\emph{Case $n=1$.} Take $g_R = \chi(x/R)$, $\chi\in C^\infty_0(-2,2)$, $\chi=1$ on
$[-1,1]$, $0\le\chi\le1$. Then \eqref{gR} holds, since $\int|g_R'|^2 \le CR^{-1}$. Choosing $s$ with
$\w{\dG_k(\vQ)s,s} = k\nu_{\max}(\vQ)$ (Lemma \ref{lem-dG}(ii)), we get from \eqref{quadform}

\[
\w{(-\Delta^{(k)}_V)u,u} \le CR^{-1} - k\nu_{\max}(\vQ)+o(1)<0
\]
 for $R$ large, and (i) applies.

\emph{Case $n=2$.} We use instead the logarithmic cut-off $g_R=1$ on $\{|x|\le R\}$,
$g_R = \log(R^2/|x|)/\log R$ for $R\le|x|\le R^2$, $g_R=0$ for $|x|\ge R^2$ (regularized near
the two circles, so that $g_R\in C^\infty_0$). Then $\int|\na g_R|^2 \le C/\log R \to 0$ and $g_R\to1$ pointwise, so that \eqref{gR} holds, and the
argument is the same.
\ef

By Theorem \ref{th-neg} (ii), when $n=1,2$, $-\Delta_V^{(k)}$ has negative eigenvalues for every $k\ge 1$ and every non-constant $V$ satisfying \eqref{ass1d} with $\rho>1$.
Another way to see the existence of negative eigenvalues is to study the large coupling limit. We replace $V$ by $gV$, where $g\ge1$ is a large coupling constant. We assume that $V$ has a
non-degenerate local minimum at $x_0$, and that $H = \mbox{\rm Hess}\,V(x_0)>0$ has the eigenvalues
$\lambda_1,\dots,\lambda_n>0$. We denote by $\lambda_{\max}(H)$ the largest of them. Set, for $m,\alpha \in \bN^n$ with $|\alpha| = k$,
\be \label{mualpha}
\mu(m,\alpha) := \sum_{l=1}^n \lambda_l\left( m_l - \alpha_l\right) \ \ge\ -k\,\lambda_{\max}(H) .
\ee
The lower bound holds because $\sum_l\alpha_l = k$. Only finitely many pairs $(m,\alpha)$ give a
negative value, since $\mu(m,\alpha)<0$ forces
\[
\sum_l\lambda_lm_l < \sum_l\lambda_l\alpha_l \le k\lambda_{\max}(H).
\]
 Finally, with $h = g^{-1}$,
\be \label{Ph}
\begin{aligned}
&P_h := h^2\big(-\Delta^{(k)}_{gV}\big) = -h^2\Delta_x + W(x) - h\,M(x), \\
&W := \f 14|\na V|^2, \qquad M := \f 12\Delta V + \dG_k\big( {\rm Hess}\,V\big) .
\end{aligned}
\ee
This is a semiclassical Schr\"odinger operator, whose principal symbol is the \emph{scalar} function
$p(x,\xi) = |\xi|^2+W(x)$. By \eqref{ass1d}, $W \ge 0$ and $|M(x)| \le C\w x^{-1-\rho}\to0$ as
$|x|\to\infty$, so that $\sigma_{\rm ess}(P_h) = [0,\infty[$.

Theorem \ref{th-largeg} below is close to the semiclassical approximation for the Witten Laplacian on $p$-forms studied in \cite{hs4}.
The method is the same: localization near the critical points and reduction to a harmonic oscillator with explicit eigenvalues.
 However, the Witten Laplacian on forms is non-negative, since it is the square of $d_V+d_V^*$. Here, we are interested in the negative eigenvalues of $-\Delta^{(k)}_{gV}$.

\begin{theorem}[large coupling, all dimensions]\label{th-largeg} Let $V$, $x_0$ and $H$ be as
above, and assume that $V$ satisfies \eqref{ass1d}.

\textup{(i)} For $g$ large, $\inf\sigma(-\Delta^{(k)}_{gV})$ is a discrete eigenvalue and
\be \label{infQk}
\inf \sigma\left( -\Delta^{(k)}_{gV}\right) \ \le\ -k\,\lambda_{\max}(H)\,g + O(1) \ <\ 0 .
\ee

\textup{(ii)} Assume in addition that $x_0$ is the \emph{only} critical point of $V$, and let
$\ep>0$ be such that $-\ep$ is \emph{not} one of the values \eqref{mualpha}. Then, for $g$ large,
the eigenvalues of $-\Delta^{(k)}_{gV}$ in $]-\infty,-\ep g]$ are, with multiplicity, exactly the
\be \label{Emalpha}
E_{m,\alpha}(g) = g\,\mu(m,\alpha) + O(g^{\f 12}) , \qquad \mu(m,\alpha) < -\ep .
\ee
\end{theorem}
\pf This is a critical semiclassical eigenvalue problem, in the sense that the energy level $0$ is also the limit of $W$ at infinity. See \cite{hs1} for non-critical energy levels.  Since $\na V(x_0) = 0$ and
${\rm Hess}\,W(x_0) = \f 12 H^2$, $W$ vanishes at $x_0$ and that well is non-degenerate. Setting
$x = x_0 + h^{\f 12}y$,
\be \label{model}
P_h = h\left( K_H + \vE_h\right) , \qquad
K_H := -\Delta_y + \f 14 |Hy|^2 - \f 12 {\rm tr}\,H - \dG_k(H) ,
\ee
where $\vE_h(y) = h^{\f 12}c_3(y) + O(h)\,O(\w y^4)$. Here $c_3$ is an \emph{odd} polynomial, given by the
Taylor terms of odd order of $W$ and of $-hM$ at $x_0$. The operator $K_H$, a harmonic oscillator
tensored with the constant matrix $-\dG_k(H)$, has spectrum
\[
\left\{ \sum_l\lambda_l\left( m_l+\tfrac 12\right) - \tfrac 12 \sum_l \lambda_l
- \sum_l \alpha_l\lambda_l \ ;\ m \in \bN^n, \ \alpha \in \bN^n, \ |\alpha| = k\right\}
= \big\{ \mu(m,\alpha)\big\} ,
\]
since ${\rm tr}\,H = \sum_l\lambda_l$ and $\sigma(\dG_k(H)) = \{\sum_l\alpha_l\lambda_l;\ |\alpha|=k\}$
by Lemma \ref{lem-dG}(ii). An eigenfunction of $K_H$ for $\mu(m,\alpha)$ is $\phi_{H,m}\otimes s_\alpha$, where
$\phi_{H,m}$ is the Hermite function of $-\Delta_y+\f 14|Hy|^2$ for $\sum_l\lambda_l(m_l+\f 12)$ and
$s_\alpha$ is an eigenvector of $\dG_k(H)$ for $\sum_l\alpha_l\lambda_l$.

\emph{(i)} Let $\chi \in C^\infty_0(\bR^n)$ be equal to $1$ near $0$, let $\alpha$ be given by Lemma
\ref{lem-dG}(ii), so that $\sum_l\alpha_l\lambda_l = k\lambda_{\max}(H)$, and put
\[
u_h(x) = h^{-\f n4}\,\chi(x-x_0)\,\phi_H\big( h^{-\f 12}(x-x_0)\big)\otimes s_\alpha ,
\]
where $\phi_H(y) = (2\pi)^{-\f n4}(\det H)^{\f14}e^{-\f14 y\cdot Hy}$ is the normalized ground state of $-\Delta_y+\f14|Hy|^2$.
Then $\|u_h\| = 1 + O(h^\infty)$ and, by \eqref{Ph} and \eqref{model},
$\w{P_hu_h,u_h} = h\,\mu(0,\alpha) + O(h^2)$. The term of order $h^{\f 32}$ vanishes because $c_3$
is odd while $|\phi_H|^2$ is even. Multiplying by $h^{-2}=g^2$ and using
$\mu(0,\alpha) = -k\lambda_{\max}(H)$, the min-max principle gives \eqref{infQk}. As
$\sigma_{\rm ess}(-\Delta^{(k)}_{gV}) = [0,\infty[$, the infimum is a discrete eigenvalue.

\emph{(ii)} \emph{A lower bound away from $x_0$.} Fix $\ep>0$ and $\delta>0$, and let $R$ be such that
$|M(x)| \le \f\ep 4$ for $|x| \ge R$. Since $x_0$ is the only critical point of $V$, $W$ vanishes
only at $x_0$, so that
\be \label{cdelta}
W(x) \ \ge\ c_\delta \ >\ 0 \qquad \mbox{ on the compact set } \
F_\delta := \left\{ \delta \le |x-x_0|,\ |x| \le 2R\right\} .
\ee
Let $w \in H^1(\bR^n;\bC^{\mu_k})$ be supported in
$\{|x-x_0|\ge\delta\}$. On $F_\delta$, one has $W - h|M| \ge c_\delta - Ch \ge 0$ by \eqref{cdelta}.
For $|x| \ge 2R$, one has $W \ge 0$ and $|M| \le \f\ep 4$, so that $W - h|M| \ge -\f{\ep h}4$. Hence
\be \label{outside}
\w{P_hw,w} \ \ge\ -\f\ep 2\,h\,\|w\|^2 , \qquad w \in H^1(\bR^n;\bC^{\mu_k}), \quad
\supp w \subset \{|x-x_0|\ge\delta\} .
\ee

\emph{Harmonic approximation.} Let $\chi_0,\chi_1 \in C^\infty(\bR^n)$ with
$\chi^2_0+\chi^2_1 = 1$, $\chi_0 = 1$ on $\{|x-x_0|\le\delta\}$ and
$\supp\chi_0 \subset \{|x-x_0|<2\delta\}$. The IMS localization formula and \eqref{outside} give
\[
P_h \ = \ \chi_0P_h\chi_0 + \chi_1P_h\chi_1 - h^2\left( |\na\chi_0|^2+|\na\chi_1|^2\right)
\ \ge\ \chi_0P_h\chi_0 - \f\ep2 h\,\chi^2_1 - C_\delta h^2 .
\]
Let $\lambda \le -\ep h$ and $h$ be small. Since $\chi_0u\in H^1_0(B(x_0,2\delta))$ for $u\in H^1(\bR^n;\bC^{\mu_k})$, $P_h$ has at most as many eigenvalues below $\lambda$
as the Dirichlet realization of $P_h$ in $B(x_0,2\delta)$ has below $\lambda + C_\delta h^2$. By
Dirichlet bracketing, it has at least as many as the Dirichlet realization in $B(x_0,\delta)$ has below
$\lambda$. For $\delta$ small, $B(x_0,2\delta)$ contains the single non-degenerate well
\eqref{model} and no other zero of $W$. Hence the harmonic approximation \cite{hs1} applies to
both realizations. It gives for their
eigenvalues below $-\f\ep 2h$ the values $h\,\mu(m,\alpha)+O(h^{\f32})$, in bijection with the
eigenvalues $\mu(m,\alpha)$ of $K_H$. Since
$-\ep$ is not an eigenvalue of $K_H$ and these are finitely many,
\[
\kappa := \min\left\{ \left| \mu(m,\alpha)+\ep\right| ;\ \mu(m,\alpha) < 0 \right\} \ >\ 0 .
\]
Hence, for $h$ small, each of these values lies at distance $\ge\f\kappa 2 h$ from the threshold
$-\ep h$, and the two counts at $\lambda = -\ep h$ agree. Multiplying by $h^{-2}=g^2$ gives
\eqref{Emalpha}.
\ef

\subsection{Resonances by phase-space deformation}\label{ssec-phase}

By Corollary \ref{cor-intpow}, obtained with the reduced determinant of Proposition \ref{prop-det}, the imaginary parts of the poles
produced by the negative eigenvalues of $-\Delta^{(k)}_V$ are of order $O(\gamma^{-\infty})$. To show that they are in fact exponentially small, we use an analytic deformation in the
frequencies $\xi$ with $|\xi|\ge c\gamma$, $0<c<1$. It is a complex deformation in phase space, generated by the escape
function $x\cdot\xi/|\xi|$. Complex deformations in phase space were introduced by Helffer and Sj\"ostrand to study
resonances \cite{hs2}. Tunneling in phase space has been studied in \cite{gm15,mns,na95b}, and escape functions in phase space were used for the eigenvalues of the KFP
operator in \cite{hss}.

In this subsection $\gamma = 1$, $P = H_0+L$, and
\[
\widehat P_0(\xi,\theta) := \widehat P_0(\xi,1,\theta) = H_0+ie^{-\theta}v\cdot\xi.
\]
 For $u\in L^2(\bR^{2n})$ let
$\hat u(\xi,v) := \int e^{-ix\cdot\xi}u(x,v)\,dx$ and $\widehat\vH := L^2(\bR^n_\xi,(2\pi)^{-n}d\xi;L^2(\bR^n_v))$. By \cite[Corollary 2.4]{w4},
$D(P) = D(H_0)\cap D(v\cdot\na_x)$, whose Fourier image is $\vD_P := \{w\in\widehat\vH;\ H_0w,\ (\xi\cdot v)w\in\widehat\vH\}$, and
\be \label{PFourier}
\widehat{Pu}(\xi) = \widehat P_0(\xi,0)\,\hat u(\xi)-\int\widehat{\na V}(\xi-\xi')\cdot\na_v\hat u(\xi')\,\f{d\xi'}{(2\pi)^n} .
\ee
Let $S_\ep := \{\zeta\in\bC^n;\ \Re\zeta\neq0,\ |\Im\zeta|<\ep|\Re\zeta|\}$, $\rho_* := \min(\rho,n-\f12)$ and $\gm(r) := 1+r^{\rho_*-n}$, which is locally
integrable on $\bR^n$ as a function of $r = |\xi|$. Under \eqref{ass1b}, $V$ is holomorphic in
$\vO_a := \{|\zeta|<C\}\cup\{|\Im\zeta|<\tan a\,|\Re\zeta|\}$. Let
\be \label{deltaV}
\delta_V := \sup\big\{\delta>0;\ V \mbox{ extends holomorphically to } \vO_a\cup\{\zeta\in\bC^n;\ |\Im\zeta|<\delta\}\big\}\ \ge\ C\sin a .
\ee

\begin{lemma}\label{lem-hatV} For every $\delta_1<\delta_V$ there is $\ep_*>0$ such that $\widehat{\na V}$ is a locally integrable function which extends
holomorphically to $S_{\ep_*}$, with $\widehat{\na V}(-\zeta) = \overline{\widehat{\na V}(\bar\zeta)}$ and
\be \label{hatVbound}
|\widehat{\na V}(\zeta)|\le C\,e^{-\delta_1|\Re\zeta|}\,\gm(|\Re\zeta|), \qquad \zeta\in S_{\ep_*} .
\ee
\end{lemma}
\pf This is a variant of the Paley--Wiener theorem \cite[Theorem IX.13]{rs2}. By  a contour deformation of $\bR^n$, one can show that there exists $b_0>0$ such that
 \beas
 |\widehat{\na V}(\zeta)| &\le & Ce^{-\delta_1|\Re\zeta|}\int\w x^{-\rho}e^{-b_0|x||\Re\zeta|/2}dx \\
 &\le & C'e^{-\delta_1|\Re\zeta|}\gm(|\Re\zeta|).
 \eeas
  The symmetry holds since $\na V$ is real on $\bR^n$. \ef

\subsubsection{The deformed operator}
Let $0<c'_\Psi<1$, and let $\chi_\Psi$ be a non-decreasing function of class $C^\infty$ on $[0,\infty[$, with values in $[0,1]$.
We assume that $\chi_\Psi(0) = 0$, $\chi_\Psi = 1$
on $[c'_\Psi,\infty[$, and $r\mapsto\chi_\Psi(r)/r$ extends to a smooth even function on $\bR$. Let $c_\Psi := \sup\{r;\ \chi_\Psi(r) = 0\}$, so that
$0\le c_\Psi<c'_\Psi$. The basic choice in this
subsection is $\chi_\Psi(r) = r/c'_\Psi$ near $0$, for which $c_\Psi = 0$.
In \S\ref{sec-flux} we use a function $\chi_\Psi$ vanishing on $[0,c_\Psi]$, $c_\Psi>0$.
Put, for $s\in\bC$,
\[
\Psi(\xi) := \chi_\Psi(|\xi|)\,\f\xi{|\xi|}, \qquad \phi_s(\xi) := \xi+s\Psi(\xi) = e^{-\theta_s(\xi)}\xi, \qquad
\theta_s(\xi) := -\log\Big(1+\f{s\chi_\Psi(|\xi|)}{|\xi|}\Big) .
\]
The radial form of the deformation is essential: it preserves the vector field $v\cdot \xi$. If it is destroyed as in
the analytic distortion in $x$ when $n\ge 2$,  the operator is no longer  subelliptic (\cite[Proposition 4.2]{w4}).

Fix $\delta_1<\delta_V$ and $\ep_*$ as in Lemma \ref{lem-hatV}. Since $\Psi$ is Lipschitz, with constant $C_\Psi$, there is $s_1>0$ such that, for $s\in\Sigma := \{|s|<s_1\}$ and
$\xi\neq\xi'$, $\phi_s(\xi)-\phi_s(\xi')\in S_{\ep_*}$, $|\Re(\phi_s(\xi)-\phi_s(\xi'))|\ge\f12|\xi-\xi'|$, $|\theta_s(\xi)|\le\f\pi4$, $J_s := \det\na\phi_s\neq0$, and, for
$s$ real, $\phi_s$ is a diffeomorphism. The point is that $\phi_s(\xi)$ is a complex multiple of $\xi$, so that
\be \label{Fdil}
\bF_s(\xi) := H_0+i\phi_s(\xi)\cdot v = \widehat P_0\big(\xi,\theta_s(\xi)\big), \qquad \sigma(\bF_s(\xi)) = \{j+\phi_s(\xi)\cdot\phi_s(\xi);\ j\in\bN\} .
\ee
Thus the deformed free fibres are dilated free fibres, with a parameter depending on $\xi$. Let $s_0$ be as in Theorem \ref{th2.1} with
$\gamma = 1$ and $\phi_0 = \f\pi4$.

\begin{lemma}[fibred subelliptic estimate]\label{lem-subfib} There is $C$ such that, for $s\in\Sigma$, $\xi\in\bR^n$, $\varkappa\ge s_0$ and
$f\in D(H_0)$,
\be \label{subfib}
\|(\varkappa+H_0)f\|+\w\xi^{\f23}\|f\|+\|(\xi\cdot v)f\|\le C\,\|(\bF_s(\xi)+\varkappa)f\| .
\ee
\end{lemma}
\pf By \eqref{Fdil}, this is \eqref{e2.2b} with $\gamma = 1$ and $\theta = \theta_s(\xi)$, $|\theta_s(\xi)|\le\f\pi4$, extended from $\vS(\bR^n_v)$ to
$D(H_0) = D(\bF_s(\xi))$ by closure. \ef

For $s$ real,
$(\vV_sw)(\xi) := J_s(\xi)^{1/2}w(\phi_s(\xi))$ is unitary on $\widehat\vH$ and preserves $\vD_P$. For $s\in\Sigma$, let $\bP_s := \bF_s+\bW_s$ on $\widehat{\vH}$, where $\bF_s$ is the multiplication by $\bF_s(\xi)$ and
\be \label{Ps}
\begin{aligned}
&(\bW_sw)(\xi) := \int K_s(\xi,\xi')\cdot\na_vw(\xi')\,d\xi', \\
&K_s(\xi,\xi') := -(2\pi)^{-n}J_s(\xi)^{\f12}\,\widehat{\na V}(\phi_s(\xi)-\phi_s(\xi'))\,J_s(\xi')^{\f12} .
\end{aligned}
\ee
With $\kappa_\delta(\zeta) := e^{-\delta|\zeta|}\gm(|\zeta|/2)$, integrable on $\bR^n$ for $\delta>0$, \eqref{hatVbound} gives $|K_s(\xi,\xi')|\le C\kappa_{\delta_1/2}(\xi-\xi')$. Let

\[
\vD_{\rm an} := {\rm span}\{x^\beta e^{-|x|^2/2}g(v);\ \beta\in\bN^n,\ g\in\vS(\bR^n)\}.
\]

\begin{prop}\label{prop-def}
\textup{(a)} For $s\in\Sigma$, $\|\bW_s(1+H_0)^{-1/2}\|\le C$ and $\|(\bW_s-\bW_{\Re s})(1+H_0)^{-1/2}\|\le C|\Im s|$. The family $\{\bP_s\}_{s\in\Sigma}$ is holomorphic
of type \textup{(A)}, with domain $\vD_P$, and $\bP_s = \vV_s\widehat P\vV_s^{-1}$ for $s$ real.

\textup{(b)} For $s\in\Sigma$ and $w\in\vD_P$, $\big|\Re\w{w,\bP_sw}-\w{w,H_0w}\big|\le C|\Im s|\,\|(1+H_0)^{1/2}w\|^2$.

\textup{(c)} There is $C_1$ such that, for $f,g\in\vD_{\rm an}$, $s\in\Sigma$ and $\Re z<-C_1$, $\bP_s-z$ is invertible and
\be \label{matrixel}
\w{g,(P-z)^{-1}f} = \w{\vV_{\bar s}\hat g,(\bP_s-z)^{-1}\vV_s\hat f}_{\widehat\vH} .
\ee
Here $\vV_s\hat f$ is defined by the same formula.
\end{prop}
\pf (a) The Schur test gives the first bound.
Cauchy's inequality on polydiscs of radius proportional to $|\xi-\xi'|$ in
$S_{\ep_*}$ gives $|K_s-K_{\Re s}|\le C|\Im s|\kappa_{\delta_1/4}(\xi-\xi')$.
Indeed, the arguments of $\widehat{\na V}$ differ by $i\Im s(\Psi(\xi)-\Psi(\xi'))$, and $J_s^{1/2}-J_{\Re s}^{1/2}$ is $O(|\Im s|)$.
Hence $\bP_s-\bP_{\Re s}$ is $H_0^{1/2}$-bounded.
Moreover, $H_0$ is $\bP_{\Re s}$-bounded by \cite[Corollary 2.4]{w4}, since $\bP_{\Re s}$ is unitarily
equivalent to $\widehat P$ by $\vV_{\Re s}$, which commutes with $H_0$.
Therefore $\{\bP_s\}$ is holomorphic of type (A).
For $s$ real, the change of variables
$\xi\mapsto\phi_s(\xi)$ in \eqref{PFourier} gives $\bP_s = \vV_s\widehat P\vV_s^{-1}$.
(b) One has $\Re\w{w,i\phi_s\cdot vw} = -\Im s\int\w{w,(\Psi\cdot v)w}d\xi/(2\pi)^n$, $\bW_{\Re s}$
is skew-adjoint, and (a) applies.
(c) By (b), $\|(1+H_0)^{1/2}w\|\le C\|(\bP_s-z)w\|$ for $\Re z<-C_1$.
The set of $s\in\Sigma$ for which $\bP_s-z$ is
invertible contains $\Sigma\cap\bR$ and is open and closed in $\Sigma$, hence is $\Sigma$.
Both sides of \eqref{matrixel} are holomorphic in $s$, since $\hat f$
is entire with Gaussian decay on $\{\phi_s(\xi)\}$. They coincide for $s$ real. \ef

Recall that $\vD_P$, defined after \eqref{PFourier}, is the Fourier image of $D(P)$, and that $\bF_s$ denotes the multiplication by $\bF_s(\xi)$
on $\widehat\vH$, with domain $\vD_P$.

\begin{cor}[subelliptic estimate]\label{th-subdef} There is $C$ such that, for $s\in\Sigma$, $\varkappa\ge s_0$ and $w\in\vD_P$,
\be \label{subdef}
\|(\varkappa+H_0)w\|+\|\w\xi^{\f23}w\|+\|(\xi\cdot v)w\|\le C\,\|(\bF_s+\varkappa)w\| .
\ee
There is $\varkappa_1$ such that the same estimate holds for $\bP_s$ in place of $\bF_s$ if $\varkappa\ge\varkappa_1$.
\end{cor}
\pf Integrate \eqref{subfib} in $\xi$. For $\bP_s$,
\[
\|\bW_sw\|\le C\|(1+H_0)^{1/2}w\|\le C\varkappa^{-1/2}\|(\varkappa+H_0)w\|
\]
 by Proposition
\ref{prop-def}(a), and this term is absorbed for $\varkappa$ large. \ef

\subsubsection{Spectrum and resonances}
\begin{figure}[ht]
\centering
\includegraphics[width=0.4\textwidth]{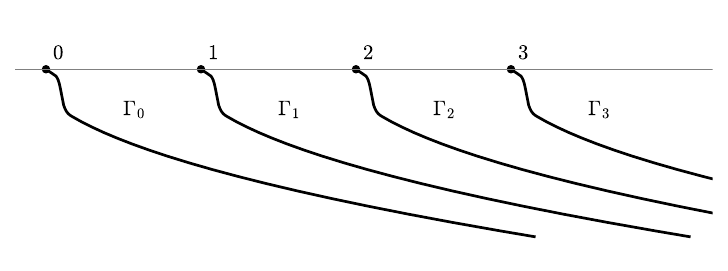}
\caption{\small The curves $\Gamma_0,\dots,\Gamma_3$ for $t = 0.3$, $\chi_\Psi(r) = r$ on $[0,\f14]$ and $\chi_\Psi = 1$ on $[\f12,\infty[$ ($c_\Psi = 0$, $c'_\Psi = \f12$). The
curve $\Gamma_j$ leaves $j$ in the direction $\arg(z-j) = -2\arctan t$.}\label{fig-Lambda}
\end{figure}

For $0<t<s_1$ and $j\in\bN$ let $\lambda_t(r) := r-it\chi_\Psi(r)$, $\Gamma_j := \{j+\lambda_t(r)^2;\ r\ge0\}$ and $\Lambda_t := \bigcup_j\Gamma_j$.
After decreasing $s_1$ if necessary, each $\Gamma_j$ is the graph $\{x-ih_j(x);\ x\ge j\}$ of a continuous function $h_j\ge0$, $h_j = 0$ on $[j,j+c_\Psi^2]$ and $h_{j+1}<h_j$ on $]j+1,\infty[$. Hence
 $\bC\setminus\Lambda_t$ is connected, and the region between $\Gamma_j$ and
$\Gamma_{j+1}$ can be reached from the upper half-plane across the interval $]j+c_\Psi^2,j+1[$; see Figure \ref{fig-Lambda}.

\begin{prop}\label{prop-essdef} Let $0<t<s_1$.
\begin{itemize}
\item[(i)] $\sigma(\bF_{-it}) = \Lambda_t$.
\item[(ii)] $\sigma_{\rm ess}(\bP_{-it}) = \Lambda_t$, and $(\bP_{-it}-z)^{-1}$ is meromorphic in $\bC\setminus\Lambda_t$, its poles being eigenvalues of $\bP_{-it}$ with
Riesz projections of finite rank.
\end{itemize}
\end{prop}
\pf (i) For $z\notin\Lambda_t$, $\sup_\xi\|(\bF_{-it}(\xi)-z)^{-1}\|<\infty$, by \eqref{subfib} for $|\xi|$
large and by continuity for $|\xi|$ bounded.
The theorem on direct integrals of \cite{w2} (see \cite[Theorem 3.1]{w4}) gives (i).
(ii) As
$(1+H_0)(\bF_{-it}(\xi)-z)^{-1}$ is bounded and $\|(\bF_{-it}(\xi)-z)^{-1}\|\to0$ as $|\xi|\to\infty$, $\bW_{-it}(\bF_{-it}-z)^{-1}$ is compact.
Since $\bC\setminus\Lambda_t$ is
connected and contains points of the resolvent set of $\bP_{-it}$ (Proposition \ref{prop-def}(c)), (ii) follows from \cite[XIII.4, Lemma 3]{rs4}, as
in \cite[Proposition 3.3]{w4}, and from the analytic Fredholm theorem. \ef

The region between $\Gamma_j$ and $\Gamma_{j+1}$ plays the role of the region between the half-lines $j+e^{-2i\varphi}[0,\infty[$ and
$j+1+e^{-2i\varphi}[0,\infty[$ for the dilation (see Figure \ref{fig-Lambda} below). 

By analogy with \eqref{resdefn}, we call
\[
{\rm Res}^{\rm ph}(P) := \bigcup_{0<t<s_1} {\rm Res}^{\rm ph}_t(P)
\]
the resonances of $P$ defined by phase-space deformation, where
\[
{\rm Res}^{\rm ph}_t(P) := \big\{z\in\bC\setminus\Lambda_t;\ \Re z>0,\ -h_0(\Re z)<\Im z\le0,\ z \mbox{ is a pole of } (\bP_{-it}-z)^{-1}\big\}.
\]

Let $0<\varphi<a$, $\Im\theta = \varphi$, and let $\vU$ be the connected component containing the upper half-plane of
$\{-2\varphi<\arg z<\f{3\pi}2\}\setminus\big(\Lambda_t\cup\bigcup_j(j+e^{-2i\varphi}[0,\infty[)\big)$. It contains the rectangles
$\{j+c'^2_\Psi<\Re z<j+1,\ -\min(2c'_\Psi t,c'^2_\Psi\tan2\varphi)<\Im z\le0\}$, $j\in\bN$. In $\vU$, $(\bP_{-it}-z)^{-1}$ and $R(z,\theta)$ are continuations of the
resolvent from the upper half-plane across the same intervals of $\bR_+$.

\begin{prop}\label{prop-equiv} In $\vU$, the poles of $(\bP_{-it}-z)^{-1}$ and those of $R(z,\theta)$ coincide, with multiplicities. In particular, they do not
depend on $t$.
\end{prop}
\pf It follows as in \cite{hm,w4} by considering the matrix elements
\[
G(z) := \w{U_{\bar\theta}g,R(z,\theta)U_\theta f}, \qquad \widetilde G(z) := \w{\vV_{it}\hat g,(\bP_{-it}-z)^{-1}\vV_{-it}\hat f} .
\]
These functions are meromorphic in $\vU$ for $f,g\in\vD_{\rm an}$, by \cite{w4} and Proposition \ref{prop-essdef}(ii). They are equal to $\w{g,(P-z)^{-1}f}$ for $\Re z\ll0$, by \cite{ac,w4} and \eqref{matrixel},
 hence equal in $\vU$. \ef

\subsection{Width of the poles}\label{sec-flux}
In this subsection, $\gamma\ge1$ is again the friction parameter, and $\gamma\to\infty$. We use the phase-space deformation of \S\ref{ssec-phase} at the
frequency scale $\gamma$. Let $0<c_\Psi<c'_\Psi<1$, with $c_\Psi$ close to $1$, and let $\chi_\Psi$ be as in \S\ref{ssec-phase}, now chosen with $\chi_\Psi = 0$ on $[0,c_\Psi]$ and $\chi_\Psi>0$ on $]c_\Psi,\infty[$. Let $t_0\in\,]0,1]$ and $\Sigma_0 := \{s\in\bC;\ |\Re s|<1,\ |\Im s|<2t_0\}$. Put
\[
\begin{aligned}
&\Psi(\xi) := \chi_\Psi(|\xi|/\gamma)\,\f\xi{|\xi|}, \qquad \phi_s(\xi) := \xi+s\Psi(\xi) = e^{-\theta_s(\xi)}\xi, \\
&\theta_s(\xi) := -\log\Big(1+\f{s\chi_\Psi(|\xi|/\gamma)}{|\xi|}\Big), \qquad s\in\Sigma_0 .
\end{aligned}
\]
Since $\|\na\Psi\|_\infty\le C_\Psi\gamma^{-1}$, the conditions imposed on $s$ in \S\ref{ssec-phase} hold for $s\in\Sigma_0$ when $\gamma\ge\gamma_0$.
The definitions and proofs of \S\ref{ssec-phase} apply, with $H_0$ replaced by $\gamma H_0$.
We use them in the following quantitative form.
The constants $c_\Psi<c'_\Psi<1$ and $\delta_2<\delta_1<\delta_V$ will be chosen at the end of the proof of Theorem \ref{th-real}, and $\gamma_0$
depends on them. For $s\in\Sigma_0$ and $\gamma\ge\gamma_0$,
\be \label{phiprop}
\begin{aligned}
&|J_s-1|\le C\gamma^{-1}, \qquad |\theta_s(\xi)|\le6(c_\Psi\gamma)^{-1} , \\
&\big|\phi_s(\xi)-\phi_s(\xi')-(\xi-\xi')\big|\le3C_\Psi\gamma^{-1}|\xi-\xi'| ,
\end{aligned}
\ee
\be \label{Fdilg}
\begin{aligned}
&\bF_s(\xi) := \gamma H_0+i\phi_s(\xi)\cdot v = \widehat P_0\big(\xi,\gamma,\theta_s(\xi)\big), \\
&E^{\theta_s(\xi)}_j(\xi,\gamma) = \gamma j+\gamma^{-1}\phi_s(\xi)\cdot\phi_s(\xi) ,
\end{aligned}
\ee
so that Theorem \ref{th2.1}, \eqref{e4.26} and Lemma \ref{lem2.3} apply to the deformed fibres, uniformly in $\xi$, with $\zeta = \gamma^{-1}\phi_s(\xi)$. The
operators $\bW_s$, $\bP_s = \bF_s+\bW_s$ and $\vV_s$ are defined by \eqref{Ps} and the same formulas. By \eqref{hatVbound} and \eqref{phiprop},
\be \label{kbound}
|K_s(\xi,\xi')|\le C\,\kappa_{\delta_2}(\xi-\xi'), \qquad K_s(\xi,\xi') = K_0(\xi,\xi') \ \mbox{ if } |\xi|,|\xi'|\le c_\Psi\gamma .
\ee
Proposition \ref{prop-def} holds, with $P(\gamma)$ in place of $P$, and with
\be \label{defg}
\begin{aligned}
&\|(\bW_s-\bW_{\Re s})(1+H_0)^{-1/2}\|\le C|\Im s|\gamma^{-1} \quad \mbox{in (a)}, \\
&\big|\Re\w{w,\bP_sw}-\gamma\w{w,H_0w}\big|\le C|\Im s|\,\|(1+H_0)^{1/2}w\|^2 \quad \mbox{in (b)}.
\end{aligned}
\ee
The curves $\Gamma_j$ become $\{\gamma j+\gamma^{-1}(r-it\chi_\Psi(r/\gamma))^2;\ r\ge0\}$. Proposition \ref{prop-equiv} holds in a region $\vU$ which contains, for $\gamma$
large depending on $\varphi$, the rectangles
\be \label{rectangle}
\big\{z\in\bC;\ \gamma j+c'^2_\Psi\gamma<\Re z<\gamma(j+1),\ -2c'_\Psi t<\Im z\le0\big\}, \qquad j\in\bN .
\ee

Let $\mathfrak b(u,u') := \iint u\,u'\,dx\,dv$ be the bilinear form on $\vH$. In Fourier variables,
\be \label{bform}
\mathfrak b(u,u') = \iint\hat u(-\xi,v)\,\hat u'(\xi,v)\,dv\,\f{d\xi}{(2\pi)^n} .
\ee
We write $\mathfrak b(w,w')$ for the right-hand side when $w,w'\in\widehat\vH$. The identity ${}^tP(\gamma) = \bJ P(\gamma)\bJ$ following \eqref{bJ} reads
$\mathfrak b(u,P(\gamma)u') = \mathfrak b(\bJ P(\gamma)\bJ u,u')$.
The complex conjugation $u\mapsto\bar u$, which commutes with $P(\gamma)$, becomes the antiunitary
operator $(\vC w)(\xi) := \overline{w(-\xi)}$.

\begin{prop}\label{prop-sym} For $s\in\Sigma_0$ and $w,w'\in\vD_P$, one has

\textup{(i)} $\vC\bP_s\vC = \bP_{\bar s}$;

\textup{(ii)} $\mathfrak b(w,\bP_sw') = \mathfrak b(\bJ\bP_s\bJ w,w')$.
\end{prop}
\pf (i) $\phi_s$ is odd, $J_s$ is even, $\overline{\phi_s} = \phi_{\bar s}$, so that $\overline{\bF_s(-\xi)} = \bF_{\bar s}(\xi)$, and $\overline{K_s(-\xi,-\xi')} = K_{\bar s}(\xi,\xi')$ by Lemma
\ref{lem-hatV}. (ii) The transpose in $v$ of $\bF_s(\xi)$ is $\bF_s(\xi) = \bJ\bF_s(-\xi)\bJ$, that of $\p_{v_l}$ is $-\p_{v_l} = \bJ\p_{v_l}\bJ$, and $K_s(-\xi',-\xi) = K_s(\xi,\xi')$.
The
change of variables $(\xi,\xi')\mapsto(-\xi',-\xi)$ gives the claim on $\vS(\bR^{2n})$, hence on $\vD_P$. \ef

 For $\zeta\in\bC^n$ let $e_\alpha(\zeta) := \psi_\alpha(\cdot+2i\zeta)$, so that
$
\Pi^\zeta_j = \sum_{|\alpha| = j}\w{\overline{e_\alpha(\zeta)},\cdot}\,e_\alpha(\zeta)
$
 by \eqref{Riesz}.
A contour shift gives, for $|\alpha| = |\beta|$,
\be \label{nodiag}
\w{\overline{e_\alpha(\zeta)},\p_{v_l}e_\beta(\zeta)} = \w{\psi_\alpha,\p_{v_l}\psi_\beta} = 0 .
\ee
Thus $\Pi^\zeta_j\p_{v_l}\Pi^\zeta_j = 0$: the coupling does not couple a channel to itself at first order. If $|\theta|\le\f12$ and
$(\widehat P_0(\xi,\gamma,\theta)-z)u = f$ with $|z|\le M\gamma$, the real part of $\w{u,e^\theta(\widehat P_0(\xi,\gamma,\theta)-z)u}$ gives, as in the proof of Lemma
\ref{lem2.3},
\be \label{energy}
\|(1+H_0)^{\f12}u\|\le C_M\big(\|u\|+\gamma^{-\f12}\|u\|^{\f12}\|f\|^{\f12}\big) .
\ee

By Corollary
\ref{cor-intpow}, the poles $z_l(\gamma)$ produced by a negative eigenvalue $\eta_0$ of $-\Delta^{(k)}_V$ lie just below the threshold $\gamma k$. They are at distance $|\eta_0|\gamma^{-1}+ O(\gamma^{-2})$
from it, in the
direction $\arg(z-\gamma k) = \pi$. Since $H_0$ is even in $v$ while the skew-adjoint operator $L$ is
odd in $v$, the parity operator
$\bJ$ of \eqref{bJ} satisfies
\be \label{JPJ}
\bJ\, P(\gamma)\, \bJ = P(\gamma)^* , \qquad \bJ\, P_0(\gamma)\, \bJ = P_0(\gamma)^* .
\ee
If $\Im z_l(\gamma) >0$, then $z_l(\gamma)$ is an eigenvalue of $P(\gamma)$. By the reality of $P(\gamma)$, that is,
$P(\gamma)\bar u = \overline{P(\gamma)u}$, $\overline{z_l(\gamma)}$ is also an eigenvalue of $P(\gamma)$, lying in
the lower half-plane. It is a pole of $R(z,\gamma,\bar\theta)$, but in general not of $R(z,\gamma,\theta)$.
With $\ep_0$ as in Theorem \ref{th-bij}, put
\[
\bD_\gamma := \{z\in\bC;\ |z-\gamma k-\eta_0\gamma^{-1}|<\ep_0\gamma^{-1}\} ,
\]
and let $\delta_V$ be
defined by \eqref{deltaV}.

\begin{theorem}\label{th-real}
Let $\eta_0<0$ be an eigenvalue of $-\Delta^{(k)}_V$ of multiplicity $k_0$
and let
$z_1(\gamma),\allowbreak\dots,z_{k_0}(\gamma)$ be the poles of $R(z,\gamma,\theta)$ given by Theorem
\ref{th-bij} in $\bD_\gamma$. Then
the $z_l(\gamma)$ are independent of $\theta$ for $0<\varphi<a$. Moreover, for every $c<2\delta_V$, there is $C_c$ such
that, for $\gamma$ large,
\be \label{expsmall}
|\Im\, z_l(\gamma)| \le C_c\,e^{-c\gamma}, \qquad 1\le l\le k_0 .
\ee
In particular, the coefficients of the asymptotic expansions \eqref{expansion} of Theorem \ref{th3.10} are real.
\end{theorem}


\subsubsection{Proof of Theorem \ref{th-real}: the deformation for large \texorpdfstring{$\gamma$}{gamma}}\label{sec-proofreal}

Throughout the proof, \eqref{ass1b} holds, and $\eta_0<0$, $k_0$, $\ep_0$ and $\bD_\gamma$ are as in
Theorem \ref{th-real}. By Proposition \ref{prop-det}(iii), $\ep_0$ may be chosen independent of $\theta$.
By \eqref{expansion}, we may
decrease it without changing the poles in $\bD_\gamma$ for $\gamma$ large. We assume $4\ep_0\le|\eta_0|$. We write $\eta := \gamma(z-\gamma k)$, so that
$z\in\bD_\gamma$ if and only if $|\eta-\eta_0|<\ep_0$, and then $\Re\eta\le-\f34|\eta_0|$ and $|z|\le(k+1)\gamma$. We fix $\sigma\in\,]\f34,1[$.

Let $\vU$ be the region of Proposition \ref{prop-equiv} at the scale $\gamma$. For $\gamma$ large, $\bD_\gamma$ is contained in the union of the upper
half-plane and of the rectangle \eqref{rectangle} with $j = k-1$, hence in $\vU$. By Proposition \ref{prop-equiv} and Proposition \ref{prop-sym}(i),
 for $0<t\le t_0$ and $\gamma$ large, depending on $t$ and $\varphi$, every pole $z$ of $R(z,\gamma,\theta)$ in $\bD_\gamma$ is an eigenvalue of $\bP_{-it}$, of the same
multiplicity. Moreover, if $\bP_{-it}w = zw$, then $\bP_{it}\vC w = \bar z\,\vC w$. In particular, the poles of $R(z,\gamma,\theta)$ in $\bD_\gamma$
are the eigenvalues of $\bP_{-it}$ there. Hence they do not depend on $\theta$, for $0<\varphi<a$ and $\gamma$ large depending on $\varphi$.

\subsubsection{The deformed free fibres}
Let $\chi_0+\chi_1 = 1$ be the partition of unity \eqref{pu} for which \eqref{e4.26} holds with $M = k+2$ and Lemma \ref{lem2.3} holds with $\ell = k$ and
$b = k+\f12$. Enlarging $R$, we assume $R\ge\sqrt k+1$. Let $R'\ge R$ be such that $\supp\chi_0\subset\{|\xi|\le R'\}$, and let $\gamma_0$ be so large that
$(R'+1)\gamma^\sigma<c_\Psi\gamma$ for $\gamma\ge\gamma_0$. For $j<k$, let
$r_j := \gamma\sqrt{k-j}$, $c_N := \min\{\f1{8\sqrt k},\f12(1-c'_\Psi)\}$,
\[
\vN_j := \big\{\xi\in\bR^n;\ \big||\xi|-r_j\big|\le c_N\gamma\big\}, \qquad \vartheta_j(\xi) := \vartheta\Big(\f{|\xi|-r_j}{c_N\gamma}\Big) ,
\]
where $\vartheta\in C^\infty_0(]-1,1[;[0,1])$ is equal to $1$ on $[-\f12,\f12]$. The shells $\vN_j$ are disjoint, at mutual distance $\ge\gamma/(4\sqrt k)$, and contained in
$\{c'_\Psi\gamma\le|\xi|\le R\gamma\}$, where $\Psi(\xi) = \xi/|\xi|$ and $\chi_0(\gamma^{-1}\xi) = 1$. For $0<t\le t_0$ and $z\in\bD_\gamma$ put
\begin{align*}
& \xi_t := \phi_{-it}(\xi) = \xi-it\Psi(\xi),  \qquad \zeta_t(\xi) := \gamma^{-1}\xi_t, \\
&R_t(\xi) := \big(\bF_{-it}(\xi)-z\big)^{-1}, \qquad D_j(\xi) := \gamma j+\gamma^{-1}\xi_t\cdot\xi_t-z, \\
&R^{\rm reg}(\xi) := R_t(\xi)-\sum_{j<k}\f{\vartheta_j(\xi)}{D_j(\xi)}\,\Pi^{\zeta_t(\xi)}_j ,
\end{align*}
so that, by \eqref{Fdilg}, $D_j(\xi) = E^{\theta_{-it}(\xi)}_j(\xi,\gamma)-z$ and $R_t(\xi)\Pi^{\zeta_t(\xi)}_j = D_j(\xi)^{-1}\Pi^{\zeta_t(\xi)}_j$.
\begin{lemma}\label{lem-fibres} There is $\gamma_0$ such that, for $0<t\le t_0$, there is $\gamma_t\ge\gamma_0$ with the following properties.
For $\gamma\ge\gamma_t$, $z\in\bD_\gamma$ and $\xi\in\bR^n$, $z\notin\sigma(\bF_{-it}(\xi))$.
Moreover, with constants independent of $t$, $\xi$, $z$ and $\gamma$:
\begin{itemize}
\item[(i)] $\|(1+H_0)^{\f12}R_t(\xi)\|\le C(\gamma^{-1}+\gamma\w\xi^{-2})$ for $|\xi|\le\gamma^\sigma$;
\item[(ii)] $|D_j(\xi)|\ge\f14\big(\big||\xi|-r_j\big|+t\big)$ and $\|(1+H_0)^{\f12}\Pi^{\zeta_t(\xi)}_j\|\le C$ for $\xi\in \vN_j$, $j<k$;
\item[(iii)] $\|(1+H_0)^{\f12}R^{\rm reg}(\xi)\|\le C\gamma^{1-2\sigma}$ and $\|(1+H_0)^{\f12}R_t(\xi)\|\le Ct^{-1}$ for $|\xi|\ge\gamma^\sigma$.
\end{itemize}
\end{lemma}
\pf Write
\[
E_j := E^{\theta_{-it}(\xi)}_j(\xi,\gamma) = \gamma j+\gamma^{-1}\xi_t\cdot\xi_t.
\]
 One has $\xi_t\cdot\xi_t = (|\xi|-it\chi_\Psi(|\xi|/\gamma))^2$, so that
$\Re\,\xi_t\cdot\xi_t\ge|\xi|^2-1$ and
\[
\Im\,\xi_t\cdot\xi_t = -2t\chi_\Psi(|\xi|/\gamma)|\xi|\le0.
\]
 We first show, for $\gamma\ge\gamma_0$ and a constant $c_F>0$ independent of $t$:
\begin{itemize}
\item[(F1)] $|E_j-z|\ge c_F\gamma$ if $j>k$, or if $j<k$ and $\big||\xi|-r_j\big|\ge\f12c_N\gamma$;
\item[(F2)] $|E_k-z|\ge c_F\gamma^{-1}(|\xi|^2+1)$.
\end{itemize}
If $j>k$, $\Re(E_j-z)\ge\gamma-\gamma^{-1}(1+|\eta|)$. If $j<k$ and $|\xi|\le c_\Psi\gamma$, then $\chi_\Psi = 0$ and $\Re(E_j-z)\le-(1-c_\Psi^2)\gamma+C\gamma^{-1}$. If $j<k$,
$|\xi|\ge c_\Psi\gamma$ and $\big||\xi|-r_j\big|\ge\f12c_N\gamma$, then $$|\Re(E_j-z)|\ge\gamma^{-1}\big||\xi|-r_j\big|(|\xi|+r_j)-C\gamma^{-1}\ge\f12c_N\gamma-C\gamma^{-1}.$$ This gives (F1). For
(F2), note that $\gamma(E_k-z) = \xi_t\cdot\xi_t-\eta$. If $\chi_\Psi(|\xi|/\gamma) = 0$, this is $|\xi|^2-\eta$, with $\Re(-\eta)\ge\f34|\eta_0|$. Otherwise, $|\xi|\ge c_\Psi\gamma$ and
$\Re(\xi_t\cdot\xi_t-\eta)\ge|\xi|^2-1$.

(ii) On $\vN_j$, $\chi_\Psi = 1$ and
\[
D_j(\xi) = \f{(|\xi|-r_j)(|\xi|+r_j)}\gamma-\f{t^2+\Re\eta}\gamma-\f i\gamma\big(2t|\xi|+\Im\eta\big) ,
\]
with $(|\xi|+r_j)/\gamma\ge1$, $2t|\xi|/\gamma\ge t$ and $|\eta|\le2|\eta_0|$. Hence $|\Im D_j|\ge\f t2$ and $|\Re D_j|\ge\big||\xi|-r_j\big|-C\gamma^{-1}$. These bounds give the
lower bound for $\gamma\ge\gamma_t := \max(\gamma_0,4C/t)$, by distinguishing the cases $\big||\xi|-r_j\big|\ge t$ and $\big||\xi|-r_j\big|<t$. The range of $\Pi^\zeta_j$ is spanned by the
$e_\alpha(\zeta)$, $|\alpha| = j$, so that
\[
\|(1+H_0)^{1/2}\Pi^\zeta_j\|\le\sum_{|\alpha| = j}\|(1+H_0)^{1/2}e_\alpha(\zeta)\|\,\|e_\alpha(\zeta)\|,
\]
which is bounded for $\zeta$ in the compact set
\be \label{Xi}
\Xi := \{\zeta\in\bC^n;\ |\Re\zeta|\le R',\ |\Im\zeta|\le1\} .
\ee
This set contains $\zeta_t(\xi)$ for $|\xi|\le R'\gamma$. The same holds for $\|\Pi^\zeta_j\|$.

By \eqref{Fdilg} and $|\theta_{-it}(\xi)|\le6(c_\Psi\gamma)^{-1}$, one has, for $\gamma$ large,
\[
|z|\le\gamma k+C\gamma^{-1}<b\cos(\Im\theta_{-it}(\xi))\gamma.
\]
Hence Lemma \ref{lem2.3} applies and gives
\be \label{R0def}
R_t(\xi) = \sum_{j=0}^k\chi_0(\gamma^{-1}\xi)\,\f{\Pi^{\zeta_t(\xi)}_j}{E_j-z}+r_k(\xi), \qquad \|r_k(\xi)\|\le C\gamma^{-1} ,
\ee
as soon as $z\neq E_j$ for all $j$. This holds by (F1), (F2) and (ii). By \eqref{energy}, with $\theta = \theta_{-it}(\xi)$, an $L^2$-bound $\|R_t(\xi)\|\le a$ with
$a\ge\gamma^{-1}$ implies $\|(1+H_0)^{1/2}R_t(\xi)\|\le Ca$. The same holds for $R_t(\xi)$ composed on the right with a bounded operator commuting with
$\bF_{-it}(\xi)$.

(i) For $|\xi|\le\gamma^\sigma$, \eqref{R0def}, (F1), (F2) and the boundedness of $\|\Pi^\zeta_j\|$ on $\Xi$ give $\|R_t(\xi)\|\le C(\gamma\w\xi^{-2}+\gamma^{-1})$, hence (i).

(iii) Let $|\xi|\ge\gamma^\sigma$. Then $|E_k-z|\ge c_F\gamma^{2\sigma-1}$ by (F2). Moreover, $\chi_0(\gamma^{-1}\xi)-\vartheta_j(\xi)$ vanishes unless $\big||\xi|-r_j\big|\ge\f12c_N\gamma$, and then
$|E_j-z|\ge c_F\gamma$ by (F1). By \eqref{R0def},
\[
\|R^{\rm reg}(\xi)\|\le C(\gamma^{1-2\sigma}+\gamma^{-1})\le C\gamma^{1-2\sigma}.
\]
 Since

\[
R^{\rm reg}(\xi) = R_t(\xi)\big(1-\sum_{j<k}\vartheta_j(\xi)\Pi^{\zeta_t(\xi)}_j\big),
\]
 where the last factor is bounded and commutes with $\bF_{-it}(\xi)$, the first bound of (iii)
follows. The second one then follows with (ii). \ef

\subsubsection{The resonant states at high frequencies}

In this paragraph, $0<t\le t_0$ is fixed, $\gamma\ge\gamma_t$, $z$ is a pole of $R(z,\gamma,\theta)$ in $\bD_\gamma$, and $w\in\vD_P$ is an eigenfunction of $\bP_{-it}$
for $z$ with $\|w\| = 1$ (see the beginning of \S\ref{sec-proofreal}). Let $\chi^\sigma_j(\xi) := \chi_j(\gamma^{-\sigma}\xi)$, as in \eqref{chisig}. For $0<\delta_3<\delta_2$ and
$N\in\bN\cup\{\infty\}$, let
\[
\omega_N(\xi) := \exp\Big(\delta_3\min\big\{(\w\xi-R''\gamma^\sigma)_+,\,N\big\}\Big), \qquad R'' := R'+1 .
\]
Then $\omega_N = 1$ on $\supp\chi^\sigma_0$, $|\log\omega_N(\xi)-\log\omega_N(\xi')|\le\delta_3|\xi-\xi'|$, $\omega_N$ is even, and it is bounded with bounded inverse for
$N<\infty$. Thus $\omega_N = e^{\delta_3\min\{d,N\}}$, where $d(\xi) := (\w\xi-R''\gamma^\sigma)_+$ is a $1$-Lipschitz function, comparable to the distance from $\xi$ to the low-frequency region $\{|\xi|\le R''\gamma^\sigma\}$. Thus $d$ is the Agmon distance in phase space, for the metric $\delta_3|d\xi|$. This metric is the relevant one because the coupling $\bW_{-it}$ has a kernel decaying like $e^{-\delta_2|\xi-\xi'|}$, with $\delta_2>\delta_3$ (Lemma \ref{lem-chan}). Proposition \ref{prop-hf} below is the corresponding Agmon estimate: $w$ decays like $e^{-\delta_3d}$ away from the low-frequency region. These multiplications commute with $\bF_{-it}$, $R_t$, $H_0$, $\Pi_j$ and $\bJ$. Put $\delta' := \delta_2-\delta_3$ and, for $N<\infty$,
\[
A := \chi^\sigma_1R_t, \qquad \bW^\omega := \omega_N\bW_{-it}\omega_N^{-1}, \qquad Y := \bW^\omega A .
\]

\begin{lemma}\label{lem-chan} Uniformly in $N<\infty$,
\begin{align*}
&\|(1+H_0)^{\f12}A\|\le Ct^{-1}, \qquad \|\bW^\omega(1+H_0)^{-\f12}\|\le C, \\
&\|Y\|\le Ct^{-1}, \qquad \|Y^2\|\le Ct^{-2}\gamma^{-\f12} .
\end{align*}
\end{lemma}
\pf The first bound is Lemma \ref{lem-fibres}(iii), since $\supp\chi^\sigma_1\subset\{|\xi|\ge\gamma^\sigma\}$.
The kernel of $\bW^\omega$ is $\omega_N(\xi)K_{-it}(\xi,\xi')\omega_N(\xi')^{-1}\cdot\na_v$,
bounded by $C\kappa_{\delta'}(\xi-\xi')$ by \eqref{kbound}. The Schur test gives the second bound, and the third follows. Write
$Y = \sum_{j<k}Y_j+Y_{\rm r}$ with $Y_j := \bW^\omega\chi^\sigma_1\vartheta_jD_j^{-1}\Pi^{\zeta_t}_j$ and $Y_{\rm r} := \bW^\omega\chi^\sigma_1R^{\rm reg}$. By Lemma \ref{lem-fibres},
$\|Y_j\|\le Ct^{-1}$ and
\[
\|Y_{\rm r}\|\le C\gamma^{1-2\sigma}\le C\gamma^{-1/2},
\]
 so that $\|Y^2-\sum_{j,j'<k}Y_jY_{j'}\|\le Ct^{-1}\gamma^{-1/2}$. Now $Y_jY_{j'} = \bW^\omega X$,
where, with $\zeta = \zeta_t(\xi)$, $\zeta' = \zeta_t(\xi')$,
\[
X(\xi,\xi') = \sum_{|\alpha| = j,\,|\beta| = j'}X_{\alpha\beta}(\xi,\xi')\,e_\alpha(\zeta)\w{\overline{e_\beta(\zeta')},\cdot} ,
\]
\[
\begin{aligned}
X_{\alpha\beta}(\xi,\xi') = {}&\f{(\chi^\sigma_1\vartheta_j)(\xi)\,\omega_N(\xi)}{D_j(\xi)}\;\sum_{l=1}^nK^l_{-it}(\xi,\xi')\,\w{\overline{e_\alpha(\zeta)},\p_{v_l}e_\beta(\zeta')} \\
&\times\f{(\chi^\sigma_1\vartheta_{j'})(\xi')}{\omega_N(\xi')\,D_{j'}(\xi')} .
\end{aligned}
\]
Since the norms of $e_\alpha(\zeta)$ and of $(1+H_0)^{1/2}e_\alpha(\zeta)$ are bounded on $\Xi$,
\[
\|\bW^\omega X\|\le C\max_{\alpha,\beta}\|X_{\alpha\beta}\|_{\rm Schur} .
\]
Here
\[
\|K\|_{\rm Schur} := \max\big\{\sup_\xi\int|K(\xi,\xi')|\,d\xi',\ \sup_{\xi'}\int|K(\xi,\xi')|\,d\xi\big\}
\]
is the Schur norm of a scalar kernel
$K$. It bounds the norm on $L^2(\bR^n)$ of the integral operator with kernel $K$ (Schur test). If $j\neq j'$, then
$|\xi-\xi'|\ge\gamma/(4\sqrt k)$ on $\vN_j\times \vN_{j'}$, and Lemma \ref{lem-fibres}(ii) gives
\[
\|X_{\alpha\beta}\|_{\rm Schur}\le Ct^{-2}\int_{|\zeta|\ge\gamma/(4\sqrt k)}\kappa_{\delta'}(\zeta)\,d\zeta\le Ct^{-2}e^{-\delta'\gamma/(8\sqrt k)} .
\]
If $j = j'$, \eqref{nodiag} and the mean value theorem on
the convex set $\Xi$ give
\[
\begin{aligned}
\big|\w{\overline{e_\alpha(\zeta)},\p_{v_l}e_\beta(\zeta')}\big| &= \big|\w{\overline{e_\alpha(\zeta)},\p_{v_l}(e_\beta(\zeta')-e_\beta(\zeta))}\big| \\
&\le C|\zeta-\zeta'|\le C\gamma^{-1}|\xi-\xi'| ,
\end{aligned}
\]
so that
\[
|X_{\alpha\beta}(\xi,\xi')|\le Ct^{-2}\gamma^{-1}|\xi-\xi'|\,\kappa_{\delta'}(\xi-\xi')
\]
Hence $\|X_{\alpha\beta}\|_{\rm Schur}\le Ct^{-2}\gamma^{-1}$, since $|\zeta|\kappa_{\delta'}(\zeta)$ is integrable on $\bR^n$.
Altogether,
\[
\|Y^2\|\le C(t^{-1}\gamma^{-1/2}+t^{-2}\gamma^{-1}+t^{-2}e^{-\delta'\gamma/(8\sqrt k)})\le Ct^{-2}\gamma^{-1/2}.
\]
 \ef

Recall that $w\in\vD_P\subset\widehat\vH$ is the normalized eigenfunction of $\bP_{-it}$ for $z$ fixed above. The weight $\omega_N$ with $N = \infty$,
\be \label{omegainfty}
\omega_\infty(\xi) = \exp\big(\delta_3(\w\xi-R''\gamma^\sigma)_+\big),
\ee
is unbounded. The second bound
below means that $\omega_\infty w\in D(H_0^{1/2})$, with the stated bound.

\begin{prop}\label{prop-hf} One has
\[
\|(1+H_0)^{\f12}w\|\le C, \qquad \|(1+H_0)^{\f12}\omega_\infty w\|\le Ct^{-2}, \qquad \||\xi|\,w\|\le Ct^{-2}\gamma^\sigma .
\]
\end{prop}
\pf Since $\Re\w{w,(\bP_{-it}-z)w} = 0$, \eqref{defg} gives
\[
(\gamma-Ct)\w{w,H_0w}\le\Re z+Ct\le\gamma k+C,
\]
 hence $\w{w,H_0w}\le k+1$ for $\gamma$ large.
Next, $(\bF_{-it}-z)w = -\bW_{-it}w$, so that $\chi^\sigma_1w = -A\bW_{-it}(\chi^\sigma_1w+\chi^\sigma_0w)$. Multiplying by $\omega_N$, $N<\infty$, and using $\omega_N\chi^\sigma_0 = \chi^\sigma_0$, we get
\[
(1+A\bW^\omega)x = y, \qquad x := \omega_N\chi^\sigma_1w, \quad y := -A\bW^\omega\chi^\sigma_0w ,
\]
with $x,y\in D(H_0^{1/2})$. Applying $\bW^\omega$ gives $(1+Y)\bW^\omega x = \bW^\omega y$. By Lemma \ref{lem-chan}, $1+Y$ is invertible for $\gamma$ large, with $(1+Y)^{-1} = (1-Y)(1-Y^2)^{-1}$ and $\|(1+Y)^{-1}\|\le Ct^{-1}$. Hence
\[
x = y-A\bW^\omega x = y-A(1+Y)^{-1}\bW^\omega y = -A(1+Y)^{-1}\bW^\omega\chi^\sigma_0w,
\]
 and
\[
\|(1+H_0)^{\f12}\omega_N\chi^\sigma_1w\|\le\|(1+H_0)^{\f12}A\|\,\|(1+Y)^{-1}\|\,\|\bW^\omega(1+H_0)^{-\f12}\|\,\|(1+H_0)^{\f12}w\|\le Ct^{-2} ,
\]
uniformly in $N$. Letting $N\to\infty$ and using $\omega_\infty\chi^\sigma_0 = \chi^\sigma_0$ gives the second bound. Finally, $\chi^\sigma_1 = 1$ on $\{|\xi|\ge R'\gamma^\sigma\}$, so that
\[
\||\xi|\,w\|\le2R''\gamma^\sigma+\sup_{|\xi|\ge2R''\gamma^\sigma}|\xi|e^{-\delta_3(|\xi|-R''\gamma^\sigma)}\,\|\omega_\infty\chi^\sigma_1w\|\le Ct^{-2}\gamma^\sigma .
\]
\ef

\begin{lemma}\label{lem-levelk} One has $\|(1-\Pi_k)w\|\le Ct^{-2}\gamma^{\sigma-1}$ and
\be \label{Jww}
\big|\w{\bJ w, w} - (-1)^k\big| \le 2\|(1-\Pi_k)w\|^2 \le Ct^{-4}\gamma^{2\sigma-2} .
\ee
In particular, $|\w{\bJ w,w}|\ge\f12$ for $\gamma$ large.
\end{lemma}
\pf Let $\vI_+ := \{k+1,k+2,\dots\}$, $\vI_- := \{0,\dots,k-1\}$, $\Pi_{\vI_\pm} := \sum_{j\in \vI_\pm}\Pi_j$ and $P^{\vI_\pm} := \Pi_{\vI_\pm}\bP_{-it}\Pi_{\vI_\pm}$. Since $v_l$ and $\p_{v_l}$ map
$\Ran\Pi_j$ into $\Ran\Pi_{j-1}\oplus\Ran\Pi_{j+1}$, one has $\Pi_{\vI_\pm}\bP_{-it}\Pi_{\vI_\mp} = 0$. Hence $\Pi_{\vI_\pm}(\bP_{-it}-z)w = 0$ gives
$(P^{\vI_\pm}-z)\Pi_{\vI_\pm}w = -\Pi_{\vI_\pm}\bP_{-it}\Pi_kw$. Here $\Pi_{\vI_\pm}w\in\vD_P$, since $w\in D(H_0)$, $|\xi|w\in\widehat\vH$ by Proposition \ref{prop-hf}, and
$[v_l,\Pi_{\vI_\pm}]$ is bounded. For $f\in\Ran\Pi_{\vI_+}\cap\vD_P$, \eqref{defg} and $H_0\ge k+1$ on $\Ran\Pi_{\vI_+}$ give

\[
\Re\w{f,(P^{\vI_+}-z)f}\ge\big((\gamma-Ct)(k+1)-Ct-\Re z\big)\|f\|^2\ge\f\gamma2\|f\|^2.
\]
For $f\in\Ran\Pi_{\vI_-}\cap\vD_P$, $H_0\le k-1$ and the upper bound in (b) give
$\Re\w{f,(z-P^{\vI_-})f}\ge\f\gamma2\|f\|^2$. Hence $\|f\|\le2\gamma^{-1}\|(P^{\vI_\pm}-z)f\|$ for such $f$. Since
$\Pi_{\vI_\pm}\bP_{-it}\Pi_k = \Pi_{\vI_\pm}(i\xi_t\cdot v+\bW_{-it})\Pi_k$, $\|v_l\Pi_k\|\le C$, $\|\bW_{-it}(1+H_0)^{-1/2}\|\le C$ and $\Pi_k$ commutes with $\xi$,
\[
\|(1-\Pi_k)w\|\le2\gamma^{-1}\big(C\||\xi|\,w\|+C\big)\le Ct^{-2}\gamma^{\sigma-1}
\]
by Proposition \ref{prop-hf}. Finally, $\bJ$ is unitary, commutes with the $\Pi_j$ and equals $(-1)^k$ on $\Ran\Pi_k$. Hence
$\w{\bJ w,w} = (-1)^k\|\Pi_kw\|^2+\w{\bJ(1-\Pi_k)w,(1-\Pi_k)w}$. Since $\|\Pi_kw\|^2 = 1-\|(1-\Pi_k)w\|^2$, this gives \eqref{Jww}. \ef

\subsubsection{The flux identity and the end of the proof}

\begin{lemma}[flux identity]\label{lem-flux} If $w\in\vD_P$ and $\bP_{-it}w = zw$, then
\be \label{fluxt}
2i\,\Im z\;\w{\bJ w,w} = \mathfrak b\big(\bJ(\bP_{-it}-\bP_{it})\vC w,\,w\big) .
\ee
\end{lemma}
\pf By \eqref{bform},
$$\mathfrak b(\bJ\vC w,w) = \iint\overline{w(\xi,-v)}\,w(\xi,v)\,dv\,(2\pi)^{-n}d\xi = \w{\bJ w,w}.$$ By Proposition \ref{prop-sym}(ii),

\[
z\,\mathfrak b(\bJ\vC w,w) = \mathfrak b(\bJ\vC w,\bP_{-it}w) = \mathfrak b(\bJ\bP_{-it}\vC w,w).
\]
Moreover, by Proposition \ref{prop-sym}(i), $\bP_{it}\vC w = \bar z\vC w$, so that

\[
\mathfrak b(\bJ\bP_{-it}\vC w,w) = \bar z\,\mathfrak b(\bJ\vC w,w)+\mathfrak b(\bJ(\bP_{-it}-\bP_{it})\vC w,w).
\]
 \ef

\noindent The identity \eqref{fluxt} is the analogue, for the deformed operators, of the flux identity \eqref{flux} in the proof of Proposition \ref{prop-flux}. The parity $\bJ$ and the
transposition ${}^tP(\gamma) = \bJ P(\gamma)\bJ$ replace the adjoint, and the reality of $P(\gamma)$ relates $\bP_{-it}$ to $\bP_{it}$. Its right-hand side involves
only the difference of the two deformed operators. Since $\phi_{-it}-\phi_{it} = -2it\Psi$, this difference is
\be \label{Pdiff}
\bP_{-it}-\bP_{it} = 2t\,\Psi(\xi)\cdot v+(K_{-it}-K_{it})\cdot\na_v .
\ee
Here $\Psi$ is supported in $\{|\xi|\ge c_\Psi\gamma\}$, and the kernel $K_{-it}-K_{it}$ vanishes on $\{|\xi|\le c_\Psi\gamma,\ |\xi'|\le c_\Psi\gamma\}$ by \eqref{kbound}.
Hence only the high frequencies, where $w$ is exponentially small by Proposition \ref{prop-hf}, contribute. This is quantified by the following lemma.

\begin{lemma}\label{lem-diff} For $u,u'$ with $(1+H_0)^{1/2}\omega_\infty u$, $(1+H_0)^{1/2}\omega_\infty u'\in\widehat\vH$,
\[
\begin{aligned}
&\big|\mathfrak b\big(u,(\bP_{-it}-\bP_{it})u'\big)\big| \\
&\qquad\le C\,e^{-2\delta_3(c_\Psi\gamma-R''\gamma^\sigma)}\,\|(1+H_0)^{\f12}\omega_\infty u\|\,\|(1+H_0)^{\f12}\omega_\infty u'\| .
\end{aligned}
\]
\end{lemma}
\pf We use \eqref{Pdiff}. For the first term, one has $|\int f\,(\vec e\cdot v)g\,dv|\le C\|(1+H_0)^{1/2}f\|\,\|(1+H_0)^{1/2}g\|$ for
$|\vec e|\le1$, since $H_0$ is real. Moreover, $\omega_\infty(\xi)^{-2}\le e^{-2\delta_3(c_\Psi\gamma-R''\gamma^\sigma)}$ on $\supp\Psi$, and $\omega_\infty$ is even. For the second term, if
$\max(|\xi|,|\xi'|)\ge c_\Psi\gamma$, then
\[
\begin{aligned}
\log\omega_\infty(-\xi)+\log\omega_\infty(\xi')+\delta_3|\xi-\xi'|&\ge\delta_3\big(|\xi|+|\xi'|+|\xi-\xi'|\big)-2\delta_3R''\gamma^\sigma \\
&\ge2\delta_3c_\Psi\gamma-2\delta_3R''\gamma^\sigma ,
\end{aligned}
\]
so that
\[
\omega_\infty(-\xi)^{-1}|K_{-it}-K_{it}|(\xi,\xi')\,\omega_\infty(\xi')^{-1}\le Ce^{-2\delta_3(c_\Psi\gamma-R''\gamma^\sigma)}\kappa_{\delta'}(\xi-\xi')
\]
 by \eqref{kbound}, and the Schur test
concludes. \ef

\noindent{\bf End of the proof of Theorem \ref{th-real}.} Let $c<2\delta_V$. Choose $\delta_3<\delta_2<\delta_1<\delta_V$ and $c_\Psi<1$ such that  $c<2\delta_3c_\Psi < 2\delta_V$, then $c'_\Psi\in\,]c_\Psi,1[$, and fix $t = t_0$.
 Let $z$ be a pole of $R(z,\gamma,\theta)$ in $\bD_\gamma$, and let $w$ be an eigenfunction of $\bP_{-it}$ for $z$ with $\|w\| = 1$.
Since $\bJ v_l = -v_l\bJ$ and $\bJ\p_{v_l} = -\p_{v_l}\bJ$, one has $\bJ(\bP_{-it}-\bP_{it}) = -(\bP_{-it}-\bP_{it})\bJ$. Hence the right-hand side of \eqref{fluxt} is
$-\mathfrak b((\bP_{-it}-\bP_{it})\bJ\vC w,w)$. We use the symmetry of $\mathfrak b$ and Lemma \ref{lem-diff}, with $u = w$ and $u' = \bJ\vC w$. Since
$\|(1+H_0)^{1/2}\omega_\infty\bJ\vC w\| = \|(1+H_0)^{1/2}\omega_\infty w\|$, we obtain
\[
2|\Im z|\,|\w{\bJ w,w}|\le Ce^{-2\delta_3(c_\Psi\gamma-R''\gamma^\sigma)}\|(1+H_0)^{\f12}\omega_\infty w\|^2 .
\]
By Proposition \ref{prop-hf} and Lemma \ref{lem-levelk},
\[
|\Im z|\le Ct_0^{-4}e^{-2\delta_3c_\Psi\gamma+2\delta_3R''\gamma^\sigma}\le C_ce^{-c\gamma},
\]
 since $\sigma<1$. This proves \eqref{expsmall}, since
the poles $z_l(\gamma)$ are independent of $\theta$. \ef

\subsubsection{Real poles below the second threshold}

It is natural to ask whether the poles produced by a negative eigenvalue of $-\Delta^{(k)}_V$ are real and,
if so, whether they are embedded eigenvalues of $P(\gamma)$. Below the second threshold $\gamma$, a
real pole is always an embedded eigenvalue.
In order to use the results of \cite{w2}, we assume in addition that $V$ is long-range and satisfies:
\be \label{ass1c}
|V(x)| + \w x |\na V(x)| \le C\w x^{-\ep_1} \qquad \mbox{ for some } \ep_1> 0.
\ee

\begin{prop}\label{prop-flux} Assume \eqref{ass1b} and \eqref{ass1c}, and let $\gamma\ge1$. If
$\lambda \in\,]0,\gamma[$ is a real pole of $R(z,\gamma,\theta)$, then $\lambda$ is an \emph{embedded
eigenvalue} of $P(\gamma)$.
\end{prop}
\pf Let $\rho_1=1+\ep_1$ with $\ep_1$ given by (\ref{ass1c}).  Let $W = -\na V\cdot\na_v$ and $R_0(z,\gamma) = (P_0(\gamma)-z)^{-1}$. By the limiting
absorption principle of \cite{w2}, the boundary values $R_\pm = R_0(\lambda\pm i0,\gamma)$ exist in
$\vL(s,-s)$, $s>\f12$. By \cite[Proposition 3.6]{w4}, there is
$u \in \vH^{-s}\setminus\{0\}$, $\f12<s<\f{\rho_1}2$, such that
\be \label{outg}
u = -R_+Wu .
\ee
Put $f = Wu$. By \eqref{ass1c}, $\w x^{\sigma}f \in \vH^{-1,0}$ with $\sigma = \rho_1-s > s$. Moreover,
$u$, $R_\pm f \in \vH^{1,-s}\subset\vH^{1,-\sigma}$ \cite[Proposition 2.1]{w2}. Hence all the
pairings below are well defined.

By \eqref{JPJ}, $R_0(z,\gamma)^* = \bJ R_0(\bar z,\gamma)\bJ$,
hence $R_- = \bJ R_+^*\bJ$. Using \eqref{outg},
\[
\w{\bJ f, R_+f} = -\w{\bJ f,u}, \qquad
\w{\bJ f, R_-f} = \w{R_+f,\bJ f} = -\overline{\w{\bJ f,u}} ,
\]
so that $\w{\bJ f, (R_+-R_-)f} = -2i\,\Im\w{\bJ Wu,u}$. Since $W$ is skew-adjoint and odd in $v$,
$\bJ W\bJ = -W = W^*$, and
\[
\w{\bJ Wu,u} = \w{\bJ u,Wu} = \overline{\w{\bJ Wu,u}}
\]
 is real.
Therefore
\be \label{flux}
\w{\bJ f, (R_+-R_-)f} = 0 .
\ee

We now compute the jump $R_+-R_-$. At $\theta = 0$, $\zeta = \gamma^{-1}\xi$ and, by Lemma \ref{lem2.3} with
$\ell=1$ and $b\in\,]1,2[$,
\[
\widehat R_0(z,\xi,\gamma) = \chi_0(\gamma^{-1}\xi)\sum_{j\le1}
(E_j(\xi,\gamma)-z)^{-1}\Pi^\zeta_j + r_1 .
\]
Here $E_j(\xi,\gamma) = \gamma j+\gamma^{-1}|\xi|^2$, and $r_1$ is
holomorphic in $\{|z|<b\gamma\}\ni\lambda$. We may assume $\chi_0 = 1$ on $\{|\xi|\le1\}$. Since
$0<\lambda<\gamma$, only $E_0$ takes the value $\lambda$, on the sphere
$S_0 = \{|\xi|^2 = \gamma\lambda\}$. On $S_0$, $\na_\xi E_0 = 2\xi/\gamma \neq 0$ and
$\chi_0(\gamma^{-1}\xi) = 1$, since $|\xi|^2 = \gamma\lambda<\gamma^2$. Hence
\[
\w{\bJ f, (R_+-R_-)f} = 2\pi i \int_{S_0}\w{\bJ\widehat f(\cdot,\xi),\,\Pi^\zeta_0\widehat f(\cdot,\xi)}_v
\,\f{\gamma\,d\sigma_0(\xi)}{2|\xi|} .
\]
The traces on $S_0$ are defined since $\widehat f \in H^{\sigma}(\bR^n_\xi)$ with values in
$\vH^{-1}_v$ and $\sigma>\f12$. By \eqref{Riesz},
$\Pi^\zeta_0 = e_0\w{e'_0,\cdot}_v$ with $e_0 = \psi_0(\cdot+2i\zeta)$ and
$e'_0 = \psi_0(\cdot-2i\zeta)$. As $\psi_0$ is even and has real Taylor coefficients,
$\bJ e_0 = e'_0$, hence
\[
\w{\bJ g,\Pi^\zeta_0 g}_v = \w{g,e'_0}_v\w{e'_0,g}_v = |\w{e'_0,g}_v|^2 \ge 0.
\]

By \eqref{flux} and the jump formula, $\w{e'_0,\widehat f(\cdot,\xi)}_v = 0$, that is,
$\Pi^\zeta_0\widehat f(\cdot,\xi) = 0$, for $d\sigma_0$-almost every $\xi\in S_0$. The jump
$(R_+-R_-)f$ is carried by this term only. Hence it vanishes, and $u = -R_0(\lambda\pm i0,\gamma)Wu$
for \emph{both} signs. By \cite[Lemma 3.2]{w2}, $u \in \vH$ and $\lambda \in \sigma_p(P(\gamma))$.
\ef

\begin{rmk}\label{rmk-width}
The physics intuition says that in the overdamped regime $\gamma\to \infty$, the diffusion dominates the drift and there are no bound states (see \cite{kr,brink}, and \cite[Chapter 10]{risk} for the high-friction limit). Our analysis partially justifies this folklore: outside $o(\f 1 \gamma)$-neighborhoods of the thresholds $\gamma k$, there are no eigenvalues of $P(\gamma)$, except possibly the poles $z_l(\gamma)$ induced by a negative eigenvalue $\eta_0$ of $-\Delta^{(k)}_V$,  for which we do not know whether they are eigenvalues or not.
 By Lemma \ref{lem-levelk}, for a pole $z = z_l(\gamma)$ and a normalized resonant state $w$, the left-hand side of \eqref{fluxt} is $2i\,(-1)^k\,\Im z\,\big(1+O(\gamma^{2\sigma-2})\big)$, $t$ being fixed.
 A rigorous evaluation of the right-hand side of \eqref{fluxt} would decide this
question, which we leave open; see Proposition \ref{prop-flux} for real poles below the second threshold.
\end{rmk}

\appendix
\section{Proof of Theorem \ref{th3.2}}\label{app-Gr}

We use the notation of Subsection \ref{ssec-inv}: $\gamma$ and $\theta$ are fixed, $R(z) = R(z,\gamma,\theta)$, $P = P(\gamma,\theta)$, and $\vP(z)$ is the Grushin problem \eqref{Gr1}.

\begin{prop}\label{prop-Gr}  Assume \eqref{ass1}. Let $\gamma\ge\gamma_0$ and $z\in\vU_{k+1}(\gamma)$. Then,

\textup{(i)}  $\vP(z) : D(P)\times\vG \to
L^2(\bR^{2n})\times\vG$ is closed and injective, and its range is dense but is not the whole space.
Consequently, $\vP(z)^{-1}$, with domain $\Ran\,\vP(z)$, is a closed, densely defined
and unbounded operator on $L^2(\bR^{2n})\times\vG$, with range $D(P)\times\vG$.

\textup{(ii)} $(f,g)\in L^2(\bR^{2n})\times\vG$ belongs to $\Ran\,\vP(z)$ if
and only if there exists $w\in L^2(\bR^{2n})\cap\Ran\Pi'$ such that
\be \label{RanPint}
(P'-z)w = \Pi'(f-L_\theta E_+g) \quad \mbox{ and } \quad
E_-\big(f-(P-z)(E_+g+w)\big)\in\vG.
\ee
 Such a $w$ is unique, and
\be \label{Pinvint0}
\vP(z)^{-1}(f,g) = \Big(E_+g+w,\ E_-\big(f-(P-z)(E_+g+w)\big)\Big) .
\ee
\end{prop}
\pf (i) Since $P$ is closed and $E_\pm$ are bounded, $\vP(z)$ is closed.
If $\vP(z)(u,u_-) = 0$, then $E_-u=0$, that is,
$u = \Pi'u$. Applying $\Pi'$ to $(P-z)u = -E_+u_- \in \Ran\Pi$ gives
$(P'-z)u = 0$, whence $u=0$ by Lemma \ref{lem3.1}. Then $E_+u_-=0$, that is,
$u_-=0$. This proves the injectivity.

{\it The range is dense.} Let $K$ be the complex conjugation and ${}^tA := KA^*K$ the transpose of a closed
operator $A$. Since $P$ is the closure of its restriction to $C_0^\infty(\bR^{2n})$ and has
maximal domain, ${}^tP$ is the maximal realization of the formal transpose, that is,
${}^tP = \bJ P\bJ$ with $\bJ$ as in \eqref{bJ}, since $H_0$ is even and $L_\theta$ is odd in $v$.
Since the Hermite functions are real, ${}^tE_\pm = E_\mp$. Moreover, $\bJ E_+ = E_+\bJ_\vG$ and $E_-\bJ = \bJ_\vG E_-$, where
$\bJ_\vG = (-1)^j$ on $\vG_j$. Hence
\[
{}^t\vP(z) = \left(\begin{array}{cc} {}^tP-z & E_+ \\ E_- & 0\end{array}\right)
= \left(\begin{array}{cc} \bJ & 0 \\ 0 & \bJ_\vG\end{array}\right)\vP(z)
\left(\begin{array}{cc} \bJ & 0 \\ 0 & \bJ_\vG\end{array}\right)
\]
is injective, and so is $\vP(z)^* = K\,{}^t\vP(z)K$. Thus
\[
(\Ran\vP(z))^\perp = \ker\vP(z)^* = \{0\}.
\]

{\it The range is not the whole space}. Let $u\in D(P)$ and $u_j := E^{(j)}_-u$. Since $L_\theta$ maps
$\Ran\Pi_j$ into $\Ran\Pi_{j-1}\oplus\Ran\Pi_{j+1}$ and $H_0 = 0$ on $\Ran\Pi_0$, one has, in the sense of
distributions,
$E^{(0)}_-(P-z)u = -zu_0-{}^t\bX^\theta_0u_1$ by \eqref{Xblocks}, so that ${}^t\bX^\theta_0u_1\in L^2(\bR^n)$.
Take $g\in\vG$ with only a component $g_1 = (h,0,\dots,0)$ on the level $1$, where $h\in L^2(\bR^n)$ and
$\p_{x_1}h\notin L^2(\bR^n)$. As the coefficients $(\p_{x_l}V)(e^\theta x)$ are bounded by \eqref{ass1},
${}^t\bX^\theta_0g_1 = {}^tX^\theta_1h\notin L^2(\bR^n)$, up to a non-zero constant. If
$\vP(z)(u,u_-) = (0,g)$, then $u_1 = g_1$ and $(P-z)u = -E_+u_-\in L^2$, a
contradiction. Hence $(0,g)\notin\Ran\vP(z)$.

Finally, $\vP(z)^{-1}$ with domain $\Ran(\vP(z))$ is closed, densely defined and  unbounded.

(ii) If $\vP(u,u_-) = (f,g)$, then $w := \Pi'u\in L^2\cap\Ran\Pi'$ satisfies \eqref{RanPint}. Indeed,
\[
\Pi u = E_+g, \quad \Pi'(P-z)u = \Pi'f, \quad \Pi'(P-z)E_+g = \Pi'L_\theta E_+g
\]
and
$E_-(f-(P-z)u) = u_-$. Conversely, if $w$ satisfies \eqref{RanPint}, set $u := E_+g+w\in L^2$ and
$u_- := E_-(f-(P-z)u)\in\vG$. Then
\[
\Pi'(f-(P-z)u) = \Pi'f-\Pi'L_\theta E_+g-(P'-z)w = 0, \]
so that
$f-(P-z)u = E_+u_-\in L^2$. Hence $u\in D(P)$, since $P$ has maximal domain. Moreover,
$(P-z)u+E_+u_- = f$ and $E_-u = g$. The uniqueness of $w$ follows from the injectivity of $\vP$.
\ef

If $\vP(z)(u,u_-) = (f,0)$, then
$E_-u = 0$, and applying $\Pi'$ to the first equation gives $u = R'f$. Hence \eqref{Dm} holds.
Under \eqref{ass1b}, $R'$ extends to $\vS'$, so that the only possible $w$ in
\eqref{RanPint} is
\[
w = R'\Pi'(f-L_\theta E_+g) = \Pi'(\vE f+\vE_+g).
\]
Moreover, by \eqref{r1}, \eqref{r2} and $E_-E_+ = 1$,
$E_-(f-(P-z)(E_+g+w)) = \vE_-f+\vE_{-+}g$. Thus Proposition \ref{prop-Gr}(ii) gives \eqref{RanP}.
\medskip

\noindent{\bf Proof of Theorem \ref{th3.2}.} (i) Let $w\in\vG$ and $u = -RE_+w\in D(P)$. Then $\vP(u,w) = (0,E_-u)$, so that
$E_-u\in\vD_+$, $\vE_+E_-u = u$ and $\vE_{-+}E_-u = w$. Hence $-E_-RE_+$ is a right inverse of $\vE_{-+}$.
Conversely, if $w\in\vD_+$, then $(P-z)\vE_+w = -E_+\vE_{-+}w$, so that
$\vE_+w = -RE_+\vE_{-+}w$ and $w = E_-\vE_+w = -E_-RE_+\vE_{-+}w$. Thus $-E_-RE_+$ is also a left inverse.
Hence $\vE_{-+}$ is invertible with bounded inverse, so that it is closed, and $\vD_+ = \Ran(E_-RE_+)$.

By the proof of Proposition \ref{prop-Gr}(i),
\[
R^* = K\,{}^tR\,K = K\bJ R\bJ K, \qquad (E_-RE_+)^* = K\bJ_\vG(E_-RE_+)\bJ_\vG K ,
\]
so that $(E_-RE_+)^*$ is injective. Hence $\vD_+$ is dense. Let
$f\in\vD_-$ and $w = \vE_{-+}^{-1}\vE_-f\in\vD_+$. Then $\vP(\vE f,\vE_-f) = (f,0)$ and
$\vP(\vE_+w,\vE_-f) = (0,w)$, so that $\vP(\vE f-\vE_+w,0) = (f,-w)$. Thus
$(P-z)(\vE f-\vE_+w) = f$, that is, $Rf = \vE f-\vE_+w$, which is \eqref{Grushindilated2}. Moreover,
$E_-(\vE f-\vE_+w) = -w$, that is, $E_-Rf = -w$, which is \eqref{EEmp}.

(ii) If $u\in\ker(P-z)$, then $\vP(u,0) = (0,E_-u)$, so that $E_-u\in\vD_+$,
$\vE_+E_-u = u$ and $\vE_{-+}E_-u = 0$. Conversely, if $w\in\ker\vE_{-+}$, then
$(P-z)\vE_+w = 0$ and $E_-\vE_+w = w$. The last assertion follows, since the spectrum of
$P$ outside $\sigma_{\rm ess}(P) = \sigma(P_0)$ consists of isolated
eigenvalues.
\ef

\section{Large-friction expansion of the free resolvent}\label{app-free}

In this Appendix, we give some results on the reduced  resolvent $R'(z, \gamma, \theta)$, used for the large-friction expansions of poles.
We use the notation of Section \ref{sect2}. Thus $P_0(\gamma,\theta) = \gamma H_0+L_0$ with
$L_0 := e^{-\theta}v\cdot\na_x$, where $\gamma\ge1$ is large and $\theta\in\bC$, with $\varphi\in\,]0,\f\pi4[$, is fixed. In the fibre
$\xi$,
\[
\widehat P_0(\xi) = \widehat P_0(\xi,\gamma,\theta) = \gamma H_0+ie^{-\theta}v\cdot\xi
\]
 has the eigenvalues
$E^\theta_j(\xi,\gamma)$, $j\in\bN$ (the channels), and $\widehat R_0(\xi) = \widehat R_0(z,\xi,\gamma,\theta)$. The operators $a_l$, $a^*_l$
are those of \eqref{creation}, so that $H_0 = \sum_la^*_la_l$ by \eqref{aa}. Finally, $\Pi_j$ is the orthogonal
projection onto $\ker(H_0-j)$ and $\mu_j = \mbox{rank}\,\Pi_j$. The condition $\varphi>0$ is not used in Proposition \ref{prop-Rprime}, which
holds for $|\theta|<a$ and is used for real $\theta$, including $\theta = 0$, in Proposition \ref{prop-det}(iii).

If the lowest levels of $H_0$ are removed, the free resolvent has an expansion in powers of $\gamma^{-1}$. Let $k\in\bN$ and $m\ge k+1$
be integers, $\Pi'_m := \sum_{j\ge m}\Pi_j$, $z = \gamma k+z'$ with $|z'|\le c_1$, and
\be \label{Rprime}
R'_0(z,\gamma,\theta):=\big(\Pi'_mP_0(\gamma,\theta)\Pi'_m-z\big)^{-1}\Pi'_m, \qquad
L' := \Pi'_mL_0\Pi'_m, \qquad S_{k,m} := (H_0-k)^{-1}\Pi'_m .
\ee
For
$m=k+2$, $R'_0$ is the free part of the operator $R'$ of Lemma \ref{lem3.1}. As in \eqref{e4.19} and
\eqref{e4.19b}, since $H_0\ge m$ on $\Ran\Pi'_m$,
\be \label{nr}
\big\|R'_0(z,\gamma,\theta)\big\| + \big\|(1+H_0)^{\f12}R'_0(z,\gamma,\theta)\big\| \le \f C\gamma .
\ee
Set, for $N \ge 1$,
\be \label{wN}
w_N:=\sum_{p=0}^{N-1}(-1)^p z'^{\,N-1-p}
\sum_{\substack{i_0,\dots,i_p\ge0\\ i_0+\cdots+i_p=N-1-p}}
S_{k,m}^{i_0+1}L'S_{k,m}^{i_1+1}L'\cdots L'S_{k,m}^{i_p+1} .
\ee
It is a polynomial in $z'$ and in $L' = e^{-\theta}\,\Pi'_mv\cdot\na_x\Pi'_m$, hence entire in $\theta$.
The first ones are, with $S = S_{k,m}$,
\be \label{w123}
w_1=S,\qquad w_2=z'S^2-SL'S,\qquad w_3=z'^2S^3-z'\big(S^2L'S+SL'S^2\big)+SL'SL'S .
\ee

\begin{prop}\label{prop-Rprime} Let $f\in\vS(\bR^{2n})$, $s\in\bR$ and $N\ge1$.

\textup{(i)} One has
\be \label{remprime}
\Big\| (H_0+1)^{\f12}\w{D_x}^s\Big(R'_0(z,\gamma,\theta)f-\sum_{q=1}^{N}\gamma^{-q}w_qf\Big)\Big\|
\ \le\ C_N\,\gamma^{-(N+1)}\,\big\|\w{D_x}^{s+N}f\big\| .
\ee

\textup{(ii)} Let $m = k+2$ and let $R'(z,\gamma,\theta)$ be the operator of Lemma \ref{lem3.1}. There
are operators $\widetilde w_q$, polynomial in $z'$ and holomorphic in $\theta\in\{|\theta|<a\}$, such
that
\be \label{remprime2}
\Big\| \w{D_x}^s\Big(R'(z,\gamma,\theta)f-\sum_{q=1}^{N}\gamma^{-q}\widetilde w_qf\Big)\Big\|
\ \le\ C_N\,\gamma^{-(N+1)}\,\big\|(H_0+1)^{N}\w{D_x}^{s+N}f\big\| .
\ee
The operators $\widetilde w_q$ are finite sums of products of the $w_j$ and of
$W'_\theta = \Pi'_mW_\theta\Pi'_m$. The constants are locally uniform in $\theta$.
\end{prop}
\pf (i) $P_0(\gamma,\theta)$ commutes with $D_x$, and in the fibre $\xi$,
\[
\Pi'_m\widehat P_0(\xi,\gamma,\theta)\Pi'_m-z = (\gamma(H_0-k)-z')\big(1+\rho_0L'(\xi)\big)
\]
with

\[
\rho_0 := (\gamma(H_0-k)-z')^{-1}\Pi'_m = \sum_{i\ge0}z'^i\gamma^{-i-1}S_{k,m}^{i+1}.
\]

One has $\|(H_0+1)^{\f12}\rho_0\|\le C\gamma^{-1}$ and $\|\rho_0L'(\xi)\|\le C\w\xi\gamma^{-1}$, since $\rho_0$
gains one power of $H_0$ while $L'(\xi)$ costs $H_0^{1/2}|\xi|$. For $|\xi|\le\rho\gamma$, with $\rho$ small,
$\big(1+\rho_0L'\big)^{-1}\rho_0 = \sum_r(-\rho_0L')^r\rho_0$. Collecting the powers of $\gamma^{-1}$
gives \eqref{wN}. The omitted terms all begin with $\rho_0$ or $S_{k,m}$ on the left. After multiplication by $(H_0+1)^{\f12}$, they are
bounded by $C_N\gamma^{-(N+1)}\w\xi^{N}$. For
$|\xi|\ge\rho\gamma$, one uses \eqref{nr} and the fact that $w_q$ has degree $\le q-1$ in $\xi$.

(ii) By \eqref{nr} and \eqref{ass1}, $\|W'_\theta R'_0\| \le C\gamma^{-1}$, also after conjugation by
$\w{D_x}^s$, and
\[
R' = \sum_{r=0}^N R'_0(-W'_\theta R'_0)^r + R'(-W'_\theta R'_0)^{N+1}.
\]
The last term
is $O(\gamma^{-N-2})$ on every $H^s$ by Lemma \ref{lem3.4}. We expand each factor $R'_0$ by (i),
from the right. The operators $W'_\theta$ and $w_q$ map $\vS$ into $\vS$, and each $W'_\theta$ costs
one power of $(H_0+1)^{\f12}$. This gives \eqref{remprime2}. Under \eqref{ass1},
$(\na V)(e^\theta x)$ is holomorphic in $\theta$ for $|\theta|<a$, and so are the $\widetilde w_q$.
\ef

\section*{Acknowledgments}

The author warmly thanks Jean-Marc Bouclet for useful discussions on the relationship between the
KFP operator and the Witten Laplacian. Research supported in part by Nanjing University of
Aeronautics and Astronautics under grant No. 90YAH25017.

\end{document}